\documentclass[12pt]{amsart}

\usepackage[english]{babel}

\usepackage[
letterpaper,
top=2cm, bottom=2cm,
left=3cm, right=3cm,
marginparwidth=1.75cm
]{geometry}
\usepackage{amsmath, amsfonts, amssymb, amsthm}

\usepackage{graphicx}
\usepackage{enumitem}
\usepackage{aliascnt}
\usepackage{tikz}
\usetikzlibrary{calc, arrows, decorations.text}

\usepackage[
colorlinks = true,
pagebackref = true,
bookmarks = true,
unicode,
pdfauthor = author
]{hyperref}
\usepackage{xcolor}

\newtheorem*{rmks}{Remarks}

\newtheorem{thm}{Theorem}[section]

\newaliascnt{lma}{thm}
\newtheorem{lma}[lma]{Lemma}
\aliascntresetthe{lma}

\newaliascnt{prop}{thm}
\newtheorem{prop}[prop]{Proposition}
\aliascntresetthe{prop}

\newaliascnt{cor}{thm}
\newtheorem{cor}[cor]{Corollary}
\aliascntresetthe{cor}

\newaliascnt{conj}{thm}
\newtheorem{conj}[conj]{Conjecture}
\aliascntresetthe{conj}

\newaliascnt{rmk}{thm}
\newtheorem{rmk}[rmk]{Remark}
\aliascntresetthe{rmk}

\newaliascnt{defn}{thm}

\aliascntresetthe{defn}

\newcommand{\R}{\mathbb{R}}

\numberwithin{equation}{section}
\tikzset{ > =stealth'}

\begin{document}

\title[Uniqueness of critical points on triangles]{Uniqueness of critical points of the second Neumann eigenfunctions on triangles}

\author[H. Chen]{Hongbin Chen}
\address{School of Mathematics and Statistics, Xi'an Jiaotong University, Xi'an, 710049, China}
\email{hbchen@xjtu.edu.cn}

\author[C. Gui]{Changfeng Gui}
\address{
Department of Mathematics, University of Macau, Taipa, Macau and
Zhuhai UM Science and Technology Research Institute,
Hengqin, Guangdong, 519031, China}
\email{changfenggui@um.edu.mo}

\author[R. Yao]{Ruofei Yao}
\address{
School of Mathematics, South China University of Technology, Guangzhou, 510641, China}
\email{yaorf5812@126.com; yaoruofei@scut.edu.cn}


\begin{abstract}
This paper investigates the second Neumann eigenfunction ${u}$ of a planar triangle ${T}$. In a recent paper by Judge and Mondal [Ann. Math., 2022], it was shown that ${u}$ has no critical points in the interior of ${T}$. In this paper, we show that ${u}$ has at most one non-vertex critical point and that ${u}$ is monotone in a certain direction in ${T}$. More precisely, when ${T}$ is not equilateral, we show that ${u}$ vanishes at some vertex if and only if ${T}$ is superequilateral, and that ${u}$ has a non-vertex critical point if and only if ${T}$ is acute and not superequilateral. These results confirm both the original theorem and Conjecture 13.6 of Judge and Mondal [Ann. Math., 2020]. We also resolve the objective of Polymath 7 (research thread 1), namely, that the extrema of ${u}$ are attained only at the endpoints of the longest side. In addition, we settle a conjecture of Siudeja [Proc. Amer. Math. Soc., 2016] on the ordering of mixed Dirichlet--Neumann Laplacian eigenvalues for triangles. Our proofs combine the continuity method, eigenvalue inequalities, the maximum principle, and the moving plane method. 
\end{abstract}

\keywords{Hot spots; Neumann eigenfunctions; Continuity method; Monotonicity; Eigenvalue inequality; Moving plane method}
\subjclass[2020]{Primary 35P05; Secondary 58J50, 35B50, 35J05, 35J25}

\maketitle


\section{Introduction} \label{Sect1intro}

Suppose that $\Omega$ is a bounded domain in Euclidean space $\R^{{n}}$ with Lipschitz boundary. We consider the Neumann eigenvalue problem
\begin{equation} 
\begin{cases}
\Delta {u} + \mu {u} = 0 & \text{ in } \Omega, 
\\ \nabla {u} \cdot \nu = 0 & \text{ on } \partial\Omega, 
\end{cases}\end{equation} 
where $\nu$ denotes the unit outward normal vector defined at $C^1$ points of $\partial\Omega$. It is well known that for nice enough bounded domains $\Omega$ (e.g., convex domains or domains with smooth boundary), there exists a nondecreasing sequence of eigenvalues satisfying $0 = \mu_{1} < \mu_{2} \leq \mu_{3} \leq \cdots \to \infty$. 

The hot spots conjecture, first proposed by Rauch at a conference in 1974 \cite{Rau74}, asserts that the second Neumann eigenfunction attains its maximum and minimum only on the boundary of the domain. The second Neumann eigenfunction of the Laplacian plays a fundamental role in the study of partial differential equations and boundary value problems; accordingly, the hot spots conjecture has attracted considerable attention. 

For a historical overview and a discussion of different versions of the hot spots conjecture, the reader may refer to \cite{BB99}. The first published positive result was obtained by Kawohl for cylindrical domains (see \cite[Corollary 2.15]{Kaw85}); he also stated that the conjecture holds for specific domains such as parallelepipeds, balls, and annuli (see \cite[p.~46]{Kaw85}). Additionally, in the same manuscript, Kawohl \cite{Kaw85} suggested that the conjecture might be false in general, as the second Neumann eigenfunction could attain its maximum in the interior of certain nonconvex domains. Counterexamples in specific multiply connected settings were constructed by Burdzy and Werner \cite{BW99}, Bass and Burdzy \cite{BB00}, Burdzy \cite{Bur05}, and Kleefeld \cite{Kle21}. More recently, de Dios \cite{deD24} disproved the hot spots conjecture for convex sets in sufficiently high dimensions. 

Motivated by these counterexamples, the conjecture is now often stated as follows.

\begin{conj}
The second Neumann eigenfunction attains its global maximum and minimum only on the boundary for convex domains in $\R^{2}$ and, more generally, for simply connected planar domains.
\end{conj}

The hot spots conjecture remains open, despite many special cases having been solved. In 1999, under certain technical assumptions, Ba\~nuelos and Burdzy \cite{BB99} were able to verify the conjecture for certain domains with a line of symmetry and for obtuse triangles. A year later, Jerison and Nadirashvili \cite{JN00} proved the conjecture for domains with two axes of symmetry. Pascu \cite{Pas02} proved the conjecture for convex domains with a single axis of symmetry, under the additional assumption that the second Neumann eigenfunction is antisymmetric. 
In a different direction, building on the work in \cite{BB99}, Atar and Burdzy \cite{AB04} proved the conjecture for lip domains, i.e., planar domains bounded by two Lipschitz graphs with Lipschitz constant $1$.
Miyamoto \cite{Miy09} proved the conjecture for convex domains close to the ball and, in \cite{Miy13}, constructed a polygon showing that there can be arbitrarily many isolated hot spots on the boundary. 
Siudeja \cite{Siu15} and Judge and Mondal \cite{JM20, JM22ar} proved the conjecture for acute triangles, while Judge and Mondal \cite{JM22c, JM25} also studied the critical points of eigenfunctions for polygons. 
The works \cite{AB02, NSY20} studied the location of the nodal line of the second Neumann eigenfunction for certain triangles. Steinerberger \cite{Ste20} showed that the global extrema must be close, in a specified sense, to the endpoints of a diameter of a convex domain. 
The papers \cite{Ste23, MPW23} showed that the global extrema are controlled by the boundary extrema. 
Recently, Rohleder \cite{Roh21} announced a non-probabilistic proof of the main result of \cite{AB04}, while the first and third authors of this paper, together with Wu, studied in \cite{CWY25} the monotonicity properties of the second even Neumann eigenfunction on a class of symmetric domains. For further developments of the conjecture in higher dimensions, we refer the reader to \cite{Yan11, CLW19, KR24}.

In another direction, Ba\~nuelos, Pang, and Pascu \cite{BPP04} proposed the ``hot spots'' property for the mixed (Dirichlet--Neumann) eigenvalue problem: namely, whether the positive eigenfunction of the mixed problem attains its maximum only on the Neumann boundary. This question is closely related to the hot spots conjecture \cite{BP04, BPP04, JN00, Pas02} and may extend to certain nonlinear contexts \cite{BP89, CLY21, CY18, MY26, YCG21, YCL18, LY24}. In particular, \cite{LY24} investigates monotonicity properties of mixed boundary problems in triangular domains, which are then applied to establish eigenvalue inequalities via the maximum principle.

The acute triangle is the first important example without symmetry or special shape assumptions. 
In 2012, the hot spots conjecture for acute triangles became a ``Polymath Project" \cite{Pol12}. 
The conjecture for obtuse or right triangles was resolved in \cite{BB99, AB04} and, more recently, in \cite{Roh21}. The conjecture for isosceles triangles follows by combining \cite{LS10} and \cite{Miy13}. 
In 2015, Siudeja \cite{Siu15} proved the conjecture for acute triangles with smallest angle at most $\pi/6$ by sharpening ideas developed by Miyamoto \cite{Miy09, Miy13}. 
For non-equilateral triangles, the second Neumann eigenvalue has been shown to be simple for obtuse and right triangles in \cite{AB04} and for acute triangles in \cite{Siu15}. 
More recently, Judge and Mondal \cite{JM20, JM22ar} made a significant contribution by showing that the second Neumann eigenfunction for a triangle has no interior critical points; consequently, the global extrema are achieved exclusively on the boundary.


\subsection{The second Neumann eigenfunction in triangles}

In this paper we focus on some intrinsic features of the second Neumann eigenfunction ${u}$ in a triangle ${T}$. 
The symbol ${T}$ usually denotes an open triangle. We call ${T}$ a \textbf{subequilateral} triangle if it is isosceles with apex angle less than $\pi/3$, and ${T}$ a \textbf{superequilateral} triangle if it is isosceles with apex angle greater than $\pi/3$. Since the gradient $\nabla {u}$ of ${u}$ always vanishes at each vertex of the triangle ${T}$, the vertices are the trivial (i.e., obvious) critical points. Throughout this paper, we focus on the non-vertex (i.e., nontrivial) critical points, and set
\begin{equation*}
\operatorname{crit}_{\mathrm{nv}}({u}) = \{{p} \in \overline{T}: \, |\nabla {u} ({p})| = 0, \; \text{${p}$ is not a vertex of ${T}$} \}. 
\end{equation*}

Our first result concerning the hot spots conjecture for triangles is as follows. 

\begin{thm} \label{thm12}
Let ${u}$ be a second Neumann eigenfunction in a triangle ${T}$ that is not equilateral. 
Then ${u}$ has at most one non-vertex critical point. 
More precisely, 
\begin{enumerate}[label = \rm(\arabic*), start = 1]
\item \label{CGY0102ita}
The non-vertex critical point of ${u}$ (if it exists) is unique, is a saddle point, and lies in the interior of the shortest side of ${T}$. 
\item \label{CGY0102itb}
A non-vertex critical point of ${u}$ exists if and only if ${T}$ is an acute triangle that is not superequilateral.
\item \label{CGY0102itc}
${u}$ vanishes at a vertex if and only if ${T}$ is superequilateral (with ${u}$ antisymmetric with respect to the symmetry axis of ${T}$); otherwise (i.e., when ${T}$ is not superequilateral), the nodal line of ${u}$ connects interior points of the two longest sides.
\item \label{CGY0102itd}
The global extrema of ${u}$ are attained at and only at the endpoints of the longest side.
\item \label{CGY0102ite}
The eigenfunction ${u}$ is strictly monotone in the direction orthogonal to the shortest side of the triangle.
\item \label{CGY0102itf}
As a byproduct, the second Neumann eigenfunction is unique (up to multiplication by a constant).
\end{enumerate}
\end{thm}

\begin{figure}[htp] \centering\vspace*{-2ex}
\begin{tikzpicture}[scale = 5]
\pgfmathsetmacro\AngleL{75}    
\pgfmathsetmacro\AngleM{60}    
\pgfmathsetmacro\AngleS{180 - \AngleL - \AngleM}  
\pgfmathsetmacro\Short{1.00}   
\pgfmathsetmacro\Mid{\Short*sin(\AngleM)/sin(\AngleS)}    
\pgfmathsetmacro\Long{\Short*sin(\AngleL)/sin(\AngleS)}   
\pgfmathsetmacro\xAA{0}        
\pgfmathsetmacro\yAA{0}
\pgfmathsetmacro\xBB{\Short}   
\pgfmathsetmacro\yBB{0}
\pgfmathsetmacro\xCC{\Mid*cos(\AngleL)}   
\pgfmathsetmacro\yCC{\Mid*sin(\AngleL)}
\draw [black]
(\xCC, \yCC) node[left] {${z}_{3}$} --
(0, 0)      node[above left] {${z}_{1}$} --
(\Short, 0) node[above right] {${z}_{2}$} -- cycle; 
\node at ({(\xAA + \xBB)/2}, {(\yAA + \yBB)/2}) [below] {\scriptsize Shortest Side}; 

\node at ({(\xBB*0.4 + \xCC*0.6) + 0.055*\Short*sin(\AngleM)}, {(\yBB*0.4 + \yCC*0.6) + 0.061*\Long*cos(\AngleM)}) [rotate = -\AngleM] {\scriptsize Longest Side};  
\node at ({(\xCC*0.6 + \xAA*0.4) - 0.055*\Short*sin(\AngleL)}, {(\yCC*0.6 + \yAA*0.4) + 0.061*\Mid*cos(\AngleL)}) [rotate = \AngleL] {\scriptsize Medium Side}; 
\pgfmathsetmacro\xPP{\xCC*0.5 + \xAA*0.5}
\pgfmathsetmacro\yPP{\yAA}
\draw [blue] (\xPP, \yPP) node[above] {\scriptsize ${P}$}; 
\fill [blue] (\xPP, \yPP) circle (0.5pt) node[below = 9pt] {Critical Point};  
\draw [very thin, dashed] (\xCC, \yCC) -- (\xCC, \yAA); 
\fill (\xCC, 0) circle (0.3pt) node[above right] {\scriptsize ${H}$};  
\draw [thin] ({\xCC}, {\yAA + \Short*0.04}) -- ({\xCC - \Short*0.04}, {\yAA + \Short*0.04}) -- ({\xCC - \Short*0.04}, {\yAA}); 
\pgfmathsetmacro\xRRa{\xCC - 0.68*\Mid*cos(\AngleL)}
\pgfmathsetmacro\yRRa{\yCC - 0.68*\Mid*sin(\AngleL)}
\pgfmathsetmacro\xSSa{\xCC + 0.55*\Long*cos(\AngleM)}
\pgfmathsetmacro\ySSa{\yCC - 0.55*\Long*sin(\AngleM)}
\pgfmathsetmacro\xRRb{\xRRa + 0.2*\Short*sin(\AngleL)}
\pgfmathsetmacro\yRRb{\yRRa - 0.2*\Short*cos(\AngleL)}
\pgfmathsetmacro\xSSb{\xSSa - 0.2*\Short*sin(\AngleM)}
\pgfmathsetmacro\ySSb{\ySSa - 0.2*\Short*cos(\AngleM)}
\draw [red, thick]
(\xRRa, \yRRa) .. controls (\xRRb, \yRRb) and (\xSSb, \ySSb) .. (\xSSa, \ySSa);  
\draw [red, thin]
({(\xRRb + \xSSb)/2}, {(\yRRb + \ySSb)/2}) node[below] {\scriptsize $\mathcal{Z}({u}) = \{{u}=0\}$}; 
\fill (\xCC, \yCC) circle (0.4pt) node[right] {\small\bf Max};  
\fill (\xBB, \yBB) circle (0.4pt) node[below = 4pt, right] {\small\bf Min};  
\fill (\xAA, \yAA) circle (0.4pt) node[below = 4pt, left] {\small\bf local Min};  
\end{tikzpicture} \vspace*{-2ex}
\caption{The result in Theorem \ref{thm12} when ${T}$ is acute and non-isosceles}
\label{fig1}
\end{figure}

The result is illustrated in \autoref{fig1}. For an equilateral triangle, the second Neumann eigenspace is two-dimensional, and Properties~\ref{CGY0102ita}-\ref{CGY0102ite} continue to hold with minor adjustments reflecting the equal side lengths and the multiplicity of the eigenvalue; see \autoref{rmk65} for details. 

The proof of \autoref{thm12} does not depend on the main result (nonexistence of an interior critical point) in \cite[Theorem A]{JM22ar}. 
Property \ref{CGY0102ita} completes the proof of the original result \cite[Theorem 1.1]{JM20}, for which the erratum \cite{JM22ar} pointed out that the original proof was incomplete in the case of acute triangles. 
Property \ref{CGY0102itb} was originally proposed as a conjecture at the end of the notable work by Judge and Mondal \cite[Conjecture 13.6]{JM20}. 
Property \ref{CGY0102itc} settles \cite[Conjecture 6.60]{LS17}; consequently, it recovers the results of \cite{NSY20}, where only the right triangle with angle $\pi/5$ was treated using a numerical proof based on finite elements. 
We note that Property \ref{CGY0102ita} immediately implies that the global extrema are achieved only at the vertices, as suggested several times before (e.g., \cite[Conjecture 6.52]{LS17}). 
Besides, Property \ref{CGY0102itd} is a stronger result, first observed numerically by Terence Tao at the beginning of \href{https://polymathprojects.org/2012/06/12/polymath7-research-thread-1-the-hot-spots-conjecture/}{\textcolor{black}{Polymath7 research thread 1}} \cite{Pol12}. 
Property \ref{CGY0102ite}, the monotonicity property, is essential in the proof of the theorem. 
Property \ref{CGY0102itf}, the simplicity of the second eigenvalue, was first proven in \cite{Siu15} via the inequality $2\mu_{2} < \mu_{2} + \mu_{3}$ obtained by estimating eigenvalue bounds. Here we give an alternative proof by directly showing the linear dependence of any pair of second eigenfunctions. 

\begin{rmk} \label{rmk13} 
Judge and Mondal \cite{JM22ar} proved that the second Neumann eigenfunction on a triangle has no interior critical points. We briefly compare their approach with ours.

To rule out interior critical points, Judge and Mondal argue by contradiction. Assuming that a second Neumann eigenfunction on some acute triangle ${T}_{0}$ has an interior critical point, they carry out a detailed analysis to produce an acute triangle ${T}_{\xi}$ whose corresponding eigenfunction ${u}_{\xi}$ possesses four non-vertex critical points on the boundary (see \cite[Proposition~3.7]{JM22ar} for their types and locations). These boundary critical points force the directional derivatives to have many nodal domains. By selecting a linear combination of two suitable directional derivatives supported in two disjoint nodal domains as a test function, they obtain a contradiction. 

In the present paper, we leverage a monotonicity property to show that the eigenfunction has at most two non-vertex critical points; see \autoref{lma35}. To establish the uniqueness of the non-vertex critical point, we again argue by contradiction, assuming there are exactly two such points. In this setting, the directional derivative has fewer nodal domains, which is a key difference from the strategy in \cite{JM22ar}. Nevertheless, by choosing an appropriate nodal domain, we show that the second Neumann eigenvalue is strictly larger than the first eigenvalue of a mixed (Dirichlet--Neumann) problem on a subdomain, leading to a contradiction; see \autoref{thm38}. In addition, we obtain two further results. First, we prove a rigidity statement: the nodal line meets a vertex if and only if the triangle is superequilateral. Second, we precisely locate the global extrema, in agreement with numerical observations due to Terence Tao (see \cite{Pol12}). 
\end{rmk}


\subsection{Inequalities for mixed eigenvalues in triangles} 

The proof of \autoref{thm12} relies on comparison inequalities between Neumann eigenvalues and eigenvalues of mixed (Dirichlet--Neumann) type. A typical example of this kind of comparison was established by Siudeja \cite{Siu15} and generalized by Lotoreichik and Rohleder \cite{LR17}. 
For a (Lipschitz) domain $\Omega$ with a portion $\Gamma_{D}$ of $\partial\Omega$, we denote by $\lambda_{1}^{\Gamma_{D}}(\Omega) = \lambda_{1}(\Omega, \Gamma_{D})$ or simply $\lambda_{1}^{\Gamma_{D}}$ the first (smallest) eigenvalue of the mixed (Dirichlet--Neumann) boundary value problem
\begin{equation} \label{CGY0104Mix}
\begin{cases}
\Delta \varphi + \lambda \varphi = 0 & \text{ in } \Omega, \\
\varphi = 0 & \text{ on } \Gamma_{D}, \\ 
\partial_{\nu}\varphi = 0 & \text{ on } \Gamma_{N}: = \partial\Omega \setminus \Gamma_{D}. 
\end{cases}
\end{equation}
An important related question is to optimize the placement of the mixed boundary conditions in a fixed domain. 
For a triangle, we denote by $\lambda_{1}^{{X}}$ the smallest eigenvalue when the Dirichlet condition is imposed on a combination of sides indicated by the symbol $X$, where ${X}$ is any combination of the letters ${L}$, ${M}$, ${S}$, standing for the longest, medium and shortest sides, respectively; Neumann conditions are imposed on any side not listed in ${X}$. Let $\lambda_{1}$ denote the first Dirichlet eigenvalue of the same triangle, i.e., $\lambda_{1} = \lambda_{1}^{LMS}$. 

Our second result concerning eigenvalue inequalities for triangles is as follows. 

\begin{thm} \label{thm14Eig}
For an arbitrary triangle ${T}$, we have
\begin{equation*}
\lambda_{1}^{S} < \lambda_{1}^{M} < \lambda_{1}^{L} < \lambda_{1}^{MS} < \lambda_{1}^{LS} < \lambda_{1}^{LM} < \lambda_{1}
\end{equation*}
as long as all the sides have different lengths. 
\end{thm}

\autoref{thm14Eig} was originally proposed as a conjecture in \cite[Conjecture~1.2]{Siu16} (see also \cite[Conjecture~6.13]{LS17}), where partial results were obtained. Here we provide a complete proof. Our argument relies on the monotonicity of the first mixed eigenfunction and on the monotonicity of the corresponding eigenvalue under deformations of triangles; see \autoref{Sect4a}.


\subsection{Sketch of the proof of \autoref{thm12}}

Over the past four decades, the monotonicity of the second Neumann eigenfunction has been studied in numerous significant works; see \cite{Kaw85, BB99, JN00, Pas02, AB04, BP04, BPP04, Siu15, Roh21, CWY25}. Jerison \cite{JN00, Jer19} conjectured that the monotonicity property holds for centrally symmetric convex domains, whereas eigenfunctions can fail to be monotone in general (e.g., \cite{Miy13}). In the same paper \cite{Jer19}, Jerison also noted that for many acute triangles the second Neumann eigenfunction may fail to be strictly monotone in any direction. In contrast, we show that the second Neumann eigenfunction on a triangle is indeed monotone in at least one direction. This inherent monotonicity property will be preserved through the utilization of continuity methods via domain deformation. 

It is well known and widely used that a constant vector field ${L}$ commutes with the Laplacian, see \cite{BB99, AB04, Siu15}. Moreover, commutativity with the Laplacian also holds for certain nonconstant vector fields such as rotational vector fields. To be more precise, for a point ${p} = ({p}_{1}, {p}_{2}) \in \R^{2}$, we set
\begin{equation*}
\begin{aligned}
{R}_{{p}} & = - ({x}_{2} - {p}_{2})\partial_{{x}_{1}} + ({x}_{1} - {p}_{1})\partial_{{x}_{2}}, 
\\
{R}_{{p}}{u} & = - ({x}_{2} - {p}_{2})\partial_{{x}_{1}} {u} + ({x}_{1} - {p}_{1})\partial_{{x}_{2}} {u}. 
\end{aligned}
\end{equation*}
Here ${R}_{{p}}$ is the rotational vector field generating the counterclockwise rotational flow about ${p}$, and ${R}_{{p}}{u}$ is the angular derivative (or rotational derivative) of ${u}$ about ${p}$. Each rotational vector field ${R}_{{p}}$ commutes with the Laplacian, and hence the angular derivative ${R}_{{p}}({u})$ satisfies the same equation as ${u}$ (without boundary conditions). Constant and rotational vector fields play an important role in studying the number of non-vertex critical points in \cite{JM20, JM22ar}. The angular derivative is also used in \cite{Put90, CLY21}. We shall exploit these commutation properties in this paper. 

Our approach follows the method of continuity in \cite{JN00} as well as in \cite{CLW19, JM20, JM22ar}. Starting from a triangle ${T}^{0}$ for which the second eigenfunction is known to be monotone, we construct a continuously varying family of triangles ${T}^{t}$ with ${T}^{1}$ being the target triangle, and aim to prove the same property for ${T}^{1}$. The initial domain we choose is a special triangle (e.g., the half-equilateral triangle or the right isosceles triangle), for which the corresponding eigenfunction ${u}^{0}$ is explicit and has the desired monotonicity. Assume the monotonicity property persists up to some $\bar{t}$. By the continuous dependence of eigenfunctions on ${t}$, the property then holds at $\bar{t}$ and for ${t} \leq \bar{t}$. For the second eigenfunction on the critical triangle ${T}^{\bar{t}}$ (and also on every ${T}^{t}$ with ${t} \leq \bar{t}$), we will establish two key facts: 
\begin{itemize}
\item
The tangential derivative is nonzero along the interiors of the two longest sides.
\item
The eigenfunction is nonzero at all vertices, except in the superequilateral case. 
\end{itemize}
These two facts play an important role in the proof of \autoref{thm12} and indicate certain stability under domain perturbations. Once they are established, the monotonicity property remains valid beyond ${t} = \bar{t}$ and eventually for all ${t} \in [0, 1]$. 

The uniqueness of a non-vertex critical point plays a pivotal role in the hot spots conjecture. 
For a triangle ${T}^{t}$ with ${t} \leq \bar{t}$, the proof proceeds in three steps. First, by applying the maximum principle to the angular derivative, we exclude critical points on a large portion of the two longest sides (see \autoref{lma31}). Second, using an eigenvalue inequality (see \eqref{CGY0327}) and taking the directional derivative as a test function, we rule out non-vertex critical points on the remaining parts of those two sides (see \autoref{thm38}). Third, a global analysis of the nodal lines of suitable directional derivatives yields the uniqueness of a non-vertex critical point on the shortest side (see \autoref{lma35}).

The rigidity result asserts that the nodal line meets a vertex if and only if the triangle is superequilateral. This is proved via a global analysis of the nodal domain of the sum of the eigenfunction and its reflection across the internal bisector (see \autoref{thm51}). Indeed, if the triangle is not isosceles and the eigenfunction vanishes at a vertex, the monotonicity property implies that the second Neumann eigenvalue $\mu$ exceeds the first mixed eigenvalue on a suitable subtriangle, which appears to contradict the well-known inequality $\mu_{2}(\Omega) < \lambda_{1}(\Omega)$ of P{\'o}lya \cite{Pol52}. To make this rigorous, we establish an eigenvalue inequality conjectured by Siudeja \cite{Siu16} (see \autoref{thm14Eig}). Combining this with elementary lower/upper eigenvalue bounds and isoperimetric inequalities for convex cones \cite{LPT88, LP90}, we obtain a comparison between the second Neumann eigenvalue and the first mixed eigenvalue on a subtriangle (see \autoref{thm47}), leading to a contradiction. 

The uniqueness and saddle-type nature of the non-vertex critical point ensure that the global extrema occur only at the vertices; the remaining task is to determine which vertex realizes the global minimum when two local minima are present. This is nontrivial. A key observation is that the non-vertex critical point must lie close to (or not far from) the vertex with the largest interior angle. This guides an argument akin to the method of moving planes. As a consequence, we determine the exact location of the global extrema: they occur at and only at the endpoints of the longest side. In particular, for non-equilateral triangles that are not subequilateral, the vertex with the largest interior angle is not a global extreme point; see \autoref{Sect7location} for details. Moreover, our method also yields an alternative proof of both the uniqueness of the non-vertex critical point and the monotonicity property, without requiring the exact endpoints of the nodal line $\mathcal{Z}(u)$ (see \autoref{thm51} and \autoref{thm47}).

The proof relies crucially on the monotonicity property, a central theme in the pursuit of the hot spots conjecture. Leveraging this monotonicity, we construct a lower solution, built from suitable directional derivatives, to a linear problem with mixed boundary conditions. This yields a quantitative lower bound for the second Neumann eigenvalue. We then invoke eigenvalue inequalities to obtain a contradiction. This strategy may be useful in future studies, particularly for problems involving monotone eigenfunctions.

The topic of eigenvalue inequalities is both vibrant and significant, and it often emerges from the study of specific spectral problems. Beyond their intrinsic interest, such inequalities help reveal qualitative features of eigenfunctions, as in the nodal line conjecture. In the other direction, our approach offers a complementary perspective. Qualitative information about eigenfunctions, such as  monotonicity, can itself lead to new eigenvalue inequalities, including monotonicity of eigenvalues under domain variations as well as more specific comparison results. This viewpoint offers a promising avenue for further exploration in spectral theory.


\subsection{Outline of the paper} 

In \autoref{Sect2pre}, we recall the maximum principle, eigenvalue inequalities, basic properties of eigenfunctions, and monotonicity for positive solutions of mixed boundary problems. In \autoref{Sect3mon}, we prove the uniqueness of non-vertex critical points for monotone eigenfunctions. In \autoref{Sect4EI}, we establish several eigenvalue inequalities, including \autoref{thm14Eig}. The necessary and sufficient condition for the nodal line to meet a vertex is proved in \autoref{Sect5nodal}. In \autoref{Sect6unique}, we prove \autoref{thm12} except for Property~\ref{CGY0102itd}. Finally, Property~\ref{CGY0102itd} is addressed in \autoref{Sect7location}, where we also provide an alternative proof of the uniqueness of non-vertex critical points and of the monotonicity property. 


\subsection{Notations and terminologies}

We conclude this section by introducing the notation and terminology used throughout the paper. 
$\mathbb{N}^{ + }$ denotes the set of all positive integers, and $\mathbb{R}^{ + }$ denotes the set of all positive real numbers. For a set or domain $\Omega$, the \textbf{closure} of $\Omega$ is denoted by $\overline{\Omega}$, and the \textbf{interior} of $\Omega$ is denoted by $\operatorname{Int}(\Omega)$. For a continuous function ${f}$ defined on the closure $\overline{\Omega}$ of a domain $\Omega$, the nodal line (zero level set) $\mathcal{Z}({f})$ is defined as the closure of $\{{x} \in \Omega: {f}({x}) = 0\}$, and the nodal domains of ${f}$ are the connected components of $\Omega \setminus \mathcal{Z}({f})$. For notational convenience, we often identify the Euclidean plane with the complex plane $\mathbb{C}$, and we write ${z} = {x}_{1} + {x}_{2}\sqrt{ - 1}$ to represent a point $({x}_{1}, {x}_{2})$ on the plane. In particular, ${x}_{1} = \operatorname{Re}({z})$ and ${x}_{2} = \operatorname{Im}({z})$; if ${z} = {r} \exp(\sqrt{ - 1} \cdot \theta)$, then $\theta = \operatorname{arg}({z})$ and ${r} = |{z}|$. For two points ${z}_{1}$ and ${z}_{2}$ in $\R^{2}$ or $\mathbb{C}$, the symbol $\overline{{z}_{1}{z}_{2}}$ denotes the closed line segment connecting ${z}_{1}$ and ${z}_{2}$, $\operatorname{Int}(\overline{{z}_{1}{z}_{2}})$ denotes the interior of $\overline{{z}_{1}{z}_{2}}$, and $|\overline{{z}_{1}{z}_{2}}|$ denotes the length of $\overline{{z}_{1}{z}_{2}}$. 

A side of a polygon is always understood as a closed, undirected line segment. The symbol $\triangle$ denotes an open triangle, and ${T} = \triangle {z}_{1}{z}_{2}{z}_{3}$ denotes an open triangle with vertices at ${z}_{1}$, ${z}_{2}$ and ${z}_{3}$. 
The notation $\angle {z}_{1}{z}_{2}{z}_{3}$ denotes the interior angle between the two sides $\overline{{z}_{1}{z}_{2}}$ and $\overline{{z}_{3}{z}_{2}}$. We use $\alpha_{1}$, $\alpha_{2}$, and $\alpha_{3}$ to denote the interior angles of the triangle at vertices ${z}_{1}$, ${z}_{2}$ and ${z}_{3}$, respectively. That is, 
\begin{equation} \label{CGY0106angle} 
\alpha_{1} = \angle {z}_{3}{z}_{1}{z}_{2}, \quad
\alpha_{2} = \angle {z}_{1}{z}_{2}{z}_{3}, \quad \text{and} \quad
\alpha_{3} = \angle {z}_{2}{z}_{3}{z}_{1}. 
\end{equation}
Moreover, 
\begin{itemize}
\item
The symbols $\boldsymbol{\tau}_{S}$, $\boldsymbol{\tau}_{M}$, and $\boldsymbol{\tau}_{L}$ denote the unit tangential vectors to the sides $\overline{{z}_{1}{z}_{2}}$, $\overline{{z}_{3}{z}_{1}}$, and $\overline{{z}_{2}{z}_{3}}$ pointing from the first vertex to the seocnd one in the above notation of each side, respectively; 
\item
The symbols $\boldsymbol{n}_{S}$, $\boldsymbol{n}_{M}$, and $\boldsymbol{n}_{L}$ denote the unit inward normal vectors to the sides $\overline{{z}_{1}{z}_{2}}$, $\overline{{z}_{3}{z}_{1}}$, and $\overline{{z}_{2}{z}_{3}}$, respectively. 
\end{itemize}
Consequently, the pairs $(\boldsymbol{\tau}_{S}, \boldsymbol{n}_{S})$, $(\boldsymbol{\tau}_{M}, \boldsymbol{n}_{M})$, and $(\boldsymbol{\tau}_{L}, \boldsymbol{n}_{L})$ form right-handed coordinate systems in $\R^{2}$. We note that $\overline{{z}_{1}{z}_{2}}$ is taken to be the shortest side in \autoref{Sect6unique}, whereas $\overline{{z}_{2}{z}_{3}}$ is taken to be the longest side in \autoref{Sect7location}. 


\section{Preliminaries} \label{Sect2pre}

In this section, we review results on the maximum principle, eigenvalue inequalities, fundamental properties of the second Neumann eigenfunction, and the monotonicity of positive elliptic solutions for mixed boundary problems. 


\subsection{Maximum principles}

Let $\Omega$ be a bounded Lipschitz domain in $\R^{n}$. Let $\Gamma_{0}$ be a closed subset of $\partial\Omega$, and define $\Gamma_{1}: = \partial\Omega \setminus \Gamma_{0}$. Consider the mixed boundary value problem
\begin{equation} \label{CGY0202a}
\mathcal{L}{u} = {f} \mbox{ in } \Omega \text{ and } 
\mathcal{B}{u} = {g} \mbox{ on } \partial\Omega, 
\end{equation}
where the differential operator $\mathcal{L}$ and the boundary operator $\mathcal{B}$ are defined by
\begin{gather*}
\mathcal{L} {u} = \Delta {u} + {c}({x}) {u}, 
\\
\mathcal{B} {u} = {u} \text{ on } \Gamma_{0} \text{ and } \mathcal{B} {u} = \partial_{\nu} {u} \text{ on } \Gamma_{1}, 
\end{gather*}
respectively, with $\|{c}\|_{L^{\infty}(\Omega)} < {c}_{0}$ for some ${c}_{0} \in \R^{ + }$. It is well known (see \cite{GT83, PW67}) that, under the condition ${c} < 0$ in $\Omega$, we have 
\begin{equation*}
\text{If } {f} \leq 0 \text{ in } \Omega \text{ and } {g} \geq 0 \text{ on } \partial \Omega, \text{ then } {u} \geq 0 \text{ in } \Omega. 
\end{equation*}
This property is called the \textit{maximum principle}. In the case of Dirichlet boundary conditions (i.e., $\Gamma_{0} = \partial\Omega$), this reduces to the classical maximum principle, which has been extensively studied (see \cite{BN91, BNV94, GT83, PW67}). Berestycki, Nirenberg, and Varadhan \cite{BNV94} introduced the \textit{refined first eigenvalue} $\lambda_{1}(\mathcal{L}, \Omega)$, characterized by
\begin{equation}
\lambda_{1}(\mathcal{L}, \Omega) = \sup\{\lambda \in \R: \, \exists \phi > 0 \text{ in } \Omega \text{ such that } (\mathcal{L} + \lambda)\phi \leq 0 \text{ in } \Omega\}
\end{equation}
where $\phi \in W_{\text{loc}}^{2, n}(\Omega) \cap C(\overline{\Omega})$. Using this definition, the following fundamental result was established. 

\begin{lma}[Berestyski et al., 1994, \cite{BNV94}] \label{lma21BNV}
Let $\Gamma_{0} = \partial\Omega$. Then 
\begin{enumerate}
\item[\rm(1)]
The maximum principle holds for $\mathcal{L}$ in $\Omega$ if and only if $\lambda_{1}(\mathcal{L}, \Omega) > 0$. 

\item[\rm(2)]
The maximum principle holds for $\mathcal{L}$ in $\Omega$ provided that there exists a function $\phi \in W_{\mathrm{loc}}^{2, n}(\Omega) \cap C(\overline{\Omega})$ which is positive in $\Omega$ and satisfies 
\begin{equation*}
\text{either } \; \mathcal{L}\phi \leq, \not\equiv 0 \text{ in } \Omega, \quad \text{or} \quad \mathcal{L}\phi \equiv 0 \ \text{in } \Omega \text{ and } \phi \not\equiv 0 \ \text{on } \partial\Omega.
\end{equation*}
\end{enumerate}
\end{lma}

For mixed boundary conditions, we obtain an analogous result.

\begin{lma}[Lemma 6 of \cite{YCL18}] \label{lma22aMP}
The maximum principle holds for $(\mathcal{L}, \mathcal{B})$ in $\Omega$ if and only if $\lambda_{1}^{\rm{Mix}} > 0$, where $\lambda_{1}^{\rm{Mix}}$ denotes the principal eigenvalue of the mixed boundary value problem
\begin{equation} \label{CGY0202b}
\begin{cases}
\mathcal{L} \varphi + \lambda \varphi = 0 & \text{in } \Omega, \\
\varphi = 0 & \text{on } \Gamma_{0}, \\
\partial_{\nu} \varphi = 0 & \text{on } \Gamma_{1}. 
\end{cases}
\end{equation}
\end{lma}

We next consider the homogeneous equation associated with \eqref{CGY0202a}: 
\begin{equation} \label{CGY0202c}
\mathcal{L} {v} = 0 \text{ in } \Omega \text{ and } \mathcal{B} {v} = 0 \text{ on } \partial\Omega. 
\end{equation}
Let $\mathcal{H}$ denote the closure of $C_{c}^{1}(\Omega \cup \Gamma_{1})$ in ${H}^{1}(\Omega)$. We define a \textit{weak upper solution} to \eqref{CGY0202c} as a function ${v} \in {H}^{1}(\Omega)$ satisfying ${v}^{ - } = \min\{{v}, 0\} \in \mathcal{H}$ and
\begin{equation*}
\int_{\Omega} (\nabla {v} \cdot \nabla \phi - {c} {v} \phi) dx \geq 0
\end{equation*}
for all for all $\phi \in \mathcal{H}$ with $0 \leq \phi$. Similarly, ${v}$ is a \textit{weak lower solution} if $ - {v}$ is a weak upper solution. A \textit{weak solution} is a function that is both a weak upper and a weak lower solution. 

\begin{lma} \label{lma22bMP}
If there exists a nonnegative weak lower solution ${v} \in {H}^{1}(\Omega) \setminus \{0\}$ of \eqref{CGY0202c}, then $\lambda_{1}^{\rm{Mix}} \leq 0$. Moreover, $\lambda_{1}^{\rm{Mix}} < 0$ if ${v}$ is not a solution to \eqref{CGY0202c}. 

If there exists a nonnegative weak upper solution ${v} \in {H}^{1}(\Omega) \setminus \{0\}$ of \eqref{CGY0202c}, then $\lambda_{1}^{\rm{Mix}} \geq 0$. Moreover, $\lambda_{1}^{\rm{Mix}} > 0$ if ${v}$ is not a solution to \eqref{CGY0202c}. 
\end{lma}

\begin{proof}
Let $\varphi > 0$ be the principal eigenfunction of \eqref{CGY0202b}. Then we have
\begin{equation*}
\lambda_{1}^{\rm{Mix}} \int_{\Omega} {v} \varphi dx = \int_{\Omega} (\nabla {v} \cdot \nabla \varphi - {c} {v} \varphi ) dx. 
\end{equation*}
The desired inequalities follow directly from the assumptions on ${v}$. 
\end{proof}

As a consequence of \autoref{lma22bMP}, we obtain preliminary bounds on the eigenvalues.
\begin{lma} \label{lma22cMP} 
If ${v} \geq, \not\equiv 0$ is a weak lower (resp., upper) solution of 
\begin{equation*} 
\Delta {v} + \mu {v} = 0 \text{ in } \Omega, \quad {v} = 0 \text{ on } \Gamma_{0} \text{ and } \partial_{\nu}{v} = 0 \text{ on } \Gamma_{1}
\end{equation*}
for some constant $\mu \in \R$, then the first mixed eigenvalue $\lambda_{1}(\Omega, \Gamma_{0})$ (as defined in \eqref{CGY0104Mix}) satisfies 
\begin{equation*}
\lambda_{1}(\Omega, \Gamma_{0}) \leq \mu \quad \text{(resp., } \lambda_{1}(\Omega, \Gamma_{0}) \geq \mu \text{)}.
\end{equation*}
\end{lma}


\subsection{Eigenvalue inequalities}

\begin{lma}[\cite{Pol52, Sze54}] \label{lma24Pol}
For any bounded Lipschitz domain $\Omega$ in $\R^{2}$, the second Neumann eigenvalue $\mu_{2}(\Omega)$ is strictly less than the first Dirichlet eigenvalue $\lambda_{1}(\Omega)$. Consequently, for any simply connected Lipschitz domain $\Omega$ in $\R^{2}$, the nodal line of the second Neumann eigenfunction is a smooth, simple curve intersecting the boundary at exactly two points. 
\end{lma}

\begin{lma}[Theorem 3.1 of \cite{LR17}] \label{lma25LR}
Let $\Omega$ be a Lipschitz domain in $\R^{n}$ such that $\Gamma_{N} = \partial\Omega \setminus \Gamma_{D}$ is flat (i.e., $\Gamma_{N}$ is contained in some hyperplane). Then the second Neumann eigenvalue $\mu_{2}(\Omega)$ is no greater than the first mixed eigenvalue $\lambda_{1}(\Omega, \Gamma_{D})$. 
\end{lma}

\begin{lma}[Corollary 5.5 of \cite{JM20}] \label{lma26aJM}
The second Neumann eigenfunction on a triangle cannot vanish at two vertices of the triangle. 
\end{lma}


\subsection{Basic properties of eigenfunctions}

We present results on the local behavior of eigenfunctions near a vertex of a polygon \cite{JN00, JM20, JM22ar, JM22c}. Following \cite{JM20}, let ${e}$ be a side of a polygon ${P}$ (viewed as a closed, undirected line segment), and let ${L}_{{e}}$ denote the constant unit tangent vector field parallel to ${e}$, oriented so that its counterclockwise rotation by $\pi/2$ points into the interior of ${P}$.

\begin{lma} \label{lma26bJM}
Let ${p}$ be an acute vertex of a polygon ${P}$ and let ${u}$ be a second Neumann eigenfunction of ${P}$. Then the following conclusions hold: 
\begin{enumerate}[label = \rm(\arabic*)]
\item
The vertex ${p}$ is a local extremum of ${u}$ if and only if ${u}({p}) \neq 0$. Moreover, the vertex ${p}$ is a local maximum (resp., minimum) point of ${u}$ if and only if ${u}({p}) > 0$ (resp., ${u}({p}) < 0$). 
\item
If ${P}$ is an acute triangle, ${e}$ is the side opposite to ${p}$ and ${e}'$ is a side adjacent to ${p}$, then there exists a small neighborhood $\mathcal{O}_{{p}}$ of ${p}$ such that
\begin{enumerate}[label = \rm(2\alph*)]
\item
If ${u}({p}) \neq 0$, then $\mathcal{Z}({L}_{{e}}{u}) \cap \mathcal{O}_{{p}}$ is a real-analytic curve that meets ${p}$, and the tangent line of $\mathcal{Z}({L}_{{e}}{u})$ at ${p}$ is not aligned with the sides adjacent to ${p}$. 
\item
If ${u}({p}) \neq 0$, then ${L}_{{e}'}{u} \neq 0$ in $\mathcal{O}_{{p}} \cap (\overline{{T}} \setminus \{{p}\})$. 
\item
If ${u}({p}) = 0$, then ${L}_{{e}}{u} \neq 0$ in $\mathcal{O}_{{p}} \cap (\overline{{T}} \setminus \{{p}\})$. 
\item
If ${u}({p}) = 0$, then $\mathcal{Z}({L}_{{e}'}{u}) \cap \mathcal{O}_{{p}}$ is a real-analytic curve that meets ${p}$, and the tangent line of $\mathcal{Z}({L}_{{e}'}{u})$ at ${p}$ is not aligned with the sides adjacent to ${p}$. 
\end{enumerate}
\end{enumerate}
\end{lma}

\begin{proof}
This result follows from the local Fourier-Bessel expansion of the eigenfunction near the vertex. The first statement is derived from \cite[Proposition 2.1]{JM22ar}, while the second statement is established in \cite[Lemma 2.2]{JM22ar} and \cite[Lemma 2.1]{JM22c}. 
\end{proof}

The fundamental results on eigenfunctions for isosceles triangles and non-acute triangles are summarized as follows. 

\begin{prop} \label{prop27iso}
Let ${T}$ be a non-equilateral isosceles triangle, and let ${u}$ be a second Neumann eigenfunction on ${T}$.
Then ${u}$ is symmetric if and only if ${T}$ is subequilateral, and antisymmetric if and only if ${T}$ is superequilateral.

{\rm(i)}
If ${T}$ is subequilateral, then ${u}$ has exactly one non-vertex critical point at the midpoint of the shortest side, is monotone in the direction perpendicular to the shortest side, and its nodal line $\mathcal{Z}({u})$ intersects the interior of each of the two longest sides.

{\rm(ii)}
If ${T}$ is superequilateral, then ${u}$ has no non-vertex critical points, is monotone in the direction perpendicular to the shortest side, and its nodal line $\mathcal{Z}({u})$ coincides with the symmetry axis. 
\end{prop}

\begin{proof}
The symmetry and antisymmetry properties of ${u}$ are established in \cite{LS10}, while properties (i) and (ii) are derived from \cite[Theorem B and Theorem C]{Miy13}. 
\end{proof}

\begin{prop}\label{prop28obt}
Let ${T} = \triangle {z}_{1}{z}_{2}{z}_{3}$ with $\angle {z}_{3}{z}_{1}{z}_{2} \geq \pi/2$. Let ${u}$ be a second Neumann eigenfunction of ${T}$. Then
\begin{enumerate}[label = \rm(\roman*)]
\item
The second Neumann eigenvalue is simple. 
\item
${u}$ does not have any non-vertex critical points, and ${u}$ is monotone in two fixed directions that are perpendicular to the shortest and the intermediate sides, respectively. 
\item
${u}$ attains its global extrema exclusively at the endpoints of the longest side of ${T}$. 
\end{enumerate}
\end{prop}

\begin{proof}
Since non-acute triangles are lip domains as defined in \cite{AB04}, the simplicity of the second Neumann eigenvalue follows from \cite[Theorem 1(i)]{AB04}. Following the proof of \cite[Theorem 1]{AB04} (see also \cite{Roh21}) and possibly after a sign change, ${u}$ has the following monotonicity properties: $\nabla {u} \cdot \boldsymbol{n}_{S} > 0$ and $\nabla {u} \cdot \boldsymbol{n}_{M} < 0$, which further imply 
\begin{equation} \label{CGY0205}
\nabla {u} \cdot \boldsymbol{\tau}_{S} < 0, \; \nabla {u} \cdot \boldsymbol{\tau}_{M} < 0 \text{ and } \nabla {u} \cdot \boldsymbol{\tau}_{L} > 0 \text{ in } {T}. 
\end{equation}
Here, the tangential and inward normal vectors are defined as in the end of \autoref{Sect1intro}. 

The absence of non-vertex critical points follows directly from \eqref{CGY0205}. In fact, suppose by contradiction that there exists a point ${p} \in \operatorname{crit}_{\mathrm{nv}}({u}) \cap \overline{{z}_{1}{z}_{2}}$. Noting that $\nabla {u} \cdot \boldsymbol{\tau}_{S}$ attains its maximum zero at ${p}$, the Hopf boundary lemma implies that the second-order mixed derivative $\partial_{\boldsymbol{\tau}_{S}\boldsymbol{n}_{S}} {u}$ is negative at ${p}$. However, the Neumann boundary condition for ${u}$ ensures that $\partial_{\boldsymbol{\tau}_{S}\boldsymbol{n}_{S}} {u}$ vanishes on $\operatorname{Int}(\overline{{z}_{1}{z}_{2}})$, leading to a contradiction. Therefore, there are no non-vertex critical points on $\overline{{z}_{1}{z}_{2}}$. 
Similarly, there are no non-vertex critical points on either of the remaining two sides. 
Thus, we conclude
\begin{equation*}
\nabla {u} \cdot \boldsymbol{\tau}_{S} < 0, \; \nabla {u} \cdot \boldsymbol{\tau}_{M} < 0 \text{ and } \nabla {u} \cdot \boldsymbol{\tau}_{L} > 0 \text{ in } \overline{{T}} \setminus \{\text{Vertices}\}. 
\end{equation*}
This implies that ${z}_{2}$ is the unique global minimum point of ${u}$ and ${z}_{3}$ is the unique global maximum point of ${u}$. This completes the proof. 
\end{proof}


\subsection{Monotonicity of positive solutions to mixed boundary problems}

\begin{prop} \label{prop29Yao}
Let $\Omega$ be a triangle in the plane and let $\Gamma_{D}$ be a side of the triangle $\Omega$. Suppose ${u}$ is a positive solution to 
\begin{equation*} 
\begin{cases}
\Delta {u} + {f}({u}) = 0 & \mbox{ in } \Omega, \\
{u} = 0 & \mbox{ on }\Gamma_{D}, \\
\partial_{\nu}{u} = 0 & \mbox{ on } \Gamma_{N} = \partial\Omega \setminus \Gamma_{D}, 
\end{cases}\end{equation*}
where ${f}$ is a locally Lipschitz continuous function on $\R$. Let $\boldsymbol{n}$ denote the unit inward normal vector to $\Gamma_{D}$. Then the following holds: 
\begin{enumerate}
\item[\rm(1)]
If the interior angle between the two Neumann sides is acute or right, then $\nabla {u} \cdot \boldsymbol{n} > 0$ in $\Omega$. 

\item[\rm(2)]
If the two Neumann sides have equal lengths, then the solution ${u}$ is symmetric with respect to the internal angle bisector of the angle between the two Neumann sides, and, in particular, $\nabla {u} \cdot \boldsymbol{n} > 0$ in $\Omega$. 
\end{enumerate}
\end{prop}

\begin{proof}
This result is proved by the method of moving planes (see \cite{LY24}).
\end{proof}

We provide a sketch of the proof for the first mixed eigenfunction in \autoref{lma41} below. This monotonicity property plays a pivotal role in establishing key eigenvalue inequalities (\autoref{thm14Eig}) in this paper. 


\section{Properties of monotone eigenfunctions} \label{Sect3mon}

In this section, we primarily focus on acute triangles and establish the uniqueness of the non-vertex critical point of the second Neumann eigenfunction, assuming an additional monotonicity condition. 

Throughout this section, the vertices of the triangle ${T} = \triangle {z}_{1}{z}_{2}{z}_{3}$ are arranged so that 
\begin{itemize}
\item
the vertices ${z}_{1}$ and ${z}_{2}$ lie on the real axis, with ${z}_{1}$ positioned to the left of ${z}_{2}$; 
\item
the vertex ${z}_{3}$ lies above the real axis, see \autoref{fig1}. 
\end{itemize}
We always denote by $\mu = \mu_{2}(T)$ the second Neumann eigenvalue of ${T}$, and by ${u}$ a second Neumann eigenfunction of ${T}$. Furthermore, we assume that ${u}$ is monotone in the direction normal to the side $\overline{{z}_{1}{z}_{2}}$, that is, the derivative $\partial_{{x}_{2}}{u}$ satisfies 
\begin{equation} \label{CGY0301}
\partial_{{x}_{2}}{u} \geq 0 \text{ in } {T}. 
\end{equation}
This monotonicity condition will be ultimately verified in \autoref{Sect6unique} when $\overline{{z}_{1}{z}_{2}}$ is the shortest side. 


\subsection{The sign of the angular derivative}

For a point ${p} = ({p}_{1}, {p}_{2}) \in \R^{2}$, the \emph{rotational derivative} about the point ${p}$ is given by
\begin{equation}
{R}_{{p}}{u}({x}) = - ({x}_{2} - {p}_{2})\partial_{{x}_{1}}{u}({x}) + ({x}_{1} - {p}_{1})\partial_{{x}_{2}}{u}({x}). 
\end{equation}
This rotational derivative satisfies the same equation as ${u}$; specifically, ${R}_{{p}}{u}$ satisfies
\begin{equation} \label{CGY0303b}
\Delta {R}_{{p}}{u} + \mu {R}_{{p}}{u} = 0 \text{ in } {T}. 
\end{equation}

\begin{lma} \label{lma31}
Suppose that ${u}$ satisfies the monotonicity condition \eqref{CGY0301}. Then the following statements hold: 
\begin{enumerate}[label = \rm(\roman*)]
\item
${u}$ is strictly monotone in the variable ${x}_{2}$, namely, 
\begin{equation} \label{CGY0304}
\partial_{{x}_{2}} {u} > 0 \text{ in } {T}. 
\end{equation}
\item
The angular derivative ${R}_{{p}}{u}$ of ${u}$ satisfies
\begin{subequations}
\begin{align}
\label{CGY0305a}
{R}_{{p}}{u} > 0 \text{ in } {T} \cap \{{x}_{2} > {p}_{2}\} \text{ whenever } {p} = ({p}_{1}, {p}_{2}) \in \overline{{z}_{3}{z}_{1}}, 
\\ \label{CGY0305b}
{R}_{{p}}{u} < 0 \text{ in } {T} \cap \{{x}_{2} > {p}_{2}\} \text{ whenever } {p} = ({p}_{1}, {p}_{2}) \in \overline{{z}_{3}{z}_{2}}. 
\end{align}
\end{subequations}
\item
The directional derivatives satisfy
\begin{subequations} \label{CGY0306}
\begin{align} \label{CGY0306a}
\nabla {u} \cdot \boldsymbol{\tau}_{M} < 0 & \text{ in } \{{x} \in \operatorname{Int}(\overline{{z}_{3}{z}_{1}}) \cup {T}: \, ({x} - {z}_{2}) \cdot \boldsymbol{\tau}_{M} \leq 0\}, 
\\ \label{CGY0306b}
\nabla {u} \cdot \boldsymbol{\tau}_{L} > 0 & \text{ in } \{{x} \in \operatorname{Int}(\overline{{z}_{3}{z}_{2}}) \cup {T}: \, ({x} - {z}_{1}) \cdot \boldsymbol{\tau}_{L} \geq 0\}, 
\end{align}
\end{subequations}
where the tangential directions $\boldsymbol{\tau}_{M}$ and $\boldsymbol{\tau}_{L}$ are defined below \eqref{CGY0106angle}. 
\end{enumerate}
\end{lma}

\begin{proof}
We establish the result in three parts. 

\textbf{Part (i): Positivity of $\partial_{{x}_{2}}{u}$}. 
First, we observe that $\partial_{{x}_{2}}{u}$ cannot be identically zero within the triangle ${T}$. Suppose, for contradiction, that $\partial_{{x}_{2}} {u} \equiv 0$ in ${T}$. This would imply that ${u}$ depends solely on the variable ${x}_{1}$. Because of the Neumann boundary condition for ${u}$, ${T}$ would have to be a rectangle, which is impossible. Therefore, we conclude that $\partial_{{x}_{2}} {u} \not\equiv 0$ in ${T}$. Since the nonnegative function $\partial_{{x}_{2}} {u}$ satisfies the same linear equation as ${u}$, namely, 
\begin{equation*} 
\Delta (\partial_{{x}_{2}} {u}) + \mu \partial_{{x}_{2}} {u} = 0 \text{ in } {T}, 
\end{equation*} 
the strong maximum principle implies that $\partial_{{x}_{2}} {u}$ is strictly positive in the interior of ${T}$. This validates inequality \eqref{CGY0304}. 

\textbf{Part (ii): The signs of the angular derivative ${R}_{{p}}{u}$}. 
Let ${p} = ({p}_{1}, {p}_{2})$ be any point on the left side $\overline{{z}_{3}{z}_{1}}$ of the triangle ${T}$. The Neumann boundary condition for ${u}$ implies that the angular derivative ${R}_{{p}}{u}$ vanishes on the left side $\overline{{z}_{3}{z}_{1}}$. The inequality \eqref{CGY0301} further guarantees the nonnegativity of ${R}_{{p}}{u}$ on $\overline{T} \cap \{{x}_{2} = {p}_{2}\}$. Moreover, the Neumann boundary condition combined with \eqref{CGY0301} tells us that
\begin{equation*}
\nabla {u} = |\nabla {u}| \boldsymbol{\tau}_{L} \text{ on } \overline{{z}_{3}{z}_{2}}. 
\end{equation*}
This ensures the nonnegativity of ${R}_{{p}}{u}$ on the right side $\overline{{z}_{3}{z}_{2}}$. Therefore, the angular derivative ${R}_{{p}}{u}$ is continuous in $\overline{T}$; it satisfies the linear equation \eqref{CGY0303b} and
\begin{equation*}
{R}_{{p}}{u} \geq 0 \text{ on } \partial ({T} \cap \{{x}_{2} > {p}_{2}\} ). 
\end{equation*}
From \autoref{lma24Pol}, 
\begin{equation} 
\lambda_{1}({T}) > \mu_{2}({T}) = \mu. 
\end{equation}
Applying the maximum principle as in \autoref{lma21BNV}, we deduce that ${R}_{{p}}{u} \geq 0$ in ${T} \cap \{{x}_{2} > {p}_{2}\}$. Furthermore, ${R}_{{p}}{u}$ cannot be identically zero. If it were, ${u}$ would be radial about the point ${p}$, which would contradict the Neumann boundary condition on $\overline{{z}_{3}{z}_{2}}$. Therefore, the strong maximum principle implies \eqref{CGY0305a}. Similarly, \eqref{CGY0305b} follows by an analogous argument.

\textbf{Part (iii): The non-vanishing of the directional derivative along $\boldsymbol{\tau}_{M}$ and $\boldsymbol{\tau}_{L}$ in a certain range}. 
The inequality \eqref{CGY0306} holds when ${T}$ is a non-acute triangle (see \autoref{prop28obt}); we therefore focus only on the acute triangle ${T}$. 
From \eqref{CGY0305b}, for every ${p} \in \overline{{z}_{3}{z}_{2}}$ we have
\begin{equation*}
{R}_{{p}}{u} < 0 \text{ on } \{{x} \in {T}: \, ({x} - {p}) \cdot \boldsymbol{\tau}_{M} = 0\}, 
\end{equation*}
which is equivalent to
\begin{equation*}
\nabla {u} \cdot \boldsymbol{\tau}_{M} < 0 \text{ on } \{{x} \in {T}: \, ({x} - {p}) \cdot \boldsymbol{\tau}_{M} = 0\}. 
\end{equation*}
Due to the arbitrariness of ${p} \in \overline{{z}_{3}{z}_{2}}$, we obtain
\begin{equation} \label{CGY0308}
\nabla {u} \cdot \boldsymbol{\tau}_{M} < 0 \text{ on } \{{x} \in {T}: \, ({x} - {z}_{2}) \cdot \boldsymbol{\tau}_{M} \leq 0\}. 
\end{equation}
Now, we reflect ${u}$ across the left Neumann boundary $\overline{{z}_{3}{z}_{1}}$ to obtain an extended function $\tilde{u}$ defined on the kite-shaped domain ${K}$ with vertices ${z}_{3}$, ${z}_{2}$, ${z}_{1}$, and ${z}_{2}'$. Combining this with \eqref{CGY0308}, we derive that $\nabla \tilde{u} \cdot \boldsymbol{\tau}_{M} \leq, \not\equiv 0$ in the triangle ${z}_{3}{z}_{2}{z}_{2}'$. Applying the strong maximum principle to $\nabla \tilde{u} \cdot \boldsymbol{\tau}_{M}$, we obtain the strict negativity of $\nabla {u} \cdot \boldsymbol{\tau}_{M}$ in $\triangle {z}_{3}{z}_{2}{z}_{2}'$. This implies that $|\nabla {u}({p})| > 0$ for a point ${p} \in \operatorname{Int}(\overline{{z}_{3}{z}_{1}})$ with $\angle {z}_{3}{p}{z}_{2} > \pi/2$. 
Based on \eqref{CGY0308} and the analysis of the nodal line of $\nabla {u} \cdot \boldsymbol{\tau}_{M}$ in \autoref{lma33}, we have $\nabla {u} \cdot \boldsymbol{\tau}_{M} < 0$ at ${p} \in \operatorname{Int}(\overline{{z}_{3}{z}_{1}})$ with $\angle {z}_{3}{p}{z}_{2} = \pi/2$. This verifies \eqref{CGY0306a}. Similarly, by reflecting the solution ${u}$ across the right side $\overline{{z}_{3}{z}_{2}}$ and applying analogous arguments, we deduce \eqref{CGY0306b}. 
This completes the proof. 
\end{proof}

\begin{rmk} 
The inequality $\lambda_{1}({T}) > \mu_{2}({T})$ can be derived from the monotonicity of ${u}$. Indeed, $\partial_{{x}_{2}}{u}>0$ satisfies the same equation as ${u}$ in ${T}$ and is not identically zero on $\partial {T}$ (see, e.g., \autoref{lma33} below); hence, by \autoref{lma21BNV}, $\lambda_{1}({T}) > \mu = \mu_{2}({T})$.
\end{rmk}


\subsection{The number of non-vertex critical points}

In this subsection, we study the number of non-vertex critical points through local and global analysis of certain directional derivatives of ${u}$. We define
\begin{equation}
\operatorname{crit}_{\mathrm{nv}}({u}) = \{{p} \in \overline{T}: \, |\nabla {u} ({p})| = 0, \; \text{${p}$ is not a vertex of ${T}$} \}. 
\end{equation}

We begin by reviewing basic local properties of nodal lines for functions associated with solutions of elliptic equations. 

\begin{lma} \label{lma33}
Assume that ${L}_{{e}}{u}({p}) = 0$ for some point ${p}$ in the interior of a side ${e}$ of ${T}$, where the unit vector field ${L}_{{e}}$ is as defined before \autoref{lma26bJM}. 
Then there exists ${n} \in \mathbb{N}^{ + }$ such that, at the point ${p}$, the tangential derivatives of ${u}$ vanish up to order ${n}$, while the tangential derivative of order $({n}+1)$ is nonzero. 
Futhermore, the nodal line $\mathcal{Z}({L}_{{e}}{u})$ has exactly ${n}$ real-analytic branches near ${p}$ for some ${n} \in \mathbb{N}^{ + }$, and the angles between the tangent line of each curve of $\mathcal{Z}({L}_{{e}}{u})$ at ${p}$ and the side ${e}$ are $(2{j} - 1)\pi/(2{n})$, ${j} = 1, 2, \ldots, {n}$. 
In particular, if ${u}$ is monotone along ${e}$ near ${p}$, then the number ${n}$ must be even. 
\end{lma}

\begin{proof}
Let $\tilde{{u}}$ be the extension of ${u}$ to the interior of the kite ${K}_{{e}}$ obtained by reflecting across ${e}$. Without loss of generality, we assume that ${e}$ lies on the ${x}_{1}$-axis, ${T}$ lies in $\{{x}_{2} > 0\}$, and ${p}$ is the origin. Then ${L}_{{e}}{u}$ coincides with $\partial_{{x}_{1}}{u}$. Thus, we have
\begin{equation*}
\Delta \partial_{{x}_{1}} \tilde{u} + \mu \partial_{{x}_{1}} \tilde{u} = 0 \text{ in } {T} \text{ and } \partial_{{x}_{2}}(\partial_{{x}_{1}} \tilde{u}) = 0 \text{ on } {e} \subset \{{x}_{2} = 0\}. 
\end{equation*}
It is clear that $\partial_{{x}_{1}} \tilde{u}$ cannot vanish identically in ${T}$, otherwise, ${u}$ would be a function of a single variable, which contradicts both the sign-changing property of ${u}$ and the fact that the domain is a triangle. 
We now use polar coordinates centered at ${p}$, with the polar axis along ${e}$; then ${T}$ is contained in $\theta \in (0, \pi)$. From \cite{HW53} or \cite{HHHO99}, it follows that
\begin{equation*}
{L}_{{e}}{u} = {c} {r}^{{n}}\cos({n}\theta) + O({r}^{{n} + 1}),
\end{equation*}
and 
\begin{equation*}
(\partial_{{x}_{1}})^{j}{u}({p}) = 0 \quad \text{for } 1 \leq {j} \leq {n}, \qquad (\partial_{{x}_{1}})^{{n}+1}{u}({p}) = {n}! \cdot {c},
\end{equation*}
for some ${n} \in \mathbb{N}^{ + }$ and ${c} \neq 0$. Consequently, the nodal line $\mathcal{Z}({L}_{{e}}{u})$ has exactly ${n}$ curves emanating from ${p}$, described by $\theta = (2{j} - 1)\pi/(2{n}) + O({r})$ as ${r} \to 0$. This implies that each curve forms an angle of $(2{j} - 1)\pi/(2{n})$ with the side ${e}$ at ${p}$. 

If ${L}_{{e}}{u} \geq 0$ on ${e}$ and near ${p}$, then both $\cos(0)$ and $\cos({n}\pi)$ must have the same sign. It follows that ${n}$ must be an even number. Hence in a neighborhood of ${p}$, there are at least two curves of $\mathcal{Z}({L}_{{e}}{u})$ that intersect the side ${e}$ transversely at ${p}$. 
\end{proof}

\begin{lma} \label{lma34}
Let ${p}$ be a point in the interior of a side ${e}$ of ${T}$. Let $(\boldsymbol{\tau}, \boldsymbol{n})$ denote the tangential and unit inward normal vectors with respect to the side ${e}$. The ordered pair $(\boldsymbol{\tau}, \boldsymbol{n})$ forms a right-handed orthonormal frame. 
Suppose that $\partial_{\boldsymbol{\tau}}{u}({p}) = 0$. Then there exists the smallest positive integer ${n}$ such that $(\partial_{\boldsymbol{\tau}})^{{j}}{u}({p}) = 0$ for every $1 \leq {j} \leq {n}$, while $(\partial_{\boldsymbol{\tau}})^{{n} + 1}{u}({p}) \neq 0$. 

Suppose further that $\partial_{\boldsymbol{n}\boldsymbol{n}}{u}({p}) \neq 0$. For ${c}_{1}, {c}_{2} \in \R \setminus \{0\}$, we set ${L}{u} = {c}_{1}\partial_{\boldsymbol{\tau}}{u} + {c}_{2}\partial_{\boldsymbol{n}}{u}$ and define the nodal line $\mathcal{Z}({L}{u})$ as the closure of $\{{x} \in {T}: {L}{u}({x}) = 0\}$. Set 
\begin{equation*}
{c}_{3} = \frac{{c}_{1}}{{c}_{2}} \cdot
\frac{(\partial_{\boldsymbol{\tau}})^{{n}+1} {u}({p})}{\partial_{\boldsymbol{n}\boldsymbol{n}} {u}({p})}.
\end{equation*}
Then in a small neighborhood $\mathcal{O}_{{p}}$ of ${p}$, the following statements hold: 
\begin{enumerate}[label = \rm(\arabic*), start = 1]
\item
If $n$ is odd, then
$\mathcal{Z}({L}{u}) \cap \mathcal{O}_{{p}}$ is a real-analytic curve with ${p}$ as an endpoint, it is tangent to ${e}$ at ${p}$ when ${n} > 1$, and intersects ${e}$ transversely when ${n} = 1$. 
\item
If $n$ is even and ${c}_{3} > 0$, then
${L}{u} \neq 0$ in $\mathcal{O}_{{p}} \setminus \{{p}\}$. 
\item
If $n$ is even and ${c}_{3} < 0$, 
then $\mathcal{Z}({L}{u}) \cap \mathcal{O}_{{p}}$ is a real-analytic curve with ${p}$ as an interior point and $\mathcal{Z}({L}{u}) \cap \mathcal{O}_{{p}} \cap \partial {T} = \{{p}\}$. 
\end{enumerate}
\end{lma}

\begin{proof}
Consider a Cartesian coordinate system $({y}_{1}, {y}_{2})$ with the origin at ${p}$ such that the ${y}_{1}$-axis is the tangential direction $\boldsymbol{\tau}$ and the positive ${y}_{2}$-axis is the inward normal direction $\boldsymbol{n}$ at ${p}$. Let $\tilde{{u}}$ be the extension of ${u}$ to the interior of the kite ${K}_{{e}}$ obtained by reflecting across ${e}$. 
Since $\partial_{{y}_{1}}{u}$ satisfies the Neumann boundary condition and $\partial_{{y}_{1}}{u}$ cannot vanish identically in ${T}$, from \cite{HW53} that, in a small neighborhood of $y = 0$, 
\begin{equation*}
\partial_{{y}_{1}}\tilde{u}(y) = {a}_{{n} + 1} \operatorname{Re}\big(({y}_{1} + {y}_{2}\sqrt{ - 1})^{n}\big) + O(|y|^{{n} + 1})
\end{equation*}
for some ${n} \in \mathbb{N}^{ + }$ and ${a}_{{n} + 1} \in \R \setminus \{0\}$. Here $O(|y|^{l})$ means that $O(|y|^{l})/|y|^{l}$ is bounded as $|y| \to 0$. Moreover, in a small neighborhood of $y = 0$, 
\begin{align*}
\tilde{u}(y) & = {b}_{0} + \frac{1}{2}{b}_{2}{y}_{2}^{2} + \frac{1}{{n} + 1}{a}_{{n} + 1}{y}_{1}^{{n} + 1} + {y}_{2}^{2} \cdot O(|y|) + O(|y|^{{n} + 2}), 
\\ 
{L}\tilde{u}(y) & = {c}_{1}{a}_{{n} + 1}{y}_{1}^{n} + {c}_{2}{b}_{2}{y}_{2} + {y}_{2} \cdot O(|y|) + O(|y|^{{n} + 1}), 
\end{align*}
where ${b}_{0} = {u} ({p})$, ${b}_{2} = \partial_{\boldsymbol{n}\boldsymbol{n}}{u}({p})$ and ${L} = {c}_{1}\partial_{{y}_{1}} + {c}_{2}\partial_{{y}_{2}}$. 
Since ${c}_{2}{c}_{1} \neq 0$ and ${b}_{2} \neq 0$, the implicit function theorem implies that there exists a real-analytic function $\varphi = \varphi({y}_{1})$ defined near ${y}_{1} = 0$ and a small neighborhood $\mathcal{O}_{0}$ of $y = 0$ such that
\begin{gather*}
\{y \in \mathcal{O}_{0}: \, \partial_{{y}_{1}}\tilde{u}(y) = 0 \} =
\{y \in \mathcal{O}_{0}: \, {y}_{2} = \varphi({y}_{1}) \}, 
\\
\varphi({y}_{1}) = - {c}_{3}{y}_{1}^{n} (1 + O(|{y}_{1}|)).
\end{gather*}
Hence, the nodal line of ${L}\tilde{u}$ is a simple real-analytic curve passing through $y = 0$. In order to find the nodal line $\mathcal{Z}({L}{u})$, it suffices to determine the condition for the positivity of $\varphi$. We obtain that $\varphi({y}_{1}) > 0$ in a sufficiently small neighborhood of ${y}_1 = 0$ if and only if one of the following holds: (i) ${n}$ is odd and ${c}_{3}{y}_{1} < 0$; (ii) ${n}$ is even and ${c}_{3} < 0$, with ${y}_{1} \neq 0$. Hence the proof is complete. 
\end{proof}

\begin{lma} \label{lma35}
Let ${T}$ be an acute triangle ${z}_{1}{z}_{2}{z}_{3}$ as before. Suppose that ${u}$ satisfies the monotonicity condition \eqref{CGY0301}. 
Then the following statements hold (see \autoref{fig3A}): 
\begin{enumerate}[label = {\rm(\roman*)}]
\item
${u}$ has at most one non-vertex critical point on the lower side $\overline{{z}_{1}{z}_{2}}$ and at most one non-vertex critical point on the union of the remaining sides $\overline{{z}_{3}{z}_{1}} \cup \overline{{z}_{3}{z}_{2}}$. 
\item
If a non-vertex critical point exists on $\overline{{z}_{3}{z}_{1}} \cup \overline{{z}_{3}{z}_{2}}$, then a non-vertex critical point exists on $\overline{{z}_{1}{z}_{2}}$. 
\item 
A non-vertex critical point exists if and only if ${u}({z}_{1}) \cdot {u}({z}_{2}) \neq 0$, that is, ${u}({z}_{1}) < 0$ and ${u}({z}_{2}) < 0$; and no non-vertex critical point exists if and only if ${u}({z}_{1}) \cdot {u}({z}_{2}) = 0$. 
\item
At the unique non-vertex critical point on $\overline{{z}_{1}{z}_{2}}$ (if it exists), the second tangential derivative is negative. Moreover, at the unique non-vertex critical point on $\overline{{z}_{3}{z}_{1}} \cup \overline{{z}_{3}{z}_{2}}$ (if it exists), the third tangential derivative along the corresponding side is nonzero. 
\end{enumerate}
\end{lma}

\begin{proof} 
The argument is based on the geometry of the nodal lines and nodal domains of the directional derivatives of the eigenfunction ${u}$. Recall that a nodal line is a union of immersed $C^{1}$ loops and properly immersed $C^{1}$ arcs; see \cite{Che76} and Section 2 of \cite{JM20}. 
Define
\begin{equation}
\begin{aligned}
\Gamma^{S} & = \overline{\{{x} \in {T}: \, \nabla {u}({x}) \cdot \boldsymbol{\tau}_{S} = 0\}}, \\
\Gamma^{M} & = \overline{\{{x} \in {T}: \, \nabla {u}({x}) \cdot \boldsymbol{\tau}_{M} = 0\}}, \\
\Gamma^{L} & = \overline{\{{x} \in {T}: \, \nabla {u}({x}) \cdot \boldsymbol{\tau}_{L} = 0\}}, 
\end{aligned}
\end{equation}
where $\boldsymbol{\tau}_{S}$, $\boldsymbol{\tau}_{M}$ and $\boldsymbol{\tau}_{L}$ denote the unit tangent directions along the sides $\overline{{z}_{1}{z}_{2}}$, $\overline{{z}_{3}{z}_{1}}$, and $\overline{{z}_{2}{z}_{3}}$, respectively; see also below \eqref{CGY0106angle}. 

\textbf{Claim 1}. 
None of $\Gamma^{S} \cup \overline{{z}_{1}{z}_{2}}$, $\Gamma^{M} \cup \overline{{z}_{1}{z}_{3}}$, or $\Gamma^{L} \cup \overline{{z}_{2}{z}_{3}}$ contains a loop. Without loss of generality, assume that $\Gamma^{S} \cup \overline{{z}_{1}{z}_{2}}$ does contain a loop. Then there exists a nodal domain ${D}$ of $\nabla {u} \cdot \boldsymbol{\tau}_{S}$ with $\partial {D} \subset \Gamma^{S} \cup \overline{{z}_{1}{z}_{2}}$. The function $\nabla {u} \cdot \boldsymbol{\tau}_{S}$ then satisfies
\begin{equation*}\begin{cases}
\Delta (\nabla {u} \cdot \boldsymbol{\tau}_{S}) + \mu (\nabla {u} \cdot \boldsymbol{\tau}_{S}) = 0 & \text{ in } {D}, \\
\nabla {u} \cdot \boldsymbol{\tau}_{S} = 0 & \text{ on } \partial {D} \setminus \overline{{z}_{1}{z}_{2}}, \\
\partial_{\nu}(\nabla {u} \cdot \boldsymbol{\tau}_{S}) = 0 & \text{ on } \partial {D} \cap \overline{{z}_{1}{z}_{2}}. 
\end{cases}\end{equation*}
Hence the first mixed eigenvalue $\lambda_{1}({D}, \partial {D} \setminus \overline{{z}_{1}{z}_{2}})$ equals $\mu$. By the variational characterization, one has $\lambda_{1}({D}, \partial {D} \setminus \overline{{z}_{1}{z}_{2}}) > \lambda_{1}({T}, \partial {T} \setminus \overline{{z}_{1}{z}_{2}})$. 
On the other hand, by the inequality between Neumann and mixed eigenvalues (see \autoref{lma25LR}), we obtain $\lambda_{1}({T}, \partial {T} \setminus \overline{{z}_{1}{z}_{2}}) \geq \mu_{2}({T}) = \mu$, which leads to a contradiction. 
Therefore, $\Gamma^{S} \cup \overline{{z}_{1}{z}_{2}}$ contains no loop, and Claim 1 follows. 

As a consequence of Claim 1, $\Gamma^{S}$ is a finite union of properly immersed $C^{1}$ arcs $\Gamma^{S}_{{j}}$, with ${j} = 1, 2, \ldots, {n}$, each having both endpoints on $\partial {T}$. Moreover, the side $\overline{{z}_{1}{z}_{2}}$ can contain at most one endpoint of any $\Gamma^{S}_{{j}}$. The analogous statements hold for $\Gamma^{M}$ and $\Gamma^{L}$. 

\textbf{Claim 2}. 
Local structure near the vertices can be described as follows. 
\begin{enumerate}
\item \label{CGY0309ita}
If ${u}({z}_{3}) = 0$, then there exists a small neighborhood $\mathcal{O}_{{z}_{3}}$ of ${z}_{3}$ such that $\nabla {u} \cdot \boldsymbol{\tau}_{S} \neq 0$ in $(\overline{{T}} \setminus \{{z}_{3}\}) \cap \mathcal{O}_{{z}_{3}}$. 
\item \label{CGY0309itb}
If ${u}({z}_{3}) \neq 0$, then there exists a small neighborhood $\mathcal{O}_{{z}_{3}}$ of ${z}_{3}$ such that $\Gamma^{S} \cap \mathcal{O}_{{z}_{3}}$ is a simple real-analytic curve passing through ${z}_{3}$ and transversal to the two sides of ${T}$ adjacent to ${z}_{3}$. Moreover, $\Gamma^{S}$ separates $\mathcal{O}_{{z}_{3}} \cap {T}$ into exactly two connected components. See red curves in \autoref{fig3A} below. 
\item \label{CGY0309itc}
If ${u}({z}_{2}) = 0$, then there exists a small neighborhood $\mathcal{O}_{{z}_{2}}$ of ${z}_{2}$ such that $\Gamma^{S} \cap \mathcal{O}_{{z}_{2}}$ is a simple real-analytic curve passing through ${z}_{2}$ and transversal to the two sides of ${T}$ adjacent to ${z}_{2}$. Moreover, $\Gamma^{S}$ separates $\mathcal{O}_{{z}_{2}} \cap {T}$ into exactly two connected components. 
\item \label{CGY0309itd}
If ${u}({z}_{2}) \neq 0$, then there exists a small neighborhood $\mathcal{O}_{{z}_{2}}$ of ${z}_{2}$ such that $\nabla {u} \cdot \boldsymbol{\tau}_{S} \neq 0$ in $(\overline{{T}} \setminus \{{z}_{2}\}) \cap \mathcal{O}_{{z}_{2}}$. 
\end{enumerate}
The same conclusions hold for $\Gamma^{M}$ and $\Gamma^{L}$; the arguments are identical and therefore omitted. Each item follows from the local Fourier--Bessel expansion of the eigenfunction near a vertex; see \autoref{lma26bJM}. 

\textbf{Claim 3}. 
The nodal sets $\Gamma^{S}$, $\Gamma^{M}$ and $\Gamma^{L}$ enjoy the following properties: 
\begin{itemize}
\item
$\Gamma^{L} \cap \operatorname{Int}(\overline{{z}_{3}{z}_{1}}) = \emptyset$ and $\Gamma^{M} \cap \operatorname{Int}(\overline{{z}_{3}{z}_{2}}) = \emptyset$.
\item 
For any ${q} \in \operatorname{crit}_{\mathrm{nv}}({u}) \cap (\overline{{z}_{3}{z}_{1}} \cup \overline{{z}_{3}{z}_{2}})$, there exists a small neighborhood $\mathcal{O}_{{q}}$ of ${q}$ such that $\Gamma^{S} \cap \mathcal{O}_{{q}}$ is a simple real-analytic curve (with ${q}$ as an interior point) contained in ${T} \cup \{{q}\}$.
\end{itemize}
Indeed, every point in $(\Gamma^{S} \cup \Gamma^{L}) \cap \operatorname{Int}(\overline{{z}_{3}{z}_{1}})$ is a critical point of ${u}$. Now let ${q} \in \operatorname{crit}_{\mathrm{nv}}({u}) \cap \overline{{z}_{3}{z}_{1}}$. 
Note that
\begin{equation*}\begin{aligned}
&\boldsymbol{\tau}_{S} = {c}_{1}^{S, M} \boldsymbol{\tau}_{M} + {c}_{2}^{S, M} \boldsymbol{n}_{M}, \quad {c}_{1}^{S, M} = - \cos\alpha_{1} < 0, \quad {c}_{2}^{S, M} = \sin\alpha_{1} > 0, 
\\
&\boldsymbol{\tau}_{L} = {c}_{1}^{L, M} \boldsymbol{\tau}_{M} + {c}_{2}^{L, M} \boldsymbol{n}_{M}, \quad
{c}_{1}^{L, M} = - \cos\alpha_{3} < 0, \quad {c}_{2}^{L, M} = - \sin\alpha_{3} < 0, 
\end{aligned}
\end{equation*}
where $\alpha_{j}$ is defined in \eqref{CGY0106angle}. 
Differentiating the Neumann boundary condition of ${u}$ along the boundary, we obtain $\partial_{\boldsymbol{\tau}_{M}\boldsymbol{n}_{M}}{u} = 0$ on $\overline{{z}_{3}{z}_{1}}$. From the monotonicity condition \eqref{CGY0301} and the Neumann boundary condition for ${u}$, we deduce that
\begin{equation} \label{CGY0310} 
\partial_{\boldsymbol{\tau}_{M}}{u} \leq 0 \text{ on } \overline{{z}_{3}{z}_{1}} \text{ and } \partial_{\boldsymbol{\tau}_{L}}{u} \geq 0 \text{ on } \overline{{z}_{2}{z}_{3}}. 
\end{equation}
Combining this with the fact that $\partial_{\boldsymbol{\tau}_{M}}{u}({q}) = 0$, we conclude that $\partial_{\boldsymbol{\tau}_{M}\boldsymbol{\tau}_{M}}{u}({q}) = 0$. Next, applying the Hopf lemma to the positive function $\partial_{{x}_{2}} {u}$ at ${q}$ yields $\partial_{{x}_{2}{x}_{2}} {u}({q}) < 0$. Therefore, $\partial_{\boldsymbol{n}_{M}\boldsymbol{n}_{M}}{u}({q}) < 0$.
Again, by the monotonicity assumption \eqref{CGY0301} and \autoref{lma34}, we have 
\begin{equation*}
(\partial_{\boldsymbol{\tau}_{M}})^{{n} + 1} {u} < 0 = (\partial_{\boldsymbol{\tau}_{M}})^{n} {u} = \cdots = (\partial_{\boldsymbol{\tau}_{M}})^{1} {u} \text{ and } \partial_{\boldsymbol{n}_{M}\boldsymbol{n}_{M}} {u} < 0 \text{ at the point } {q}
\end{equation*}
for some even integer ${n} \geq 2$. 
Now, applying \autoref{lma34}, we infer that, in a neighborhood $\mathcal{O}_{{q}}$ of ${q}$, the nodal line $\Gamma^{S}$ is a real-analytic curve while $\Gamma^{L}$ is empty.
The same analysis applies to a non-vertex critical point on the right side $\overline{{z}_{2}{z}_{3}}$. This completes Claim 3.

\begin{figure}[htp]\centering 
\begin{tikzpicture}[scale = 2.8]
\pgfmathsetmacro\LenBase{1.00}; 
\pgfmathsetmacro\AngleL{79}; 
\pgfmathsetmacro\AngleM{59}; 
\pgfmathsetmacro\AngleS{180-\AngleL-\AngleM}; 
\pgfmathsetmacro\LenLeft{\LenBase*sin(\AngleM)/sin(\AngleS)}; 
\pgfmathsetmacro\LenRight{\LenBase*sin(\AngleL)/sin(\AngleS)}; 
\pgfmathsetmacro\xSS{\LenLeft*cos(\AngleL)}; 
\pgfmathsetmacro\ySS{\LenLeft*sin(\AngleL)}; 
\pgfmathsetmacro\xLL{0}; 
\pgfmathsetmacro\yLL{0}; 
\pgfmathsetmacro\xMM{\LenBase}; 
\pgfmathsetmacro\yMM{0}; 
\pgfmathsetmacro\xPP{\LenBase/5};  \pgfmathsetmacro\yPP{0}; 
\pgfmathsetmacro\xQQ{cos(\AngleL)*\LenLeft*0.33}; 
\pgfmathsetmacro\yQQ{sin(\AngleL)*\LenLeft*0.33}; 
\pgfmathsetmacro\xIa{\xPP+cos(90)*\LenBase/5}; 
\pgfmathsetmacro\yIa{\yPP+sin(90)*\LenBase/5}; 
\pgfmathsetmacro\xIb{\xQQ+cos(\AngleL-180)*\LenLeft/5}; 
\pgfmathsetmacro\yIb{\yQQ+sin(\AngleL-180)*\LenLeft/5}; 
\pgfmathsetmacro\xIc{\xQQ+cos(\AngleL)*\LenLeft/6}; 
\pgfmathsetmacro\yIc{\yQQ+sin(\AngleL)*\LenLeft/6}; 
\pgfmathsetmacro\xId{\xSS+cos(-88)*\LenLeft/2}; 
\pgfmathsetmacro\yId{\ySS+sin(-88)*\LenLeft/2}; 
\pgfmathsetmacro\xHa{\xPP+cos(55)*\LenBase/5}; 
\pgfmathsetmacro\yHa{\yPP+sin(55)*\LenBase/5}; 
\pgfmathsetmacro\xHb{\xQQ+cos(\AngleL-135)*\LenLeft/4}; 
\pgfmathsetmacro\yHb{\yQQ+sin(\AngleL-135)*\LenLeft/4}; 
\pgfmathsetmacro\xHc{\xQQ+cos(\AngleL-45)*\LenLeft/2}; 
\pgfmathsetmacro\yHc{\yQQ+sin(\AngleL-45)*\LenLeft/2}; 
\pgfmathsetmacro\xHd{\xMM+cos(\AngleL+90)*\LenLeft/3}; 
\pgfmathsetmacro\yHd{\yMM+sin(\AngleL+90)*\LenLeft/3}; 
\draw[green] (\xQQ, \yQQ) .. controls (\xHc, \yHc) and (\xHd, \yHd) .. (\xMM, \yMM); 
\draw[green] (\xQQ, \yQQ) .. controls (\xHb, \yHb) and (\xHa, \yHa) .. (\xPP, \yPP); 
\draw[red, thick] (\xQQ, \yQQ) .. controls (\xIb, \yIb) and (\xIa, \yIa) .. (\xPP, \yPP); 
\draw[red, thick] (\xQQ, \yQQ) .. controls (\xIc, \yIc) and (\xId, \yId) .. (\xSS, \ySS); 
\draw (\xQQ, \yQQ) node[left]{\small ${q}$}; 
\draw (\xPP, \yPP) node[below]{\small ${p}$}; 
\draw[black] (\xSS, \ySS) node[below = 2pt, right]{\small ${z}_{3}$} -- (\xLL, \yLL) node[below]{\small ${z}_{1}$} -- (\xMM, \yMM) node[below]{\small ${z}_{2}$} -- cycle; 
\draw[->] ({(\xSS+\xMM)*0.65}, \ySS*0.7) -- ++ ({0} : 0.3*\LenBase) node [above] {$\boldsymbol{\tau}_{S}$}; 
\draw[->] ({(\xSS+\xMM)*0.65}, \ySS*0.7) -- ++ ({180+\AngleL} : 0.3*\LenBase) node [right] {$\boldsymbol{\tau}_{M}$}; 
\draw[->] ({(\xSS+\xMM)*0.65}, \ySS*0.7) -- ++ ({180-\AngleM} : 0.3*\LenBase) node [right] {$\boldsymbol{\tau}_{L}$}; 
\end{tikzpicture} \hspace{2ex}
\begin{tikzpicture}[scale = 2.8]
\pgfmathsetmacro\AngleL{79}; 
\pgfmathsetmacro\AngleM{59}; 
\pgfmathsetmacro\AngleS{180-\AngleL-\AngleM}; 
\pgfmathsetmacro\LenBase{1.00}; 
\pgfmathsetmacro\LenLeft{\LenBase*sin(\AngleM)/sin(\AngleS)}; 
\pgfmathsetmacro\LenRight{\LenBase*sin(\AngleL)/sin(\AngleS)}; 
\pgfmathsetmacro\xSS{\LenLeft*cos(\AngleL)}; 
\pgfmathsetmacro\ySS{\LenLeft*sin(\AngleL)}; 
\pgfmathsetmacro\xLL{0}; 
\pgfmathsetmacro\yLL{0}; 
\pgfmathsetmacro\xMM{\LenBase}; 
\pgfmathsetmacro\yMM{0}; 
\pgfmathsetmacro\xPP{\LenBase/5}; 
\pgfmathsetmacro\yPP{0}; 
\pgfmathsetmacro\xIc{\xPP+cos(90)*\LenLeft/2}; 
\pgfmathsetmacro\yIc{\yPP+sin(90)*\LenLeft/2}; 
\pgfmathsetmacro\xId{\LenLeft*cos(\AngleL)+cos(-90)*\LenLeft/3}; 
\pgfmathsetmacro\yId{\LenLeft*sin(\AngleL)+sin(-90)*\LenLeft/3}; 
\pgfmathsetmacro\xHc{\xPP+cos(70)*\LenBase}; 
\pgfmathsetmacro\yHc{\yPP+sin(70)*\LenBase}; 
\pgfmathsetmacro\xHd{\xMM+cos(\AngleL+90)*\LenLeft/3}; 
\pgfmathsetmacro\yHd{\yMM+sin(\AngleL+90)*\LenLeft/3}; 
\draw[green] (\xPP, \yPP) .. controls (\xHc, \yHc) and (\xHd, \yHd) .. (\xMM, \yMM); 
\draw[red, thick] (\xPP, \yPP) .. controls (\xIc, \yIc) and (\xId, \yId) .. (\xSS, \ySS); 
\draw[black] (\xSS, \ySS) node[below = 2pt, right]{\small ${z}_{3}$} -- (\xLL, \yLL) node[below]{\small ${z}_{1}$} -- (\xMM, \yMM) node[below]{\small ${z}_{2}$} -- cycle; 
\draw (\xPP, \yPP) node[below]{\small ${p}$}; 
\end{tikzpicture} \hspace{2ex}
\begin{tikzpicture}[scale = 2.8]
\pgfmathsetmacro\AngleL{79}; 
\pgfmathsetmacro\AngleM{59}; 
\pgfmathsetmacro\AngleS{180-\AngleL-\AngleM}; 
\pgfmathsetmacro\LenBase{1.00}; 
\pgfmathsetmacro\LenLeft{\LenBase*sin(\AngleM)/sin(\AngleS)}; 
\pgfmathsetmacro\LenRight{\LenBase*sin(\AngleL)/sin(\AngleS)}; 
\pgfmathsetmacro\xSS{\LenLeft*cos(\AngleL)}; 
\pgfmathsetmacro\ySS{\LenLeft*sin(\AngleL)}; 
\pgfmathsetmacro\xLL{0}; 
\pgfmathsetmacro\yLL{0}; 
\pgfmathsetmacro\xMM{\LenBase}; 
\pgfmathsetmacro\yMM{0}; 
\pgfmathsetmacro\xPP{0}; 
\pgfmathsetmacro\yPP{0}; 
\pgfmathsetmacro\xIc{\xPP+cos(\AngleL*(90+\AngleL)/(180+\AngleL))*\LenLeft/3}; 
\pgfmathsetmacro\yIc{\yPP+sin(\AngleL*(90+\AngleL)/(180+\AngleL))*\LenLeft/3}; 
\pgfmathsetmacro\xId{\LenLeft*cos(\AngleL)+cos(-90)*\LenLeft/3}; 
\pgfmathsetmacro\yId{\LenLeft*sin(\AngleL)+sin(-90)*\LenLeft/3}; 
\pgfmathsetmacro\xHc{\xPP+cos(90*\AngleL/(180+\AngleL))*\LenBase/2}; 
\pgfmathsetmacro\yHc{\yPP+sin(90*\AngleL/(180+\AngleL))*\LenBase/2}; 
\pgfmathsetmacro\xHd{\LenBase+cos(\AngleL+90)*\LenBase/2}; 
\pgfmathsetmacro\yHd{0+sin(\AngleL+90)*\LenBase/2}; 
\draw[green] (\xPP, \yPP) .. controls (\xHc, \yHc) and (\xHd, \yHd) .. (\xMM, \yMM); 
\draw[red, thick] (\xPP, \yPP) .. controls (\xIc, \yIc) and (\xId, \yId) .. (\xSS, \ySS); 
\draw[black] (\xSS, \ySS) node[below = 2pt, right]{\small ${z}_{3}$} -- (\xLL, \yLL) node[below]{\small ${z}_{1}$} -- (\xMM, \yMM) node[below]{\small ${z}_{2}$} -- cycle; 
\end{tikzpicture}
\vspace*{-2ex}

\caption{
The nodal line of $\color{red}\nabla {u}\cdot \mathbf{\tau}_{S}$ (in red) and nodal line of $\color{green}\nabla {u}\cdot \mathbf{\tau}_{M}$ (in green)
}
\label{fig3A}
\end{figure}

\textbf{Claim 4}. 
The set $\Gamma^{S}$ is a simple real-analytic arc joining the vertex ${z}_{3}$ and a point ${p} \in \overline{{z}_{1}{z}_{2}}$. Moreover, ${u}$ has at most one non-vertex critical point on $\overline{{z}_{1}{z}_{2}}$, at this critical point the second tangential derivative of ${u}$ is negative, and  
\begin{equation} \label{CGY0312}
\operatorname{crit}_{\mathrm{nv}}({u}) \cap \operatorname{Int}(\overline{{z}_{1}{z}_{2}}) \neq \emptyset \Longleftrightarrow {u}({z}_{1}) \cdot {u}({z}_{2}) > 0, \text{ i.e., } {u}({z}_{1}) < 0, {u}({z}_{2}) < 0.
\end{equation} 
In fact, by \autoref{lma26bJM} and \eqref{CGY0310}, we have
\begin{equation*}
{u}({z}_{1}) \leq 0 \text{ and } {u}({z}_{2}) \leq 0, 
\end{equation*}
while ${u}({z}_{3}) = \sup_{{T}} {u} > 0$. Claim 3 and Claim 1 imply that every arc of $\Gamma^{S}$ has its endpoints in $\overline{{z}_{1}{z}_{2}} \cup \{{z}_{3}\}$. Property \eqref{CGY0309itb} of Claim 2 ensures that exactly one arc of $\Gamma^{S}$ emanates from ${z}_{3}$. Again by Claim 1, $\Gamma^{S}$ is a properly immersed ${C}^{1}$ arc connecting ${z}_{3}$ to a point ${p}$ on $\overline{{z}_{1}{z}_{2}}$; see \autoref{fig3A}. By \autoref{lma33}, there exists at least one nodal line of $\nabla {u} \cdot \boldsymbol{\tau}_{S}$ emanating from each non-vertex critical point on the lower side $\overline{{z}_{1}{z}_{2}}$. Hence ${u}$ has no non-vertex critical point on $\overline{{z}_{1}{z}_{2}} \setminus \{{p}\}$. Using properties \eqref{CGY0309itc} and \eqref{CGY0309itd} of Claim 2 together with \autoref{lma33}, there are two possibilities: 
\begin{itemize}
\item
${p} \in \operatorname{Int}(\overline{{z}_{1}{z}_{2}})$, ${u}({z}_{1}) < 0$, ${u}({z}_{2}) < 0$, $\operatorname{crit}_{\mathrm{nv}}({u}) \cap \operatorname{Int}(\overline{{z}_{1}{z}_{2}}) = \{{p}\}$, and the second tangential derivative of ${u}$ at ${p}$ is nonzero (hence negative); 
\item
${p} \in \{{z}_{1}, {z}_{2}\}$, ${u}({p}) = 0$ and $\operatorname{crit}_{\mathrm{nv}}({u}) \cap \operatorname{Int}(\overline{{z}_{1}{z}_{2}}) = \emptyset$. 
\end{itemize}

\textbf{Claim 5}. 
The function ${u}$ has at most one non-vertex critical point on each of the two oblique sides $\overline{{z}_{3}{z}_{1}}$ and $\overline{{z}_{3}{z}_{2}}$, and the third-order tangential derivative of ${u}$ at any such point is nonzero.  Moreover, if $\operatorname{crit}_{\mathrm{nv}}({u}) \cap (\overline{{z}_{3}{z}_{1}} \cup \overline{{z}_{3}{z}_{2}}) \neq \emptyset$, then necessarily $\operatorname{crit}_{\mathrm{nv}}({u}) \cap \overline{{z}_{1}{z}_{2}} \neq \emptyset$. 

To prove this, we analyze the nodal line $\Gamma^{M}$ together with $\operatorname{crit}_{\mathrm{nv}}({u}) \cap \overline{{z}_{3}{z}_{1}}$; the case of $\Gamma^{L}$ and $\operatorname{crit}_{\mathrm{nv}}({u}) \cap \overline{{z}_{3}{z}_{2}}$ is analogous.
By Claims 1, 3, and 4, the endpoints of every immersed ${C}^{1}$ arc of $\Gamma^{M}$ lie in the set $\{{z}_{2}, {p}\} \cup \overline{{z}_{3}{z}_{1}}$, where ${p}$ is the point given in Claim 4.
By the local behavior near the vertices (see \autoref{lma26bJM}), three mutually exclusive cases occur: 

Case 1: ${p} \in \operatorname{Int}(\overline{{z}_{1}{z}_{2}})$. 
Then ${p}$ is a nondegenerate critical point, and ${u}({z}_{1}) < 0$ and ${u}({z}_{2}) < 0$. 
By \autoref{lma34}, there are exactly two local branches of $\Gamma^{M}$ emanating from $\partial {T} \setminus \operatorname{Int}(\overline{{z}_{3}{z}_{1}})$, namely, one branch from each of the two points ${p}$ and ${z}_{2}$. Note that these two branches may meet and merge into a single global ${C}^{1}$ arc connecting ${p}$ and ${z}_{2}$; see the middle panel in \autoref{fig3A}.

Case 2: ${p} = {z}_{1}$. 
Then ${u}({z}_{1}) = 0$ and ${u}({z}_{2}) < 0$, and there are exactly two local branches of $\Gamma^{M}$ emanating from $\partial {T} \setminus \operatorname{Int}(\overline{{z}_{3}{z}_{1}})$, namely, one branch from each of the two points ${z}_{1}$ and ${z}_{2}$. Note that these two branches may meet and merge into a single global $C^{1}$ arc connecting ${z}_{1}$ and ${z}_{2}$; see the right panel in \autoref{fig3A}.

Case 3: ${p} = {z}_{2}$.
Then ${u}({z}_{1}) < 0$ and ${u}({z}_{2}) = 0$, and $\Gamma^{M}$ has no branch that meets $\partial {T} \setminus \operatorname{Int}(\overline{{z}_{3}{z}_{1}})$.

Combining the three cases, we conclude that $\Gamma^{M}$ contains at most two local branches meeting $\partial {T} \setminus \overline{{z}_{3}{z}_{1}}$. Together with Claim 1, it follows that $\Gamma^{M}$ consists of at most two
properly immersed $C^{1}$ arcs, which implies that at most two local branches emanate from the side $\overline{{z}_{3}{z}_{1}}$. 
\autoref{lma33} and the monotonicity condition~\eqref{CGY0301} guarantee that at least two branches of $\Gamma^{M}$ emanate
from any non-vertex critical point on $\overline{{z}_{3}{z}_{1}}$. Consequently, $\operatorname{{crit}}_{\mathrm{nv}}({u}) \cap \overline{{z}_{3}{z}_{1}}$ has at most one element. 

Assume now that a non-vertex critical point actually exists on $\overline{{z}_{3}{z}_{1}}$, denote it by ${q}$. Then $\Gamma^{M}$ is the union of two real-analytic arcs, one joining ${q}$ to ${z}_{2}$, the other joining
${q}$ to ${p}$. Additionally, \autoref{lma33} implies that the third-order tangential derivative of ${u}$ at ${q}$ is nonzero. If ${p} \in \{{z}_{1}, {z}_{2}\}$, then those two arcs together with the segment $\overline{{z}_{1}{z}_{2}}$ would form a loop, contradicting
Claim 1. Therefore ${p}$ must lie in $\operatorname{Int}(\overline{{z}_{1}{z}_{2}})$; see the left panel in \autoref{fig3A}. By definition of ${p}$ in Claim 4, this means that ${p} \in \operatorname{crit}_{\mathrm{nv}}({u}) \cap \overline{{z}_{1}{z}_{2}}$. This completes the proof of Claim 5. 
By \eqref{CGY0312} and Claim~5, we conclude parts~(ii) and~(iii) of \autoref{lma35}. Part~(iv) of \autoref{lma35} follows from Claims~4 and~5. 

\textbf{Claim 6}. 
$\operatorname{Int}(\overline{{z}_{3}{z}_{1}}) \cup \operatorname{Int}(\overline{{z}_{3}{z}_{2}})$ contains at most one critical point.
Suppose, for contradiction, that there exist two distinct non-vertex critical points ${q}_{1} \in \operatorname{crit}_{\mathrm{nv}}({u}) \cap \operatorname{Int}(\overline{{z}_{3}{z}_{2}})$ and ${q}_{2} \in \operatorname{crit}_{\mathrm{nv}}({u}) \cap \operatorname{Int}(\overline{{z}_{3}{z}_{1}})$.
By Claim 5, there exists a real-analytic arc ${\ell}_{1}$ of $\Gamma^{L}$ connecting ${q}_{1}$ to ${z}_{1}$, and a real-analytic arc ${\ell}_{2}$ of $\Gamma^{M}$ connecting ${q}_{2}$ to ${z}_{2}$.
It is clear that the arc ${\ell}_{2}$ separates ${T}$ into two subdomains, with ${q}_{1}$ and ${z}_{1}$ lying in different subdomains and on different sides of ${\ell}_{2}$. By Claims 3 and 2, $\Gamma^{L}$, and hence ${\ell}_{1}$, does not meet ${q}_{2}$ or ${z}_{2}$. Therefore, the arc ${\ell}_{1}$ must intersect the arc ${\ell}_{2}$ at some interior point $\bar{x}$ of ${T}$. At this point $\bar{x}$ we have $\nabla {u} \cdot \boldsymbol{\tau}_{L} = \nabla {u} \cdot \boldsymbol{\tau}_{M} = 0$, which implies $|\nabla {u}(\bar{x})| = 0$, contradicting the positivity of $\partial_{{x}_{2}} {u}$ in ${T}$.
Hence our assumption is false, and $\operatorname{Int}(\overline{{z}_{3}{z}_{1}}) \cup \operatorname{Int}(\overline{{z}_{3}{z}_{2}})$ contains at most one critical point. This completes Claim~5, and therefore part~(i) of \autoref{lma35} follows from Claims~4 and~6. 

This completes the proof of \autoref{lma35}.
\end{proof}

For acute triangles whose smallest angle is at most $\pi/6$, as a direct consequence of \cite{Siu15} and the local behavior near the vertices (\autoref{lma26bJM}), one can easily deduce that the monotonicity property \eqref{CGY0301} and hence \autoref{thm12} hold, except for the statement regarding the exact location of the global extrema.

\begin{cor}
Let ${u}$ be a second Neumann eigenfunction of an acute triangle ${T} = \triangle {z}_{1}{z}_{2}{z}_{3}$ with $\angle {z}_{2}{z}_{3}{z}_{1} \leq \pi/6$. Then ${u}$ must satisfy the following properties: 
\begin{enumerate} [label = \rm(\roman*)]
\item
$\nabla {u} \cdot \boldsymbol{n}_{S} \neq 0$ in ${T}$; 
\item
${u}({z}_{3}) \cdot {u}({z}_{1}) < 0$ and ${u}({z}_{3}) \cdot {u}({z}_{2}) < 0$; 
\item
${u}$ has exactly one non-vertex critical point, which lies in the interior of the shortest side $\overline{{z}_{1}{z}_{2}}$ and is a saddle point; 
\item
the second Neumann eigenvalue of ${T}$ is simple. 
\end{enumerate}
\end{cor}

\begin{proof} 
The assumptions ensure that $\overline{{z}_{1}{z}_{2}}$ is the shortest side of ${T}$. According to \cite[Lemmas 1, 2, and 3]{Siu15} and the condition $\angle {z}_{2}{z}_{3}{z}_{1} \leq \pi/6$, we conclude that ${u}$ has no critical points on $\operatorname{Int}(\overline{{z}_{3}{z}_{1}}) \cup \operatorname{Int}(\overline{{z}_{3}{z}_{2}})$. 

We claim that the nodal line $\mathcal{Z}({u})$ intersects both $\operatorname{Int}(\overline{{z}_{3}{z}_{1}})$ and $\operatorname{Int}(\overline{{z}_{3}{z}_{2}})$. Consequently, ${u}({z}_{1}) \cdot {u}({z}_{2}) > 0$ and ${u}({z}_{3}) \cdot {u}({z}_{1}) < 0$. In fact, we assume by contradiction that $\mathcal{Z}({u})$ does not intersect $\operatorname{Int}(\overline{{z}_{3}{z}_{1}})$. Reflect ${u}$ across the left side $\overline{{z}_{3}{z}_{1}}$ to obtain a new function $\tilde{u}$ on the kite ${K}$ formed by this reflection. 
By \cite[Lemma 1]{Siu15}, $\tilde{{u}}$ must be a second Neumann eigenfunction on ${K}$. However, the symmetry of $\tilde{{u}}$ and the assumption that $\mathcal{Z}(\tilde{u}) \cap \operatorname{Int}(\overline{{z}_{3}{z}_{1}}) = \emptyset$ imply that $\tilde{{u}}$ has exactly three nodal domains. This contradicts the Courant nodal domain theorem, which states that the second eigenfunction can have exactly two nodal domains. Therefore, $\mathcal{Z}({u})$ must intersect both $\operatorname{Int}(\overline{{z}_{3}{z}_{1}})$ and $\operatorname{Int}(\overline{{z}_{3}{z}_{2}})$. 

Since ${u}({z}_{3}) \neq 0$, the local expansion of ${u}$ near ${z}_{3}$ (see \autoref{lma26bJM}) indicates that ${z}_{3}$ is a strict extremum point and $\partial_{{x}_{2}}{u}$ has a fixed sign on both $\operatorname{Int}(\overline{{z}_{3}{z}_{1}})$ and $\operatorname{Int}(\overline{{z}_{3}{z}_{2}})$. Without loss of generality, assume $\partial_{{x}_{2}}{u} > 0$ on these two sides. Hence, 
\begin{equation*} 
\partial_{{x}_{2}}{u} \geq, \not\equiv 0 \text{ on } \partial {T}. 
\end{equation*} 
It is known that $\lambda_{1}( {T} ) > \mu_{2}( {T} ) = \mu$ (see \autoref{lma24Pol}). Applying the maximum principle (see \autoref{lma21BNV}) to the linear equation of $\partial_{{x}_{2}}{u}$, one deduces the positivity of $\partial_{{x}_{2}}{u}$ in ${T}$. Now \autoref{lma35} implies that ${u}$ has a unique critical point (a saddle point) on $\operatorname{Int}(\overline{{z}_{1}{z}_{2}})$. 

Lastly, the simplicity of the second eigenvalue follows from conclusion (ii). Suppose that the eigenvalue is not simple. Then, by taking a suitable linear combination of two linearly independent eigenfunctions, we can construct an eigenfunction that vanishes at a vertex. This leads to a contradiction, thereby establishing conclusion (iv). 
\end{proof}

The second eigenfunction ${u}$ is monotone in multiple directions (see \eqref{CGY0205}) for non-acute triangles. This multiplicity of monotonicity directions is valid for acute triangles only if ${u}$ vanishes at a vertex.

\begin{lma} \label{lma37}
Let ${T}$ be an (acute or non-acute) triangle $\triangle {z}_{1}{z}_{2}{z}_{3}$. If a second Neumann eigenfunction ${u}$ of ${T}$ satisfies ${u}({z}_{1}) = 0$ and
\begin{equation} \label{CGY0315}
\nabla {u} \cdot \mathbf{a} > 0 \text{ in } {T}
\end{equation}
for some unit vector $\mathbf{a} \in \R^{2}$, then
\begin{equation*}
\nabla {u} \cdot \boldsymbol{n}_{S} \neq 0, \quad
\nabla {u} \cdot \boldsymbol{n}_{M} \neq 0, \quad
\nabla {u} \cdot \boldsymbol{\tau}_{L} \neq 0 \text{ in } {T}. 
\end{equation*}
\end{lma}

\begin{proof}
For a non-acute triangle ${T}$, the monotonicity property holds non-trivially as in \autoref{prop28obt}, without requiring the assumptions ${u}({z}_{1}) = 0$ and \eqref{CGY0315}. 

We now consider the case where ${T}$ is an acute triangle. The monotonicity condition \eqref{CGY0315} and the Neumann boundary condition of ${u}$ imply that ${u}$ is monotone along any side of the triangle that is not perpendicular to $\mathbf{a}$. Clearly, $\mathbf{a}$ cannot be the normal direction to the side $\overline{{z}_{2}{z}_{3}}$; otherwise, monotonicity would imply that ${z}_{1}$ must be a global extremum of ${u}$, which contradicts the fact that ${u}({z}_{1}) = 0$ and \autoref{lma26bJM}. Thus, ${u}$ must be monotone along $\overline{{z}_{2}{z}_{3}}$. Since there are at least two sides that are not perpendicular to $\mathbf{a}$, we know that ${u}$ is monotone along at least two sides of the triangle. Without loss of generality, assume that ${u}$ is monotone along the side $\overline{{z}_{2}{z}_{1}}$. Note that ${u}$ cannot vanish at two vertices simultaneously (see \autoref{lma26aJM}), which implies that ${u}({z}_{2}) \neq 0$. After possibly multiplying ${u}$ by $ - 1$, we may assume that ${u}({z}_{2}) < 0$. Combining this with \autoref{lma26bJM}, we derive that $\nabla {u} \cdot \boldsymbol{n}_{M} \leq, \not\equiv 0$ on $\overline{{z}_{2}{z}_{1}} \cup \overline{{z}_{2}{z}_{3}}$ and hence on $\partial{T}$. Using the fact that $\lambda_{1}(T) > \mu_{2}(T) = \mu$ (see \autoref{lma24Pol}) and applying the maximum principle as in \autoref{lma21BNV}, we conclude that
\begin{equation*}
\nabla {u} \cdot \boldsymbol{n}_{M} < 0 \text{ in } {T}. 
\end{equation*}
Based on this and the fact that ${u}({z}_{1}) = 0$, \autoref{lma35} implies that ${u}$ does not have any non-vertex critical points. Additionally, it follows that $\nabla {u} \cdot \boldsymbol{n}_{S} \geq, \not\equiv0$ on $\partial{T}$. Again, applying the maximum principle from \autoref{lma21BNV}, we conclude that $\nabla {u} \cdot \boldsymbol{n}_{S} > 0$ in ${T}$. Consequently, it follows that $\nabla {u} \cdot \boldsymbol{\tau}_{L} > 0$ in ${T}$. This completes the proof. 
\end{proof}


\subsection{Uniqueness of the critical point}

\begin{thm} \label{thm38}
Let ${T}=\Delta {z}_{1}{z}_{2}{z}_{3}$ be an acute triangle with vertices ${z}_{1}$, ${z}_{2}$, and ${z}_{3}$, and let ${u}$ be a second Neumann eigenfunction of ${T}$. Suppose that the lower (horizontal) side $\overline{{z}_{1}{z}_{2}}$ is the shortest side of ${T}$, and that ${u}$ satisfies
\begin{equation} 
\partial_{{x}_{2}}{u} \geq 0 \text{ in } {T}. 
\end{equation}
Then ${u}$ has at most one non-vertex critical point. Moreover, if such a critical point exists, it is a saddle point lying in the interior of the shortest side $\overline{{z}_{1}{z}_{2}}$. This unique non-vertex critical point exists if and only if ${u}$ does not vanish at any vertex of ${T}$. 
\end{thm}

\begin{proof}
From \autoref{lma35}, the theorem holds if and only if $\operatorname{crit}_{\mathrm{nv}}({u}) \cap (\overline{{z}_{3}{z}_{1}} \cup \overline{{z}_{3}{z}_{2}}) = \emptyset$. It therefore suffices to prove $\operatorname{crit}_{\mathrm{nv}}({u}) \cap \overline{{z}_{3}{z}_{1}} = \emptyset$, since the case of $\overline{{z}_{3}{z}_{2}}$ is entirely analogous. 

We claim that $ \operatorname{crit}_{\mathrm{nv}}({u}) \cap \overline{{z}_{3}{z}_{1}} = \emptyset$ provided that $|\overline{{z}_{1}{z}_{2}}| \leq |\overline{{z}_{2}{z}_{3}}|$. 
We argue by contradiction and suppose that there is a point ${q} \in \operatorname{crit}_{\mathrm{nv}}({u}) \cap \overline{{z}_{3}{z}_{1}}$. By \autoref{lma35}, it follows that $\operatorname{crit}_{\mathrm{nv}}({u}) \cap \operatorname{Int}(\overline{{z}_{1}{z}_{2}})$ contains exactly one element, and ${u}({z}_{1}) < 0$, ${u}({z}_{2}) < 0$. 
Denote the derivative of ${u}$ along the side $\overline{{z}_{3}{z}_{1}}$ by
\begin{equation*}
{w} = \nabla {u} \cdot \boldsymbol{\tau}_{M} = {c}_{1}\partial_{{x}_{1}}{u} + {c}_{2}\partial_{{x}_{2}}{u}, 
\end{equation*}
where ${c}_{1} = - \cos\alpha_{1} < 0$, ${c}_{2} = - \sin\alpha_{1} < 0$ and $\alpha_{1}$ is the interior angle at the vertex ${z}_{1}$ of the triangle ${T}$. 
Let ${z}_{4}$ be the point on the side $\overline{{z}_{3}{z}_{1}}$ such that the line segment $\overline{{z}_{4}{z}_{2}}$ is perpendicular to the side $\overline{{z}_{3}{z}_{1}}$; see \autoref{fig3B}. 

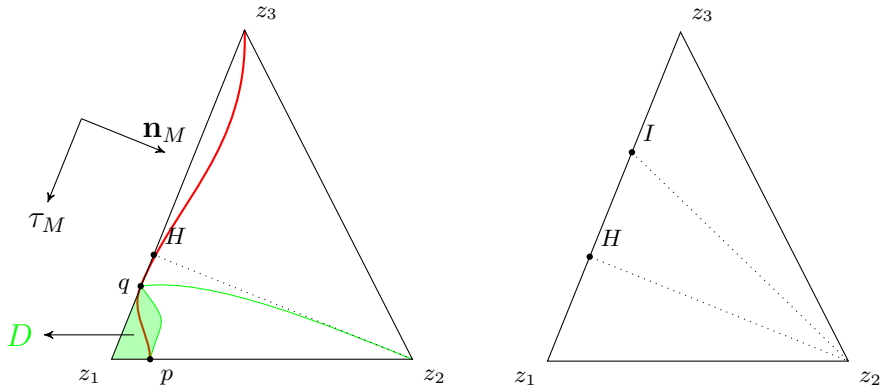
\begin{figure}[htp]\centering  
\begin{tikzpicture}[scale = 3.5]
\pgfmathsetmacro\LenBase{1.00}; 
\pgfmathsetmacro\AngleL{68};  
\pgfmathsetmacro\AngleM{63};  
\pgfmathsetmacro\AngleS{180-\AngleL-\AngleM};  
\pgfmathsetmacro\LenLeft{\LenBase*sin(\AngleM)/sin(\AngleS)}; 
\pgfmathsetmacro\LenRight{\LenBase*sin(\AngleL)/sin(\AngleS)}; 
\pgfmathsetmacro\xSS{\LenLeft*cos(\AngleL)};   
\pgfmathsetmacro\ySS{\LenLeft*sin(\AngleL)}; 
\pgfmathsetmacro\xLL{0};   
\pgfmathsetmacro\yLL{0}; 
\pgfmathsetmacro\xMM{\LenBase};   
\pgfmathsetmacro\yMM{0};  
\pgfmathsetmacro\xHH{\LenBase*cos(\AngleL)*cos(\AngleL)}; 
\pgfmathsetmacro\yHH{\LenBase*cos(\AngleL)*sin(\AngleL)}; 
\pgfmathsetmacro\xQQ{\LenBase*cos(\AngleL)*cos(\AngleL)*0.6901};  
\pgfmathsetmacro\yQQ{\LenBase*cos(\AngleL)*sin(\AngleL)*0.7001}; 
\pgfmathsetmacro\xPP{\xQQ*1.3231};  
\pgfmathsetmacro\yPP{0}; 
\pgfmathsetmacro\xIa{\xPP+cos(90)*\LenBase/12};  
\pgfmathsetmacro\yIa{\yPP+sin(90)*\LenBase/12}; 
\pgfmathsetmacro\xIb{\xQQ+cos(\AngleL-180)*\LenLeft/12};  
\pgfmathsetmacro\yIb{\yQQ+sin(\AngleL-180)*\LenLeft/12};  
\pgfmathsetmacro\xIc{\xQQ+cos(\AngleL)*\LenLeft/4};  
\pgfmathsetmacro\yIc{\yQQ+sin(\AngleL)*\LenLeft/4}; 
\pgfmathsetmacro\xId{\LenLeft*cos(\AngleL)+cos(-88)*\LenLeft/3};  
\pgfmathsetmacro\yId{\LenLeft*sin(\AngleL)+sin(-88)*\LenLeft/3}; 
\pgfmathsetmacro\xHa{\xPP+cos(75)*\LenBase/8};  
\pgfmathsetmacro\yHa{\yPP+sin(75)*\LenBase/8};  
\pgfmathsetmacro\xHb{\xQQ+cos(\AngleL-135)*\LenLeft/6};  
\pgfmathsetmacro\yHb{\yQQ+sin(\AngleL-135)*\LenLeft/6};  
\pgfmathsetmacro\xHc{\xQQ+cos(\AngleL-45)*\LenLeft/6};  
\pgfmathsetmacro\yHc{\yQQ+sin(\AngleL-45)*\LenLeft/6}; 
\pgfmathsetmacro\xHd{\LenBase+cos(\AngleL+95)*\LenLeft/4};  
\pgfmathsetmacro\yHd{0+sin(\AngleL+95)*\LenLeft/4}; 
\draw[green, thick] (\xQQ, \yQQ) .. controls (\xHc, \yHc) and (\xHd, \yHd) .. (\LenBase, 0); 
\draw[green, thick] (\xQQ, \yQQ) .. controls (\xHb, \yHb) and (\xHa, \yHa) .. (\xPP, \yPP); 
\draw[red] (\xQQ, \yQQ) .. controls (\xIb, \yIb) and (\xIa, \yIa) .. (\xPP, \yPP); 
\draw[red] (\xQQ, \yQQ) .. controls (\xIc, \yIc) and (\xId, \yId) .. (\xSS, \ySS); 
\fill[fill = green, fill opacity = 0.3]
(\xQQ, \yQQ) .. controls (\xHb, \yHb) and (\xHa, \yHa) .. (\xPP, \yPP) -- (0, 0) -- cycle; 
\draw[black] 
(\xSS, \ySS) node[below = 3pt, right]{\small ${z}_{3}$}
-- (\xMM, \yMM) node[below]{\small ${z}_{2}$}
-- (\xLL, \yLL) node[below]{\small ${z}_{1}$} -- cycle; 
\fill (\xPP, \yPP) circle (0.3pt) node[below]{\small ${p}$}; 
\fill (\xQQ, \yQQ) circle (0.3pt) node[left]{\small ${q}$}; 
\fill (\xHH, \yHH) circle (0.3pt) node[above left]{\small ${z}_{4}$}; 
\draw[->] ({(0+\xPP+\xQQ)/3}, {(0+\yPP+\yQQ)/3}) -- ++ ({-\LenBase*0.21}, {0}) node[left] { \color{green}$\mathcal{D}_{ - }$}; 
\node at ({(\xMM+\xPP+\xQQ)/3}, {(\yMM+\yPP+\yQQ)/3}) {$\mathcal{D}_{ + }$}; 
\draw[dotted] (\xMM, \yMM) -- (\xHH, \yHH); 
\draw[->] (-\LenBase*0.1, \LenBase*0.8) -- ++ ({180+\AngleL} : 0.3*\LenBase) node [below] {$\mathbf{\tau}_{M}$}; 
\draw[->] (-\LenBase*0.1, \LenBase*0.8) -- ++ ({270+\AngleL} : 0.3*\LenBase) node [above] {$\mathbf{n}_{M}$}; 
\end{tikzpicture} \hspace{2ex}
\begin{tikzpicture}[scale = 3.5]
\pgfmathsetmacro\LenBase{1.00}; 
\pgfmathsetmacro\AngleL{68};  
\pgfmathsetmacro\AngleM{63};  
\pgfmathsetmacro\AngleS{180-\AngleL-\AngleM};  
\pgfmathsetmacro\LenLeft{\LenBase*sin(\AngleM)/sin(\AngleS)}; 
\pgfmathsetmacro\LenRight{\LenBase*sin(\AngleL)/sin(\AngleS)}; 
\pgfmathsetmacro\xSS{\LenLeft*cos(\AngleL)}; 
\pgfmathsetmacro\ySS{\LenLeft*sin(\AngleL)}; 
\pgfmathsetmacro\xLL{0}; 
\pgfmathsetmacro\yLL{0}; 
\pgfmathsetmacro\xMM{\LenBase}; 
\pgfmathsetmacro\yMM{0}; 
\pgfmathsetmacro\xHH{\LenBase*cos(\AngleL)*cos(\AngleL)}; 
\pgfmathsetmacro\yHH{\LenBase*cos(\AngleL)*sin(\AngleL)}; 
\draw[dotted] (2*\xHH, 2*\yHH) -- (\LenBase, 0) -- (\xHH, \yHH); 
\draw[black] 
(\xSS, \ySS) node[below = 3pt, right]{\small ${z}_{3}$}
-- (\xMM, \yMM) node[below]{\small ${z}_{2}$}
-- (\xLL, \yLL) node[below]{\small ${z}_{1}$} -- cycle; 
\fill (2*\xHH, 2*\yHH) circle (0.3pt) node[left]{\small ${z}_{5}$}; 
\fill (\xHH, \yHH) circle (0.3pt) node[left]{\small ${z}_{4}$}; 
\end{tikzpicture} \vspace*{-2ex}
\caption{Nodal domain $\color{green}\mathcal{D}_{ - }$ of $\nabla {u} \cdot \boldsymbol{\tau}_{M}$ near ${z}_{1}$} 
\label{fig3B}
\end{figure}

Let us denote
\begin{equation} 
{T}_{ \star } = \triangle {z}_{1}{z}_{2}{z}_{4} \text{ and } {T}^{ \star } = \triangle {z}_{2}{z}_{3}{z}_{4}. 
\end{equation}
From \eqref{CGY0306a}, we know that
\begin{equation} \label{CGY0323}
{w} < 0 \text{ in } {T}^{ \star }
\end{equation}
and ${w} < 0$ on $\operatorname{Int}(\overline{{z}_{3}{z}_{4}}) \cup \{{z}_{4}\}$. Hence ${q} \in \operatorname{Int}(\overline{{z}_{4}{z}_{1}})$; see \autoref{fig3B}. By applying the Hopf lemma to the positive function $\partial_{{x}_{2}}{u}$, we get
$\partial_{{x}_{2}{x}_{2}}{u} > 0$ on $\operatorname{Int}(\overline{{z}_{1}{z}_{2}})$, which implies
\begin{equation*}
- \partial_{{x}_{2}}{w} = - {c}_{1}\partial_{{x}_{1}{x}_{2}}{u} - {c}_{2}\partial_{{x}_{2}{x}_{2}}{u} = - {c}_{2}\partial_{{x}_{2}{x}_{2}}{u} > 0 \text{ on } \operatorname{Int}(\overline{{z}_{1}{z}_{2}}). 
\end{equation*}
By the local expansion near the vertex ${z}_{1}$ and the fact that ${u}({z}_{1}) < 0$, it follows (see e.g., part (2b) of \autoref{lma26bJM}) that ${w} = \nabla {u} \cdot \boldsymbol{\tau}_{M}$ is nonzero in the intersection of $\overline{{T}}$ with some punctured neighborhood of the vertex ${z}_{1}$. Combining this with the first part of \autoref{lma26bJM}, we conclude that ${w}$ is negative in this intersection. Now let $\mathcal{D}_{ - }$ be the nodal domain of ${w}$ such that $\partial\mathcal{D}_{ - }$ contains the vertex ${z}_{1}$. Then ${w}$ is negative in $\mathcal{D}_{ - }$. 
From the local behavior of ${w}$ near the critical point ${q}$, there exists a nodal domain $\mathcal{D}_{ + }$ of ${w}$ in ${T}$ whose boundary $\partial \mathcal{D}_{ + }$ contains ${q}$, and ${w} > 0$ in $\mathcal{D}_{+}$. By \eqref{CGY0323}, we have $\mathcal{D}_{ + } \subset {T}_{ \star }$. Thanks to \autoref{lma25LR}, the set $\mathcal{Z}({w}) \cup \overline{{z}_{1}{z}_{3}}$ cannot contain a loop. Consequently, $\partial \mathcal{D}_{ + }$ meets both the point ${q}$ and some point on $\overline{{z}_{1}{z}_{2}}$. Therefore, $\mathcal{D}_{ + }$ separates $\mathcal{D}_{ - }$ and ${T}^{ \star }$, and hence $\mathcal{D}_{ - }$ is a strict subset of triangle ${T}_{ \star }$ (see \autoref{fig3B}), and ${w}$ satisfies 
\begin{equation} \label{CGY0324}
\begin{cases}
\Delta {w} + \mu {w} = 0 & \text{ in } \mathcal{D}_{ - }, 
\\ 
{w} = 0 & \text{ on } \partial\mathcal{D}_{ - } \cap {T}, 
\\ 
\partial_{\nu}{w} \geq, \not\equiv 0 & \text{ on } \partial\mathcal{D}_{ - } \cap \partial {T}. 
\end{cases}\end{equation} 
Recall that $\lambda_{1}(\mathcal{D}_{ - }, \partial\mathcal{D}_{ - } \cap {T})$ denotes the first eigenvalue of the Laplacian operator on the domain $\mathcal{D}_{ - }$ under mixed boundary conditions -- homogeneous Dirichlet boundary conditions on $\partial\mathcal{D}_{ - } \cap {T}$ and homogeneous Neumann conditions elsewhere. 
From \autoref{lma22cMP}, we obtain that 
\begin{equation*}
\lambda_{1}(\mathcal{D}_{ - }, \partial\mathcal{D}_{ - } \cap {T}) < \mu. 
\end{equation*}
By the variational characterization of the first eigenvalue, we deduce that
\begin{equation*}
\begin{aligned}
\lambda_{1}(\mathcal{D}_{ - }, \partial\mathcal{D}_{ - } \cap {T}) & = 
\inf\bigg\{\frac{\int_{\mathcal{D}_{ - }} |\nabla\varphi|^{2} dx}{\int_{\mathcal{D}_{ - }} |\varphi|^{2} dx}: \, 0 \neq \varphi \in W^{1, 2}(\mathcal{D}_{ - }), \; \varphi = 0 \text{ on } \partial\mathcal{D}_{ - } \cap {T}
\bigg\}
\\ & = 
\inf\bigg\{\frac{\int_{{T}_{ \star }} |\nabla\varphi|^{2} dx}{\int_{{T}_{ \star }} |\varphi|^{2} dx}: \, 0 \neq \varphi \in W^{1, 2}({T}_{ \star }), \; \varphi = 0 \text{ in } {T}_{ \star } \setminus \mathcal{D}_{ - }
\bigg\}
\\ & > 
\inf\bigg\{\frac{\int_{{T}_{ \star }} |\nabla\varphi|^{2}dx}{\int_{{T}_{ \star }} |\varphi|^{2} dx}: \, 0 \neq \varphi \in W^{1, 2}({T}_{ \star }), \; \varphi = 0 \text{ on } \overline{{z}_{2}{z}_{4}}
\bigg\}
\\ & = 
\lambda_{1}({T}_{ \star }, \overline{{z}_{4}{z}_{2}}).
\end{aligned}
\end{equation*} 
Thus, we have
\begin{equation} \label{CGY0326}
\lambda_{1}({T}_{ \star }, \overline{{z}_{4}{z}_{2}}) < \mu. 
\end{equation}

In order to reach a contradiction, we will show that
\begin{equation} \label{CGY0327}
\lambda_{1}({T}_{ \star }, \overline{{z}_{4}{z}_{2}}) \geq \mu \text{ if } |\overline{{z}_{1}{z}_{2}}| \leq |\overline{{z}_{3}{z}_{2}}|. 
\end{equation}
Indeed, from \cite[Lemma 2.1]{BB99}, we have
\begin{equation*}
\mu = \mu_{2}({T}) \leq \max\{\lambda_{1}({T}_{ \star }, \overline{{z}_{4}{z}_{2}}), \; \lambda_{1}({T}^{ \star }, \overline{{z}_{4}{z}_{2}})\}. 
\end{equation*}
Now we let ${z}_{5}$ be the point such that ${z}_{4}$ is the midpoint of the line segment $\overline{{z}_{1}{z}_{5}}$; see the right picture in \autoref{fig3B}. The condition $|\overline{{z}_{1}{z}_{2}}| \leq |\overline{{z}_{3}{z}_{2}}|$ implies that either ${z}_{5} = {z}_{3}$ or ${z}_{5} \in \operatorname{Int}(\overline{{z}_{3}{z}_{4}})$. 
Denote by $\varphi^{ \star } > 0$ the eigenfunction corresponding to the first eigenvalue $\lambda^{ \star } = \lambda_{1}({T}^{ \star }, \overline{{z}_{4}{z}_{2}})$, that is, 
\begin{equation*}
\Delta \varphi^{ \star } + \lambda^{ \star }\varphi^{ \star } = 0 \text{ in } {T}^{ \star }, \quad
\varphi^{ \star } = 0 \text{ on } \overline{{z}_{4}{z}_{2}}, \text{ and } \partial_{\nu}\varphi^{ \star } = 0 \text{ on } \overline{{z}_{4}{z}_{3}} \cup \overline{{z}_{2}{z}_{3}}. 
\end{equation*}
From \cite[Theorem 1.1]{JN00} (see also \autoref{lma41} below), $\varphi^{ \star }$ is monotone along the normal direction to the Dirichlet boundary $\overline{{z}_{4}{z}_{2}}$, i.e., $\nabla \varphi^{ \star } \cdot \boldsymbol{\tau}_{M} < 0$ in ${T}^{ \star }$. Using the same process as in the proof of \autoref{lma31}, we know that the function ${R}_{{z}_{2}}\varphi^{ \star }$ is negative in ${T}^{ \star }$. 
Thus, $\nabla \varphi^{ \star } \cdot \nu_{\overline{{z}_{5}{z}_{2}}} \geq 0 $ on $\overline{{z}_{5}{z}_{2}}$, where $\nu_{\overline{{z}_{5}{z}_{2}}} = (\sin(2\alpha_{1}), - \cos(2\alpha_{1}))$ denotes the unit normal vector to the line segment $\overline{{z}_{5}{z}_{2}}$. 
Therefore, 
\begin{equation*}
\Delta \varphi^{ \star } + \lambda^{ \star }\varphi^{ \star } = 0 \text{ in } \triangle {z}_{5}{z}_{4}{z}_{2}, \quad
\varphi^{ \star } = 0 \text{ on } \overline{{z}_{4}{z}_{2}}, \quad \partial_{\nu}\varphi^{ \star } \geq 0 \text{ on } \overline{{z}_{4}{z}_{5}} \cup \overline{{z}_{2}{z}_{5}}. 
\end{equation*}
From \autoref{lma22cMP}, we obtain that the first eigenvalue $\lambda_{1}(\triangle {z}_{5}{z}_{4}{z}_{2}, \overline{{z}_{4}{z}_{2}})$ is greater than or equal to $\lambda^{ \star }$, that is, 
\begin{equation*}
\lambda_{1}({T}_{ \star }, \overline{{z}_{4}{z}_{2}}) = \lambda_{1}(\triangle {z}_{5}{z}_{4}{z}_{2}, \overline{{z}_{4}{z}_{2}}) \geq \lambda_{1}({T}^{ \star }, \overline{{z}_{4}{z}_{2}}). 
\end{equation*}
Therefore, we obtain \eqref{CGY0327}. 

From \eqref{CGY0326} and \eqref{CGY0327}, we reach a contradiction, thus establishing the nonexistence of a critical point on $\operatorname{Int}(\overline{{z}_{3}{z}_{1}})$. This completes the proof. 
\end{proof}


\section{Eigenvalue inequalities} \label{Sect4EI}

In this section, we establish eigenvalue inequalities used to determine whether the second Neumann eigenfunction vanishes at a vertex of the triangle in \autoref{Sect5nodal} below.


\subsection{Inequalities for the first mixed eigenvalue on triangles} \label{Sect4a}

In this subsection we prove the eigenvalue inequalities conjectured by Siudeja \cite{Siu16}. 
Let ${T} = {T}_{{a}, {b}}$ be a triangle ${z}_{1}{z}_{2}{z}_{3}$ with vertices ${z}_{1} = (0, 0)$, ${z}_{2} = (1, 0)$ and ${z}_{3} = ({a}, {b})$ where ${a} \in \R$ and ${b} \in \R^{ + }$. Let $\lambda = \lambda_{{a}, {b}}$ and $\varphi = \varphi_{{a}, {b}}$ be the first eigenvalue and eigenfunction of the mixed boundary problem
\begin{equation} \label{CGY0401}
\begin{cases}
\Delta \varphi + \lambda\varphi = 0 & \text{ in } {T} = \triangle {z}_{1}{z}_{2}{z}_{3}, \\
\varphi = 0 & \text{ on } \Gamma_{D} = \overline{{z}_{3}{z}_{2}}, \\
\partial_{\nu}\varphi = 0 & \text{ on } \Gamma_{N} = \partial {T} \setminus \Gamma_{D}. 
\end{cases}\end{equation}
We assume the eigenfunction $\varphi_{{a}, {b}}$ is positive and normalized: 
\begin{equation*}
\int_{{T}_{{a}, {b}}} |\varphi_{{a}, {b}}|^{2} dx = 1. 
\end{equation*}

In a recent study by the third author \cite{LY24}, the moving plane method was successfully employed to derive monotonicity results for positive solutions of semilinear elliptic equations in triangles with mixed boundary conditions. The linear case of that result can be summarized as follows. 

\begin{lma} \label{lma41}
Let ${T}$ be the triangle ${z}_{1}{z}_{2}{z}_{3}$ and $\varphi$ denote the positive eigenfunction defined above. Suppose that the interior angle at Neumann vertex ${z}_{1}$ of ${T}$ is non-obtuse, i.e., $\angle {z}_{3}{z}_{1}{z}_{2} \leq \pi/2$. Then $\varphi$ has no non-vertex critical points and is monotone in the normal direction to the Dirichlet side $\overline{{z}_{3}{z}_{2}}$. Moreover, if, in addition, $\angle {z}_{1}{z}_{2}{z}_{3} \leq \pi/2$, then
\begin{equation*}
\partial_{{x}_{1}} \varphi < 0 \text{ and } \partial_{{x}_{2}} \varphi < 0 \text{ in } {T}. 
\end{equation*}
\end{lma}

\begin{proof}
This result follows from \autoref{prop29Yao}. We provide an alternative proof using the continuity method with domain deformation. For simplicity, fix ${b} > 0$ and omit the subscript ${b}$. Let ${T}_{{t}}$ denote the triangle ${T}_{{t}, {b}}$ with vertices ${z}_{1} = (0, 0)$, ${z}_{2} = (1, 0)$ and ${z}_{3}^{{t}} = ({t}, {b})$. Define $\lambda_{{t}} = \lambda_{{t}, {b}}$ and $\varphi_{{t}} = \varphi_{{t}, {b}}$. The boundaries $\Gamma_{D}^{{t}}$ and $\Gamma_{N}^{{t}}$ represent the corresponding Dirichlet and Neumann sides, respectively. Since $\lambda_{{t}}$ is a simple eigenvalue, the mappings ${t} \mapsto \lambda_{{t}}$ and ${t} \mapsto \varphi_{{t}}$ are continuous. The fact that the interior angle at the vertex ${z}_{1}$ of ${T}_{t}$ is non-obtuse implies that ${t} \geq 0$. Define
\begin{equation*}
{w}_{t}^{ - } = \nabla \varphi_{t} \cdot (1, 0), \quad
{w}_{t}^{ + } = \nabla \varphi_{t} \cdot ({t}, {b}), \quad
{w}_{t} = \nabla \varphi_{t} \cdot ({b}, 1 - {t}). 
\end{equation*}
These functions are the directional derivatives of $\varphi_{t}$ along directions parallel to the two Neumann sides and along the outward unit normal on the Dirichlet side, respectively. They satisfy the same linear equation as $\varphi_{t}$ without boundary conditions.

Recall that ${T}_{1}$ is a right triangle, and the interior angle at the mixed boundary point ${z}_{2}$ is $\pi/2$. By applying \cite[Theorem 1.1]{JN00}, we deduce that
\begin{equation} \label{CGY0402a}
{w}_{t}^{ - } < 0, \quad {w}_{t}^{ + } < 0, \quad {w}_{t} < 0 \text{ in } \overline{{T}_{t}} \setminus \{\text{vertices}\}
\end{equation}
holds when ${t} = 1$. Assume for contradiction that $\tau > 0$, where
\begin{equation*} 
\tau = \inf\{{t}' \in [0, 1]: \, \eqref{CGY0402a} \text{ holds for every } {t} \in [{t}', 1]\}. 
\end{equation*}
By continuity, we have ${w}_{\tau}^{\pm} \leq 0$ and ${w}_{\tau} \leq 0$ in ${T}_{\tau}$. The strong maximum principle implies
\begin{equation*} 
{w}_{\tau}^{ - } < 0, \quad {w}_{\tau}^{ + } < 0, \quad {w}_{\tau} < 0 \text{ in } {T}_{\tau} \cup \operatorname{Int}(\Gamma_{D}^{\tau}). 
\end{equation*}
Following the arguments in the proof of \autoref{lma31}(iii), we deduce that the tangential derivatives along both Neumann boundaries are nonzero, thereby ensuring that \eqref{CGY0402a} is valid for ${t} = \tau$. 
Note that $\varphi_{t}$ is positive. Utilizing the same arguments as in \cite[page 752]{JN00} or \cite{JM22c, Hat24}, we analyze the local behavior at all three vertices. We can then find two small positive constants $\bar{\delta} > 0$ and $\varepsilon_{1} \in (0, \tau)$ such that
\begin{equation*}
{w}_{t}^{ - } < 0, \; {w}_{t}^{ + } < 0, \; {w}_{t} < 0 \text{ in } \mathcal{V}_{\bar{\delta}}^{t} = \{{x} \in \overline{{T}_{t}}: \, 0 < \operatorname{dist}({x}, \{{z}_{1}, {z}_{2}, {z}_{3}^{t}\}) < \bar{\delta}\}
\end{equation*}
for $|{t} - \tau| < \varepsilon_{1}$. From the continuity of the solution and \eqref{CGY0402a} at ${t} = \tau$, there exists a constant $\varepsilon_{2} \in (0, \varepsilon_{1})$ such that for $|{t} - \tau| < \varepsilon_{2}$, 
\begin{equation*}
{w}_{t}^{ - } < 0, \; {w}_{t}^{ + } < 0, \; {w}_{t} < 0 \text{ in } \overline{{T}_{t}} \setminus \mathcal{V}_{\bar{\delta}}^{t}. 
\end{equation*}
Thus, \eqref{CGY0402a} is valid for $|{t} - \tau| < \varepsilon_{2}$, contradicting the definition of $\tau$. Therefore, $\tau = 0$, and \eqref{CGY0402a} holds for ${t} \in (0, 1]$. Additionally, $\varphi_{0}$ has nonzero tangential derivatives on the interior of both Neumann boundaries; hence
\begin{equation*} 
\partial_{{x}_{1}}\varphi_{0} < 0, \; \partial_{{x}_{2}}\varphi_{0} < 0, \; \nabla \varphi_{t} \cdot ({b}, 1) < 0 \text{ in } {T}_{0} \cup \operatorname{Int}(\Gamma_{D}^{0}). 
\end{equation*}
By similar arguments, we can also establish that \eqref{CGY0402a} holds for ${t} \geq 1$. 

We now show that $\partial_{{x}_{2}} \varphi_{t}$ is negative in ${T}_{t}$ for any fixed ${t} \in [0, 1]$. Indeed, since ${t}\in[0,1]$, we have $\partial_{{x}_{2}} \varphi_{t} \leq 0$ on the interior of the Dirichlet boundary, while \eqref{CGY0402a} ensures the negativity of $\partial_{{x}_{2}} \varphi_{t}$ on the interior of the upper Neumann boundary. From the local behavior near vertices, one has $\partial_{{x}_{2}} \varphi_{t} < 0$ in $\mathcal{V}^{t} \cap {T}_{t}$ for some neighborhood $\mathcal{V}^{t}$ of the three vertices of ${T}_{t}$. 
Therefore, $\partial_{{x}_{2}} \varphi_{t}$ satisfies the same linear equation as $\varphi_{t}$ and
\begin{equation*}
\partial_{{x}_{2}} \varphi_{t} \leq, \not\equiv 0 \text{ on } \partial ({T}_{t} \setminus \mathcal{V}^{t}).
\end{equation*}
Recalling that $\lambda_{1}({T}_{t} \setminus \mathcal{V}^{t}) > \lambda_{1}({T}_{t}) > \lambda_{t}$, we apply the maximum principle as stated in \autoref{lma21BNV} to $\partial_{{x}_{2}} \varphi_{t}$. This yields the negativity of $\partial_{{x}_{2}} \varphi_{t}$ in ${T}_{t}$. Thus, the proof is complete. 
\end{proof}

\begin{rmks}
Let $\varphi > 0$ be the first mixed Dirichlet-Neumann eigenfunction in a triangle ${T}$ with Dirichlet boundary $\Gamma_{D}$ being a side of ${T}$. 
\begin{enumerate}
\item[\rm(1)]
If the Neumann vertex is non-obtuse, then $\varphi$ is monotone in the normal direction to the Dirichlet side (see \autoref{lma41}). 
\item[\rm(2)]
If the Neumann vertex is obtuse, then $\varphi$ is monotone in the normal direction to the longer Neumann side. This can be proven using a similar approach to that employed in \autoref{lma41}. For further details, please refer to \cite{LY24}. 
\end{enumerate}

In both cases, the global maximum point is unique and lies on the longer (closed) Neumann side of ${T}$. This answers a question posed in \cite{Pol12}. When the Dirichlet boundary comprises two sides of the triangle, a similar property has been studied in a recent paper \cite{Hat24}. 
\end{rmks}

\begin{lma} \label{lma42}
Fix ${b} > 0$ and denote by $\lambda_{{a}, {b}}$ the first eigenvalue of \eqref{CGY0401}. Then the function ${a} \in [0, 1] \mapsto \lambda_{{a}, {b}}$ is strictly decreasing. 
\end{lma}

\begin{proof}
The lemma immediately follows from the following claim: for any ${a} \in [0, 1]$, there exists a small positive constant $\bar{\delta} = \bar{\delta}_{{a}, {b}}$ such that 
\begin{equation} \label{CGY0404}
\lambda_{{a} + {b}\delta, {b}} < \lambda_{{a}, {b}} \text{ for every } 0 < \delta < \bar{\delta}_{{a}, {b}}. 
\end{equation}

To establish \eqref{CGY0404}, let us denote $\Omega = {T}_{{a} + {b}\delta, {b}}$ and ${T} = {T}_{{a}, {b}}$. Let $\varphi = \varphi_{{a}, {b}}$ and define
\begin{equation*}
\tilde{\varphi}({y}_{1}, {y}_{2}) = \varphi({x}_{1}, {x}_{2}) \text{ with } {x}_{1} = {y}_{1} - \delta {y}_{2}, \; {x}_{2} = {y}_{2}. 
\end{equation*}
Since $\varphi = \varphi_{{a}, {b}}$ is an eigenfunction on ${T}_{{a}, {b}}$, the function $\tilde{\varphi}$ is well-defined on $\Omega = T_{{a} + {b}\delta, {b}}$ and vanishes on the right side. By direct computation, we have
\begin{equation*} 
\begin{aligned}
\int_{\Omega} |\nabla\tilde{\varphi}|^{2} dy & = 
\int_{\Omega} \big\{|\partial_{{x}_{1}}\varphi|^{2} + (\partial_{{x}_{2}}\varphi - \delta \partial_{{x}_{1}}\varphi)^{2}\big\} dy 
\\ & = 
\int_{{T}} \big\{|\nabla\varphi|^{2} - 2\delta \partial_{{x}_{1}}\varphi \partial_{{x}_{2}}\varphi + \delta^{2}(\partial_{{x}_{1}}\varphi)^{2}\big\} dx, 
\\
\int_{{T}} |\nabla\varphi|^{2}dx & = 
\lambda_{{a}, {b}}\int_{{T}}|\varphi|^{2}dx = \lambda_{{a}, {b}} \int_{\Omega} |\tilde{\varphi}|^{2} dy. 
\end{aligned}\end{equation*}
Hence, 
\begin{equation} \label{CGY0405}
\int_{\Omega} |\nabla\tilde{\varphi}|^{2} dy - \lambda_{{a}, {b}} \int_{\Omega} |\tilde{\varphi}|^{2} dy
= - 2\delta \int_{{T}} \partial_{{x}_{1}}\varphi \partial_{{x}_{2}}\varphi dx + \delta^{2} \int_{{T}} (\partial_{{x}_{1}}\varphi)^{2} dx. 
\end{equation}
Define 
\begin{equation*} 
\bar{\delta} = \bar{\delta}_{{a}, {b}}
= \frac{2\int_{{T}} \partial_{{x}_{1}}\varphi\partial_{{x}_{2}}\varphi dx}{\int_{{T}} (\partial_{{x}_{1}}\varphi)^{2} dx}. 
\end{equation*}
\autoref{lma41} implies that $\bar{\delta}$ is positive. 
Combining this with \eqref{CGY0405}, we obtain
\begin{equation*} 
\int_{\Omega} |\nabla\tilde{\varphi}|^{2} dy < \lambda_{{a}, {b}} \int_{\Omega} |\tilde{\varphi}|^{2} dy
\end{equation*}
for $0 < \delta < \bar{\delta}$. This inequality implies \eqref{CGY0404}, thus completing the proof. 
\end{proof}

As a direct consequence, we complete the proof of Conjecture 1.2 in \cite{Siu16}.

\begin{thm} \label{thm43Siu}
For any triangle ${T}$ with distinct side lengths, the following eigenvalue inequalities hold: 
\begin{equation*}
\lambda_{1}^{S} < \lambda_{1}^{M} < \lambda_{1}^{L} < \lambda_{1}^{MS} < \lambda_{1}^{LS} < \lambda_{1}^{LM} < \lambda_{1}. 
\end{equation*}
Here the notations are defined before \autoref{thm14Eig}. 
\end{thm}

\begin{proof}
By \cite[Part 2 of Theorem 1.1]{Siu16}, the following chain of inequalities holds: 
\begin{equation*}
\lambda_{1}^{MS} < \lambda_{1}^{LS} < \lambda_{1}^{LM} < \lambda_{1}. 
\end{equation*}
Note that the angles between the longest side and both the medium and shortest sides are strictly less than $\pi/2$. Using \cite[Theorem 3.1]{Roh20} or \cite[Theorem 3.1]{AR23}, we conclude that $\lambda_{1}^{L} < \lambda_{1}^{MS}$. Our remaining goal is to establish
\begin{equation} \label{CGY0407}
\lambda_{1}^{S} < \lambda_{1}^{M} < \lambda_{1}^{L}
\end{equation}
as long as appropriate sides have different lengths. 
From \autoref{lma42}, we have
\begin{equation} \label{CGY0408}
\lambda_{{a}, {b}} > \lambda_{1 - {a}, {b}} \text{ when } {a} \in [0, 1/2). 
\end{equation}
Here, $\lambda_{{a}, {b}}$ and $\lambda_{1 - {a}, {b}}$ denote the first eigenvalues of mixed boundary problems corresponding to the triangles ${T}_{{a}, {b}}$ and ${T}_{1 - {a}, {b}}$, respectively, as defined above \autoref{lma41}. We note the following: 
\begin{itemize}
\item
The triangles ${T}_{{a}, {b}}$ and ${T}_{1 - {a}, {b}}$ are congruent and symmetric with respect to the line ${x}_{1} = 1/2$. 
\item
The Dirichlet side corresponding to $\lambda_{{a}, {b}}$ is longer than that corresponding to $\lambda_{1 - {a}, {b}}$. 
\item
Both the left and right interior angles of ${T}_{{a}, {b}}$ (and hence ${T}_{1 - {a}, {b}}$) are non-obtuse because the condition ${a} \in [0, 1/2)$. 
\end{itemize}
Based on these observations and inequality \eqref{CGY0408}, it immediately follows that
\begin{equation*} 
\lambda_{1}^{S} < \lambda_{1}^{M}
\end{equation*}
since the interior angles between the longest side and the medium side, as well as between the longest side and the shortest side, are acute. Furthermore, if the largest interior angle of the triangle is non-obtuse, then
\begin{equation*} 
\lambda_{1}^{M} < \lambda_{1}^{L}. 
\end{equation*} 
This establishes inequality \eqref{CGY0407} for non-obtuse triangles.

\begin{figure}[h]\centering  \vspace*{-2ex}
\begin{tikzpicture}[scale = 3.5]
\draw[dashed] (0, 0.6) -- (0, 0) node[below]{${z}_{4}$} -- (1, 0); 
\draw[] (0, 0.6) node[below = 4pt, left]{${z}_{3}$} -- (1.8, 0)node[below]{${z}_{2}$} -- (0.7, 0)node[below]{${z}_{1}$} -- cycle; 
\end{tikzpicture} \vspace*{-2ex}
\caption{The eigenvalue inequalities on an obtuse triangle
}
\label{fig4}
\end{figure}

The remaining task is to show that $\lambda_{1}^{M} < \lambda_{1}^{L}$ when the largest interior angle of the triangle is obtuse. Indeed, let ${T}$ be a triangle ${z}_{1}{z}_{2}{z}_{3}$ with $\angle {z}_{3}{z}_{1}{z}_{2} > \pi/2$. We construct the point ${z}_{4}$ on the ray ${z}_{2}{z}_{1}$ such that the line segment $\overline{{z}_{3}{z}_{4}}$ is perpendicular to the side $\overline{{z}_{1}{z}_{2}}$, as illustrated in \autoref{fig4}. Let $\mathcal{T} = \triangle {z}_{2}{z}_{3}{z}_{4}$ and let $\varphi_{*} > 0$ and $\varphi^{*} > 0$ be the eigenfunctions corresponding to the eigenvalues
\begin{equation*}
\lambda_{*} = \lambda_{1}(\mathcal{T}, \overline{{z}_{4}{z}_{2}}) \text{ and } \lambda^{*} = \lambda_{1}(\mathcal{T}, \overline{{z}_{3}{z}_{2}}), 
\end{equation*}
respectively. From \autoref{lma41}, $\varphi_{*}$ (resp., $\varphi^{*}$) is monotonically increasing in the inward normal direction to the corresponding Dirichlet side $\overline{{z}_{4}{z}_{2}}$ (resp., $\overline{{z}_{3}{z}_{2}}$) of $\mathcal{T} = \triangle {z}_{2}{z}_{3}{z}_{4}$. Applying the maximum principle to the angular derivative about the point ${z}_{3}$ (refer to the second part of \autoref{lma31}), we deduce that the normal derivative (pointing outward from ${T}$) of $\varphi_{*}$ (resp., $\varphi^{*}$) on $\operatorname{Int}(\overline{{z}_{3}{z}_{1}})$ is negative (resp., positive). Therefore, 
\begin{gather*}
\Delta \varphi_{*} + \lambda_{*} \varphi_{*} = 0 \text{ in } {T}, \quad
\varphi_{*} = 0 \text{ on } \overline{{z}_{1}{z}_{2}}, \text{ and } 
\partial_{\nu} \varphi_{*} \leq, \not\equiv 0 \text{ on } \overline{{z}_{3}{z}_{1}} \cup \overline{{z}_{3}{z}_{2}}, 
\\
\Delta \varphi^{*} + \lambda^{*} \varphi^{*} = 0 \text{ in } {T}, \quad
\varphi^{*} = 0 \text{ on } \overline{{z}_{3}{z}_{2}}, \text{ and } 
\partial_{\nu} \varphi^{*} \geq, \not\equiv 0 \text{ on } \overline{{z}_{3}{z}_{1}} \cup \overline{{z}_{1}{z}_{2}}. 
\end{gather*}
It then follows \autoref{lma22cMP} that
\begin{equation*} 
\lambda_{1}({T}, \overline{{z}_{1}{z}_{2}}) < \lambda_{*} = \lambda_{1}(\mathcal{T}, \overline{{z}_{4}{z}_{2}}) \text{ and } \lambda_{1}({T}, \overline{{z}_{3}{z}_{2}}) > \lambda^{*} = \lambda_{1}(\mathcal{T}, \overline{{z}_{3}{z}_{2}}). 
\end{equation*}
Since \eqref{CGY0407} is established for any non-obtuse triangle, it follows that
\begin{equation*} 
\lambda_{1}(\mathcal{T}, \overline{{z}_{4}{z}_{2}}) < \lambda_{1}(\mathcal{T}, \overline{{z}_{3}{z}_{2}}). 
\end{equation*}
Combining the above inequalities, we obtain $\lambda_{1}({T}, \overline{{z}_{1}{z}_{2}}) < \lambda_{1}({T}, \overline{{z}_{3}{z}_{2}})$. By a similar argument, we also deduce that $\lambda_{1}({T}, \overline{{z}_{1}{z}_{3}}) < \lambda_{1}({T}, \overline{{z}_{3}{z}_{2}})$. Hence, 
\begin{equation*}
\max\{\lambda_{1}^{S}, \lambda_{1}^{M}\} < \lambda_{1}^{L}, 
\end{equation*}
for any obtuse triangle ${T}$. Therefore, \eqref{CGY0407} is always valid. This completes the proof. 
\end{proof}


\subsection{An eigenvalue inequality for two different triangles}

In this subsection, we always assume that ${T}$ is the triangle ${z}_{1}{z}_{2}{z}_{3}$, $\overline{{z}_{1}{z}_{4}}$ is the internal bisector of the angle of ${T}$ at vertex ${z}_{1}$ with the point ${z}_{4}$ belonging to the side $\overline{{z}_{2}{z}_{3}}$, and ${z}_{5}$ is the reflection point of ${z}_{2}$ across the line $\overline{{z}_{1}{z}_{4}}$. Denote by
\begin{equation}
{T}_{ - } = \triangle {z}_{1}{z}_{2}{z}_{4} \text{ and } {T}_{ + } = \triangle {z}_{1}{z}_{3}{z}_{4}. 
\end{equation}

\begin{lma} \label{lma44}
The following inequalities hold: 

{\rm(i)}
$\mu_{2}({T}) \leq \max\{\lambda_{1}({T}_{ + }, \overline{{z}_{1}{z}_{4}}), \lambda_{1}({T}_{ - }, \overline{{z}_{1}{z}_{4}})\}$. 

{\rm(ii)} 
$\lambda_{1}({T}_{ + }, \overline{{z}_{1}{z}_{4}}) < \lambda_{1}({T}_{ - }, \overline{{z}_{1}{z}_{4}})$ if and only if $|\overline{{z}_{1}{z}_{3}}| > |\overline{{z}_{1}{z}_{2}}|$. 
\end{lma}

\begin{proof}
The first inequality is straightforward yet useful; see the proof in \cite[Lemma 2.1]{BB99}. It is clear that $\lambda_{1}({T}_{ + }, \overline{{z}_{1}{z}_{4}}) = \lambda_{1}({T}_{ - }, \overline{{z}_{1}{z}_{4}})$ if $|\overline{{z}_{1}{z}_{3}}| = |\overline{{z}_{1}{z}_{2}}|$. We only need to prove $\lambda_{1}({T}_{ + }, \overline{{z}_{1}{z}_{4}}) < \lambda_{1}({T}_{ - }, \overline{{z}_{1}{z}_{4}})$ under the condition $|\overline{{z}_{1}{z}_{3}}| > |\overline{{z}_{1}{z}_{2}}|$. Thus, we assume $|\overline{{z}_{1}{z}_{3}}| > |\overline{{z}_{1}{z}_{2}}|$, and hence ${z}_{5} \in \operatorname{Int}(\overline{{z}_{3}{z}_{1}})$. By rotating and translating the coordinate system, if necessary, we may assume that ${z}_{1}$ is at the origin and the internal bisector $\overline{{z}_{1}{z}_{4}}$ lies along the positive ${x}_{1}$-axis. See \autoref{fig5} below. Note that $|\overline{{z}_{1}{z}_{3}}| > |\overline{{z}_{1}{z}_{2}}|$ implies that $\angle {z}_{2}{z}_{3}{z}_{1} < \pi/2$. Let $\varphi > 0$ be the eigenfunction corresponding to the eigenvalue $\lambda_{1}({T}_{ + }, \overline{{z}_{1}{z}_{4}})$ of the mixed boundary value problem. 
From \autoref{lma41}, we know that $\varphi$ has no non-vertex critical points and satisfies $\partial_{{x}_{2}}{\varphi} > 0$ in ${T}_{+}$. Combining this with the boundary conditions for $\varphi$, we know that $\nabla \varphi / |\nabla \varphi| = - {\boldsymbol{\tau}}_{M}$ on $\operatorname{Int}(\overline{{z}_{3}{z}_{1}})$, where ${\boldsymbol{\tau}}_{M}$ denotes the unit tangent vector along $\overline{{z}_{3}{z}_{1}}$, oriented from ${z}_{3}$ to ${z}_{1}$. It follows that the angular derivative ${R}_{{z}_{4}}{\varphi}$ is negative on $\operatorname{Int}(\overline{{z}_{3}{z}_{1}})$. Moreover, by the local behavior near vertices, we have ${R}_{{z}_{4}}{\varphi} < 0$ in $\mathcal{V} \cap {T}_{+}$ for some neighborhood $\mathcal{V}$ of the three vertices of ${T}_{+}$. Thus, ${R}_{{z}_{4}}{\varphi} \leq, \not\equiv 0$ in $\partial ({T}_{+} \setminus \mathcal{V})$. Recalling that $\lambda_{1}({T}_{+} \setminus \mathcal{V}) > \lambda_{1}({T}_{+}) > \lambda_{1}({T}_{+}, \overline{{z}_{1}{z}_{4}})$, we apply the maximum principle in \autoref{lma21BNV} to conclude that ${R}_{{z}_{4}}{\varphi}$ is negative in ${T}_{+}$. In particular, on $\operatorname{Int}(\overline{{z}_{4}{z}_{5}})$, the normal derivative of $\varphi$ (pointing to the right and upward) is positive, 
\begin{equation*} 
(\sin\angle{z}_{1}{z}_{4}{z}_{5}, \cos\angle{z}_{1}{z}_{4}{z}_{5}) \cdot \nabla \varphi > 0 \text{ on } \operatorname{Int}(\overline{{z}_{4}{z}_{5}}). 
\end{equation*}
This leads to
\begin{equation*}\begin{cases}
-\Delta \varphi = \lambda_{1}({T}_{ + }, \overline{{z}_{1}{z}_{4}})\varphi > 0 \text{ in } {T}_{ + }, 
\\
\varphi = 0 \text{ on } \overline{{z}_{1}{z}_{4}}, 
\\
\partial_{\nu}\varphi = 0 \text{ on } \overline{{z}_{1}{z}_{5}}, 
\\
\partial_{\nu}\varphi > 0 \text{ on } \operatorname{Int}(\overline{{z}_{4}{z}_{5}}). 
\end{cases}\end{equation*}
It follows from \autoref{lma22cMP} that $\lambda_{1}(\triangle {z}_{1}{z}_{4}{z}_{5}, \overline{{z}_{1}{z}_{4}}) > \lambda_{1}({T}_{ + }, \overline{{z}_{1}{z}_{4}})$. Hence $\lambda_{1}({T}_{ - }, \overline{{z}_{1}{z}_{4}}) > \lambda_{1}({T}_{ + }, \overline{{z}_{1}{z}_{4}})$. This completes the proof. 
\end{proof}

Let ${j}_{\nu}$ denote the first positive zero of the first kind Bessel function, ${J}_{\nu}$, of order $\nu \geq 0$. In particular, ${j}_{0} \approx 2.4048$. 

\begin{lma} \label{lma45}
For a convex planar domain $\Omega$, the second Neumann eigenvalue $\mu_{2}(\Omega)$ satisfies 
\begin{equation*} 
\mu_{2}(\Omega) < \Big(\frac{2{j}_{0}}{\operatorname{diam}(\Omega)}\Big)^{2}. 
\end{equation*} 
\end{lma}

\begin{proof}
This is obtained in \cite[Corollary 2.1]{BB99}. 
\end{proof}

The isoperimetric inequalities for convex cones in $\R^{n}$ play a fundamental role in the symmetrization of mixed boundary value problems. We state the result in the two-dimensional case only. For $\alpha \in (0, 2\pi)$ and $\rho > 0$, we define the infinite and finite spherical sectors in $\R^{2}$ as follows: 
\begin{align*}
\Sigma_{\alpha} & = \{({r}\cos\theta, {r}\sin\theta) \in \R^{2}: \, \theta \in (0, \alpha), \; 0 < {r} \}, 
\\
\Sigma_{\alpha, \rho} & = \{({r}\cos\theta, {r}\sin\theta) \in \R^{2}: \, \theta \in (0, \alpha), \; 0 < {r} < \rho \}. 
\end{align*}

\begin{lma} \label{lma46LP}
For $\alpha \in (0, \pi]$ and any bounded subdomain $\Omega$ of $\Sigma_{\alpha}$, the first mixed eigenvalue $\lambda_{1}(\Omega, \partial\Omega \cap \Sigma_{\alpha})$, defined in \eqref{CGY0104Mix}, satisfies
\begin{equation*}
\lambda_{1}(\Omega, \partial\Omega \cap \Sigma_{\alpha}) |\Omega| \geq \frac{1}{2} {j}_{0}^{2} \alpha, 
\end{equation*}
and the equality holds if and only if $\Omega = \Sigma_{\alpha, \rho}$ for $\rho > 0$. 
\end{lma}

\begin{proof}
The classical proof employs $\alpha$-symmetrization \cite{Ban80, PT85}. Under the assumption that $\alpha \leq \pi$, the domain $\Sigma_{\alpha}$ is convex. According to \cite[Theorem 1.1]{LP90}, one has 
\begin{equation*}
\sup_{{E}} \frac{|{E}|^{1/2}}{{P}_{\Sigma_{\alpha}}({E})} = \frac{1}{\sqrt{2\alpha}}, 
\end{equation*}
where the supremum is taken over all measurable subsets ${E}$ of $\Sigma_{\alpha}$. Moreover, this supremum is attained if and only if ${E}$ is the spherical sector $\Sigma_{\alpha, \rho}$. Here, $|{E}|$ denotes the Lebesgue measure of ${E}$, and ${P}_{\Sigma}({E})$ denotes the De Giorgi perimeter of ${E}$ relative to $\Sigma_{\alpha}$ (i.e., the measure of $\partial {E} \cap \Sigma_{\alpha}$; cf. \cite{PT85}). Let $\Omega$ be a bounded subdomain of $\Sigma_{\alpha}$, and define $\Omega^{\sharp} = \Sigma_{\alpha, \rho}$ with $\rho = \sqrt{2|\Omega|/\alpha}$. Let ${v}$ denote the positive eigenfunction associated with $\lambda_{1}(\Omega, \partial\Omega \cap \Sigma_{\alpha})$. By \cite[Proposition 1.2]{LPT88}, the $\alpha$-symmetrization ${v}^{\sharp} = {C}_{\alpha}{v}$ satisfies 
\begin{equation*}
\int_{\Omega^{\sharp}} |{v}^{\sharp}|^{2} dx = \int_{\Omega} |{v}|^{2} dx \text{ and } \int_{\Omega^{\sharp}} |\nabla {v}^{\sharp}|^{2} dx \leq \int_{\Omega} |\nabla {v}|^{2} dx. 
\end{equation*}
Owing to the spherical symmetry of ${v}^{\sharp}$, the Rayleigh quotient yields 
\begin{equation*}
\lambda_{1}(\Omega, \partial\Omega \cap \Sigma_{\alpha}) = \frac{\int_{\Omega} |\nabla {v}|^{2} dx}{\int_{\Omega} |{v}|^{2} dx} \geq
\frac{\int_{\Omega^{\sharp}} |\nabla {v}^{\sharp}|^{2} dx}{\int_{\Omega^{\sharp}} |{v}^{\sharp}|^{2} dx} \geq \left( \frac{{j}_{0}}{\rho} \right)^{2} = \frac{\alpha}{2|\Omega|} {j}_{0}^{2}. 
\end{equation*}
This concludes the proof.
\end{proof}

\begin{thm} \label{thm47}
Let ${T} = \triangle {z}_{1}{z}_{2}{z}_{3}$ be a triangle. Let ${z}_{4}$ be a point on the side $\overline{{z}_{3}{z}_{2}}$ such that $\overline{{z}_{1}{z}_{4}}$ is the internal bisector of the angle of ${T}$ at ${z}_{1}$. Suppose that $|\overline{{z}_{1}{z}_{2}}| \leq |\overline{{z}_{1}{z}_{3}}|$. 
Then
\begin{equation} \label{CGY0422}
\lambda_{1}(\triangle {z}_{1}{z}_{2}{z}_{4}, \overline{{z}_{1}{z}_{2}}) > \mu_{2}({T}), 
\end{equation}
where the left-hand side represents the first mixed eigenvalue as defined in \eqref{CGY0104Mix}. 
\end{thm}

\begin{proof}
Let $\alpha_{1}$, $\alpha_{2}$, and $\alpha_{3}$ denote the interior angles of triangle ${T}$ at vertices ${z}_{1}$, ${z}_{2}$, and ${z}_{3}$, respectively. The assumption $|\overline{{z}_{1}{z}_{2}}| \leq |\overline{{z}_{1}{z}_{3}}|$ implies $\alpha_{3} \leq \alpha_{2}$ and hence $\alpha_{3} < \pi/2$. We split the proof into two parts.

\textbf{Part 1}. 
\eqref{CGY0422} is valid when 
\begin{equation} 
\alpha_{3} \leq \alpha_{2} \leq (\pi+\alpha_{3})/3 \text{ and } 0 < \alpha_{3} < \pi/2. 
\end{equation}
In fact, since $\alpha_{3} \leq \alpha_{2}$, it follows from \autoref{lma44} that $\mu_{2}({T}) \leq \lambda_{1}({T}_{ - }, \overline{{z}_{1}{z}_{4}})$ with equality holding if and only if $\alpha_{3} = \alpha_{2}$. Note that $\alpha_{2} \leq (\pi+\alpha_{3})/3$ means $|\overline{{z}_{1}{z}_{4}}| \leq |\overline{{z}_{1}{z}_{2}}|$. By \autoref{thm43Siu}, $\lambda_{1}({T}_{ - }, \overline{{z}_{1}{z}_{4}}) < \lambda_{1}({T}_{ - }, \overline{{z}_{1}{z}_{2}})$ if and only if $|\overline{{z}_{1}{z}_{4}}| < |\overline{{z}_{1}{z}_{2}}|$. Combining these inequalities, we deduce that \eqref{CGY0422} holds when $|\overline{{z}_{1}{z}_{4}}| \leq |\overline{{z}_{1}{z}_{2}}| \leq |\overline{{z}_{1}{z}_{3}}|$, which corresponds precisely to $\alpha_{3} \leq \alpha_{2} \leq (\pi+\alpha_{3})/3$. This completes Part 1. 

\textbf{Part 2}. 
\eqref{CGY0422} is valid for 
\begin{equation} \label{CGY0426}
(\pi + \alpha_{3})/3 \leq \alpha_{2} < \pi - \alpha_{3} \text{ and } 0 < \alpha_{3} < \pi/2. 
\end{equation}
This part consists of a basic and tedious computation. Indeed, from \autoref{lma46LP} and \autoref{lma45}, 
\begin{equation} 
\lambda_{1}({T}_{ - }, \overline{{z}_{1}{z}_{2}}) > \frac{{j}_{0}^{2}}{|\overline{{z}_{4}{z}_{1}}| \cdot |\overline{{z}_{4}{z}_{2}}|} \frac{\angle {z}_{1}{z}_{4}{z}_{2}}{\sin \angle {z}_{1}{z}_{4}{z}_{2}} \text{ and } \mu_{2}({T}) < \Big(\frac{2{j}_{0}}{ |\overline{{z}_{3}{z}_{1}}|}\Big)^{2} . 
\end{equation}
It suffices to show that the ratio ${F} = {F}(\alpha_{2}, \alpha_{3})$ between the two right-hand sides of the above inequalities is greater than $1$. A direct computation shows that 
\begin{equation}
\begin{aligned}
{F}(\alpha_{2}, \alpha_{3}) & = \frac{|\overline{{z}_{3}{z}_{1}}|^{2}}{4|\overline{{z}_{4}{z}_{1}}| \cdot |\overline{{z}_{4}{z}_{2}}|} \cdot \frac{\angle {z}_{1}{z}_{4}{z}_{2}}{\sin \angle {z}_{1}{z}_{4}{z}_{2}}
\\ & = 
\frac{\cos^{2}(\alpha_{2}/2 - \alpha_{3}/2)\sin\alpha_{2}}{4\cos(\alpha_{2}/2 + \alpha_{3}/2)\sin^{2}\alpha_{3} } \cdot \frac{ (\pi - \alpha_{2} + \alpha_{3})/2}{ \sin ((\pi - \alpha_{2} + \alpha_{3})/2) }.
\end{aligned} 
\end{equation}

\textbf{Claim 1}: 
${F}(\alpha_{2}, \alpha_{3}) > 1$ when $\alpha_{2}=(\pi+\alpha_{3})/3$, that is, 
\begin{equation*} 
{F}(\alpha_{2}, \alpha_{3}) > 1 \text{ for } \alpha_{3} = 3\alpha_{2} - \pi \text{ and } \alpha_{2} \in (\pi/3, \pi/2). 
\end{equation*}

In fact, from the expression for ${F}$, we have 
\begin{equation*}
{F}(\alpha_{2}, 3\alpha_{2} - \pi) = \frac{\sin^{3}\alpha_{2}}{4\sin(2\alpha_{2})\sin^{2}(3\alpha_{2}) } \frac{ \alpha_{2}}{ \sin\alpha_{2} } = \frac{1}{8\cos\alpha_{2}(1 - 4\cos^{2}\alpha_{2})^{2} } \frac{ \alpha_{2}}{ \sin\alpha_{2} }. 
\end{equation*}
Since the function $\alpha_{2} \mapsto \alpha_{2}/\sin\alpha_{2}$ is increasing in $(0, \pi)$, it follows that 
\begin{equation*}
{F}(\alpha_{2}, 3\alpha_{2} - \pi) \geq \frac{25\sqrt{5}}{64}\frac{\pi/3}{\sin(\pi/3)} = \frac{25\sqrt{5}\pi}{96\sqrt{3}} > 1 \text{ for } \alpha_{2} \in (\frac{\pi}{3}, \frac{\pi}{2}). 
\end{equation*}

In order to simplify the expression ${F}$, we use the change of variables 
\begin{equation*}
{s} = \tan(\alpha_{2}/2) \text{ and } {t} = \tan(\alpha_{3}/2)
\end{equation*}
and denote ${G}({s}, {t}) = {F}(2\arctan {s}, 2\arctan {t})$. Then 
\begin{align} \label{CGY0434a}
{G}({s}, {t}) = \frac{(1 + {t}^{2})^{2}}{8{t}^{2}}\frac{(1 + {s}{t}){s}}{(1 - {s}{t})(1 + {s}^{2})} \Big( \frac{\pi}{2} - \arctan {s} + \arctan {t} \Big).
\end{align} 
Note that \eqref{CGY0426} implies that $(\pi + 2\arctan {t})/6 \leq {s} < 1/{t}$ and $0 < {t} < 1$.

\textbf{Claim 2}: 
${F}(\alpha_{2}, \alpha_{3}) > 1$ for $\alpha_{2}\geq \pi/2$. In fact, we can rewrite ${G}$ as follows: 
\begin{align} \label{CGY0434b}
{G}({s}, {t}) = 
\frac{(1 + {s}{t})^{2}}{8(1 - {s}{t})({s}{t})^{2}} \cdot \left( \frac{{s}^{2} + ({s}{t})^{2}}{{s}^{2} + 1} \right)^{3/2} \cdot \frac{ \arctan(1/{s}) + \arctan {t}}{ \sin ( \arctan(1/{s}) + \arctan {t}) }. 
\end{align}
From ${s} \geq 1$ and $0 < {t} < 1/{s}$, we have $({s}^{2} + {s}^{2}{t}^{2})/({s}^{2} + 1) \geq (1 + {s}^{2}{t}^{2})/2$. 
Substituting this into \eqref{CGY0434b}, we obtain 
\begin{equation*}
{G}({s}, {t}) > \frac{(1 + {s}{t})^{2}}{8(1 - {s}{t})({s}{t})^{2}} \cdot \bigg(\frac{1 + ({s}{t})^{2}}{2}\bigg)^{3/2} \cdot 1
= \sqrt{\frac{(1 + {s}{t})^{4} (1 + ({s}{t})^{2})^{3}}{512(1 - {s}{t})^{2}({s}{t})^{4}}}
> 1, 
\end{equation*}
where we have used the fact that for ${x} > 0$
\begin{equation*}
\begin{aligned}
(1 + {x})^{4}(1 + {x}^{2})^{3} - 512(1 - {x})^{2}{x}^{4}
=\, & (1 - {x})^{10} + 14{x}(1 - {x})^{8} + 76{x}^{2}(1 - {x})^{6} \\[1mm]
& + 72{x}^{3}(1 - {x})^{4} + 128{x}^{3}(1 - 2{x})^{2} > 0. 
\end{aligned}
\end{equation*}
This completes the proof of the assertion. 

\textbf{Claim 3}: 
For any $\alpha_{3} \in (0, \pi/2)$, we have 
\begin{equation} \label{CGY0436}
{F}(\alpha_{2}, \alpha_{3})
\geq \min\bigl\{{F}(\tfrac{\pi+\alpha_{3}}{3}, \alpha_{3}), {F}(\tfrac{\pi}{2}, \alpha_{3})\bigr\}
\text{ for } \alpha_{2}\in \bigl[\tfrac{\pi+\alpha_{3}}{3}, \tfrac{\pi}{2}\bigr].
\end{equation}
Indeed, this can be shown by proving that for each fixed ${t} \in (0, 1)$, the function 
\begin{equation} \label{CGY0437}
{s} \in [0, 1] \mapsto {G}({s}, {t})
\end{equation}
is either strictly increasing or unimodal. From the expression of ${G}$ in \eqref{CGY0434a}, one can compute that the derivatives of ${G}$ with respect to ${s}$ satisfy 
\begin{equation} \label{CGY0438a}
\partial_{{s}{s}}{G}({s}, {t}) = - {H}_{1}({s}, {t}) {H}_{2}({s}, {t}) + {H}_{3}({s}, {t})\partial_{{s}}{G}({s}, {t}), 
\end{equation}
where the functions ${H}_{i}$ ($i = 1, 2, 3$), defined for ${s} \in [0, 1]$ and ${t} \in (0, 1)$, are given by 
\begin{equation*}
\begin{aligned}
{H}_{1}({s}, {t}) & = {c}_{0} + {c}_{1} {s}^{2} + {c}_{2} {s}^{4}, 
\\
{H}_{2}({s}, {t}) & = 
\frac{(1 + {t}^{2})^{2}}{4{t}^{2}(1 - {s}{t})(1 + {s}^{2})^{3}[ 2{s}{t}(1 + {s}^{2}) + (1 - {s}^{2})(1 - {s}^{2}{t}^{2})]}, 
\\
{H}_{3}({s}, {t}) & = 
\partial_{{s}} \Bigl(\ln \frac{ 2{s}{t}(1 + {s}^{2}) + (1 - {s}^{2})(1 - {s}^{2}{t}^{2})}{ (1 - {s}{t})^{2}(1 + {s}^{2})^{2} } \Bigr). 
\end{aligned}
\end{equation*}
Here, ${x} = {s}{t} \in [0, 1)$, and 
\begin{equation*}
{c}_{0} = ({x} + 1)^{3} - 2{x}^{3} > 0, \quad {c}_{1} = - 2{x}(1 - {x}) \leq 0, \text{ and } {c}_{2} = 2 - ({x} + 1)^{3}.
\end{equation*}
Moreover, ${H}_{1}$ and ${H}_{2}$ satisfy
\begin{equation} \label{CGY0438b}
{H}_{1}({s}, {t}) > 0 \text{ and } {H}_{2}({s}, {t}) > 0 \text{ for } {s} \in [0, 1] \text{ and } {t} \in (0, 1). 
\end{equation}
The positivity of ${H}_{2}$ is immediate. To prove the positivity of ${H}_{1}$, define 
\begin{equation*}
\widetilde{H}_{1}({s}, {x}) = \bigl( ({x} + 1)^{3} - 2{x}^{3}\bigr) - 2{x}(1 - {x}) {s}^{2} + \bigl( 2 - ({x} + 1)^{3} \bigr) {s}^{4}. 
\end{equation*}
Note that ${H}_{1}({s}, {t}) = \widetilde{H}_{1}({s}, {s}{t})$. Observe also that $\partial_{{s}}\widetilde{H}_{1}(1, {x}) = 4(1 - 4{x} - 2{x}^{2} - {x}^{3})$ and
\begin{equation*}
\widetilde{H}_{1}(0, {x}) = ({x} + 1)^{3} - 2{x}^{3} > 0, \quad \widetilde{H}_{1}(1, {x}) = 2(1 - {x})(1 + {x}^{2}) > 0. 
\end{equation*}
In the case where $1 - 4{x} - 2{x}^{2} - {x}^{3} \leq 0$, the function ${s} \in [0, 1] \mapsto \widetilde{H}_{1}({s}, {x})$ is strictly decreasing, and hence $\widetilde{H}_{1}({s}, {x}) > 0$. Conversely, if
$1 - 4{x} - 2{x}^{2} - {x}^{3} > 0$, then one deduces that ${x} < 1/4$ and $({x}^{2} + 1)/{x} > 17/4$. Consequently, 
\begin{equation*}
4{c}_{0}{c}_{2} - {c}_{1}^{2} = 4({x} + 1)^{2}[{x}^{2} + 1 + (\sqrt{10} - 1){x}][{x}^{2} + 1 - (\sqrt{10} + 1){x}] > 0, 
\end{equation*}
which implies $\widetilde{H}_{1}({s}, {x}) > 0$. This verifies \eqref{CGY0438b}. 

As a direct consequence of \eqref{CGY0438a} and \eqref{CGY0438b}, every critical point of the function \eqref{CGY0437} is a strict local maximum. Furthermore, combining this with $\partial_{{s}}{G}(0, {t}) > 0$, we deduce that for any ${t} \in (0, 1)$, the function \eqref{CGY0437} is either strictly increasing or increasing-then-decreasing (i.e., it first increases and then decreases). This completes the proof of the assertion.

By combining these three assertions, we establish Part 2, thereby completing the proof.
\end{proof}


\section{The nodal line of the second Neumann eigenfunction} \label{Sect5nodal}

In this section, we prove the rigidity result that the second Neumann eigenfunction does not vanish at any vertex for non-isosceles triangles, using the eigenvalue inequality \autoref{thm47} established in \autoref{Sect4EI}. 

Following \autoref{Sect3mon}, we also assume that 
\begin{equation} \label{CGY0502a}
{u} \text{ is monotone in a unit direction } \mathbf{a} \in \R^{2}. 
\end{equation} 
This monotonicity condition holds for non-obtuse triangles (\autoref{prop28obt}). Furthermore, this monotonicity condition will be established for acute triangles by the continuity method via domain deformation, as detailed in \autoref{Sect6unique} below.

\begin{thm} \label{thm51}
Let ${u}$ be a second Neumann eigenfunction of a triangle ${T}$ with vertices ${z}_{1}$, ${z}_{2}$ and ${z}_{3}$.
Suppose that ${u}$ satisfies \eqref{CGY0502a} and 
\begin{equation} \label{CGY0502b}
{u}({z}_{1}) = 0. 
\end{equation}
Then ${T}$ is an isosceles triangle with $|\overline{{z}_{3}{z}_{1}}| = |\overline{{z}_{1}{z}_{2}}|$ and ${u}$ is antisymmetric with respect to the internal bisector of the angle of ${T}$ at the vertex ${z}_{1}$. Moreover, $\alpha_{1} = \angle {z}_{3}{z}_{1}{z}_{2} \geq \pi/3$. 
\end{thm}

\begin{proof}
Without loss of generality, we assume that
\begin{equation*}
|\overline{{z}_{3}{z}_{1}}| \geq |\overline{{z}_{1}{z}_{2}}|. 
\end{equation*}
Let $\overline{{z}_{1}{z}_{4}}$ be the internal bisector of the angle of triangle ${z}_{1}{z}_{2}{z}_{3}$ at vertex ${z}_{1}$, where ${z}_{4}$ lies on the side $\overline{{z}_{3}{z}_{2}}$. Let ${z}_{5}$ be the reflection of ${z}_{2}$ across the internal bisector $\overline{{z}_{1}{z}_{4}}$. Let ${T}_{ - }$ denote the triangle with vertices ${z}_{1}$, ${z}_{2}$, and ${z}_{4}$. The assumption $|\overline{{z}_{3}{z}_{1}}| \geq |\overline{{z}_{1}{z}_{2}}|$ implies that ${z}_{5}$ is either ${z}_{3}$ or lies in the interior of the line segment $\overline{{z}_{3}{z}_{1}}$; see \autoref{fig5}. By rotating and translating the coordinate system, if necessary, we may assume that ${z}_{1}$ is at the origin and the internal bisector $\overline{{z}_{1}{z}_{4}}$ lies along the positive ${x}_{1}$-axis; see \autoref{fig5}. 

By \autoref{lma37}, after a possible sign change, the assumptions \eqref{CGY0502a} and \eqref{CGY0502b} become ${u}({z}_{1}) = 0$ and
\begin{equation*} 
\nabla {u} \cdot \boldsymbol{n}_{S} > 0 \text{ in } {T}.
\end{equation*}
By the same argument as in \autoref{lma31}, it follows that
\begin{equation} 
{R}_{{p}} {u} < 0 \text{ in } \{{x} \in {T}: \, ({x}- {p}) \cdot \boldsymbol{n}_{S} > 0\} \text{ with ${p} = {z}_{4}$.}
\end{equation}
Since $\angle {z}_{3}{z}_{4}{z}_{5} \leq \angle {z}_{3}{z}_{2}{z}_{1}$, it follows in particular that ${u}$ satisfies
\begin{equation} \label{CGY0502d}
\nabla {u} \cdot \nu_{\overline{{z}_{4}{z}_{5}}} \geq 0 \text{ on } \operatorname{Int}(\overline{{z}_{4}{z}_{5}}), 
\end{equation}
and the inequality is strict whenever $|\overline{{z}_{3}{z}_{1}}| > |\overline{{z}_{1}{z}_{2}}|$. Here $\nu_{\overline{{z}_{4}{z}_{5}}}$ denotes the outward-pointing unit normal to the segment $\overline{{z}_{4}{z}_{5}}$ with respect to the kite ${z}_{1}{z}_{2}{z}_{4}{z}_{5}$.
Let $\angle {z}_{3}{z}_{1}{z}_{2} = \alpha_{1} \in (0, \pi)$. Consider the function
\begin{equation} \label{CGY0504}
{w}({x}) = {u}({x}_{1}, {x}_{2}) + {u}({x}_{1}, - {x}_{2})
\end{equation}
which is well-defined in the kite ${K}$ with vertices ${z}_{1}$, ${z}_{2}$, ${z}_{4}$ and ${z}_{5}$. Since ${u}$ is a Neumann eigenfunction satisfying \eqref{CGY0502d}, we derive that ${w}$ satisfies
\begin{equation} \label{CGY0506a} \begin{cases}
\Delta {w} + \mu {w} = 0 \text{ in } {K}, \\
\partial_{\nu}{w} = 0 \text{ on } \overline{{z}_{1}{z}_{2}} \cup \overline{{z}_{1}{z}_{5}} \subset \partial {K}, \\
\partial_{\nu}{w} \geq 0 \text{ on } \overline{{z}_{4}{z}_{2}} \cup \overline{{z}_{4}{z}_{5}} \subset \partial {K}
\end{cases}\end{equation}
where the equality in the boundary condition on $\overline{{z}_{4}{z}_{2}} \cup \overline{{z}_{4}{z}_{5}}$ occurs if and only if $\overline{{z}_{4}{z}_{5}}$ lies on $\overline{{z}_{3}{z}_{2}}$, i.e., $|\overline{{z}_{3}{z}_{1}}| = |\overline{{z}_{1}{z}_{2}}|$. 

The main objective is to show that ${w} \equiv 0$. We argue by contradiction and suppose that ${w}\not\equiv0$ in the kite ${K}$. Since both ${u}$ and ${w}$ are analytic, it follows that ${w} \not\equiv 0$ in any subdomain of ${K}$. Using the fact that ${w}({z}_{1}) = 0$ and the local analysis in a neighborhood of ${z}_{1} = O$ (cf. \cite[Lemma 4.1]{JM20}), the nodal line $\mathcal{Z}({w})$ consists of exactly $2{n}$ local curves $\Gamma_{\pm {j}}$ near ${z}_{1}$, and these curves $\Gamma_{\pm {j}}$ are tangent to the lines $\{{z}: \text{arg}({z}) = \pm \frac{2{j} - 1}{4n}\alpha_{1}\}$ at ${z}_{1}$, for ${j} = 1, 2, \ldots, {n}$. Let $\Gamma_{ - ({n} + 1)}$ and $\Gamma_{ + ({n} + 1)}$ denote the line segments $\overline{{z}_{1}{z}_{2}}$ and $\overline{{z}_{1}{z}_{5}}$, respectively. For $1 \leq {j} \leq {n}$, denote by ${D}_{\pm {j}}$, ${D}_{0}$ the connected components of the nodal domain of ${w}$ such that, in a small neighborhood $\mathcal{O}_{{z}_{1}}$ of ${z}_{1}$, ${D}_{\pm {j}} \cap \mathcal{O}_{{z}_{1}}$ is enclosed by $\Gamma_{\pm {j}}$, $\Gamma_{\pm({j} + 1)}$ and $\partial \mathcal{O}_{{z}_{1}}$, while ${D}_{0} \cap \mathcal{O}_{{z}_{1}}$ is enclosed by $\Gamma_{ - 1}$, $\Gamma_{ + 1}$ and $\partial \mathcal{O}_{{z}_{1}}$; see the case ${n} = 1$ in \autoref{fig5}. Consequently, ${w}$ has different signs in ${D}_{j}$ and ${D}_{{j} + 1}$. Let $\tilde{\Gamma}_{j}$ denote the maximal $C^{1}$ curve of $\mathcal{Z}({w})$ that contains $\Gamma_{j}$.  

\begin{figure}[h]\centering  \vspace*{-2ex}
\begin{tikzpicture}[scale = 2.5]
\pgfmathsetmacro\LenShort{1.0000};      
\pgfmathsetmacro\LenMid{1.3*\LenShort}; 
\pgfmathsetmacro\AngleLarge{80};        
\pgfmathsetmacro\Lenll{2*cos(\AngleLarge/2)*\LenShort*\LenMid/(\LenShort+\LenMid)}; 

\draw[] (0, 0)node[left]{${z}_{1}$} -- 
({-\AngleLarge/2}:\LenShort) node[above = 2pt, right]{${z}_{2}$} -- ({\AngleLarge/2}:\LenMid) node[below = 2pt, right]{${z}_{3}$} -- cycle; 
\draw[] ({\AngleLarge/2}:\LenShort) node[above]{${z}_{5}$} -- (\Lenll, 0) node[right]{${z}_{4}$};  
\fill[green] (0, 0) .. controls ({-\AngleLarge/4}:{\LenShort*0.2}) and ({-\AngleLarge/5}:{\LenShort*0.9}) .. ({-\AngleLarge/2}:{\LenShort*0.9}) -- cycle;  
\fill[green] (0, 0) .. controls ({\AngleLarge/4}:{\LenShort*0.2}) and ({\AngleLarge/5}:{\LenShort*0.9}) .. ({\AngleLarge/2}:{\LenShort*0.9}) -- cycle;   
\draw[very thin] (0, 0) -- (\Lenll, 0);  
\node at ({+\AngleLarge*0.37}:{\LenShort*0.6}) {\footnotesize ${D}_{+1}$}; 
\node at ({-\AngleLarge*0.37}:{\LenShort*0.6}) {\footnotesize ${D}_{-1}$}; 
\end{tikzpicture} \hspace{2ex}
\begin{tikzpicture}[scale = 2.5]
\pgfmathsetmacro\AngleLarge{80}; 
\pgfmathsetmacro\LenShort{1.0000}; 
\pgfmathsetmacro\LenMid{1.3*\LenShort}; 
\pgfmathsetmacro\xA{\LenMid*cos(\AngleLarge/2)}; 
\pgfmathsetmacro\yA{\LenMid*sin(\AngleLarge/2)}; 
\pgfmathsetmacro\xC{\LenShort*cos(\AngleLarge/2)}; 
\pgfmathsetmacro\yC{-\LenShort*sin(\AngleLarge/2)}; 
\pgfmathsetmacro\Lenll{2*cos(\AngleLarge/2)*\LenShort*\LenMid/(\LenShort+\LenMid)}; 
\pgfmathsetmacro\xiaoshu{0.4*\LenShort/(\LenShort+\LenMid)}; 
\pgfmathsetmacro\xR{\xC+\xiaoshu*(\xA-\xC)}; 
\pgfmathsetmacro\yR{\yC+\xiaoshu*(\yA-\yC)}; 
\fill[green] (0, 0) 
.. controls ({-\AngleLarge/4}:{\LenShort*0.2}) and ({-\AngleLarge/5}:{\LenShort*0.9}) .. (\xR, \yR) -- (\xC, \yC) -- cycle; 
\fill[green] (0, 0) .. controls ({\AngleLarge/4}:{\LenShort*0.2}) and ({\AngleLarge/5}:{\LenShort*0.9}) .. (\xR, -\yR) -- (\xC, -\yC) -- cycle; 
\draw[] (0, 0) node[left]{${z}_{1}$} -- ({-\AngleLarge/2}:\LenShort) node[above = 2pt, right]{${z}_{2}$} -- ({\AngleLarge/2}:\LenMid) node[below = 2pt, right]{${z}_{3}$} -- cycle; 
\draw[] ({\AngleLarge/2}:\LenShort) node[above]{${z}_{5}$} -- (\Lenll, 0) node[right]{${z}_{4}$}; 
\draw[very thin] (0, 0) -- (\Lenll, 0); 
\node at ({+\AngleLarge*0.36}:{\LenShort*0.65}) {\footnotesize ${D}_{+1}$}; 
\node at ({-\AngleLarge*0.36}:{\LenShort*0.65}) {\footnotesize ${D}_{-1}$}; 
\end{tikzpicture} \hspace{2ex}
\begin{tikzpicture}[scale = 2.5]
\pgfmathsetmacro\AngleLarge{80}; 
\pgfmathsetmacro\LenShort{1.0000}; 
\pgfmathsetmacro\LenMid{1.3*\LenShort}; 
\pgfmathsetmacro\xA{\LenMid*cos(\AngleLarge/2)}; 
\pgfmathsetmacro\yA{\LenMid*sin(\AngleLarge/2)}; 
\pgfmathsetmacro\xC{\LenShort*cos(\AngleLarge/2)}; 
\pgfmathsetmacro\yC{-\LenShort*sin(\AngleLarge/2)}; 
\pgfmathsetmacro\Lenll{2*cos(\AngleLarge/2)*\LenShort*\LenMid/(\LenShort+\LenMid)}; 
\pgfmathsetmacro\xiaoshu{0.4*\LenShort/(\LenShort+\LenMid)}; 
\pgfmathsetmacro\xR{\xC+\xiaoshu*(\xA-\xC)}; 
\pgfmathsetmacro\yR{\yC+\xiaoshu*(\yA-\yC)}; 
\fill[green] (0, 0) .. controls ({-\AngleLarge/4}:{\LenShort*0.2}) and ({-\AngleLarge/5}:{\LenShort*0.9}) .. (\xR, \yR) -- (\Lenll, 0) -- cycle; 
\fill[green] (0, 0) .. controls ({\AngleLarge/4}:{\LenShort*0.2}) and ({\AngleLarge/5}:{\LenShort*0.9}) .. (\xR, -\yR) -- (\Lenll, 0) -- cycle; 
\draw[] (0, 0) node[left]{${z}_{1}$} -- 
({-\AngleLarge/2}:\LenShort) node[above = 2pt, right]{${z}_{2}$} -- ({\AngleLarge/2}:\LenMid) node[below = 2pt, right]{${z}_{3}$} -- cycle; 
\draw[] ({\AngleLarge/2}:\LenShort) node[above]{${z}_{5}$} -- (\Lenll, 0) node[right]{${z}_{4}$}; 
\draw[very thin] (0, 0) -- (\Lenll, 0); 
\node at ({0}:{\LenShort*0.65}) {${D}_{0}$}; 
\end{tikzpicture} \vspace*{-2ex}
\caption{Nodal domain of $w(x) = u(x_{1}, x_{2}) + u(x_{1}, -x_{2})$}
\label{fig5}
\end{figure}
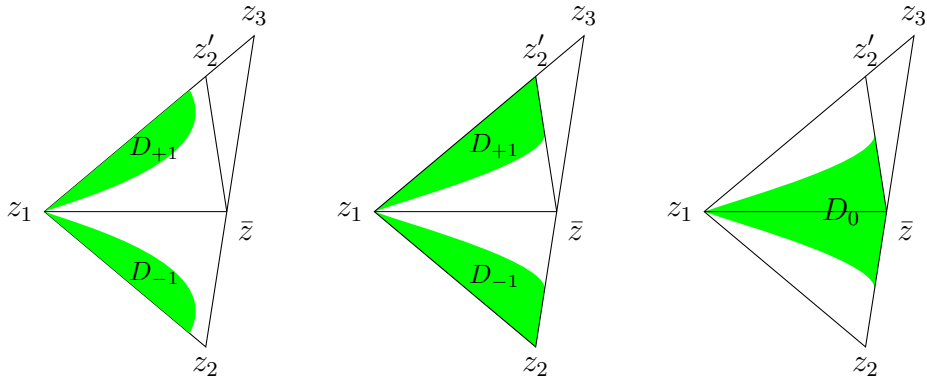

\textbf{Claim 1}. 
$\mathcal{Z}({w}) \cup \overline{{z}_{1}{z}_{2}}$ does not contain a loop. As a direct consequence of this claim and the symmetry of $\tilde{\Gamma}_{{j}}$ and $\tilde{\Gamma}_{-{j}}$ with respect to the internal bisector $\overline{{z}_{1}{z}_{4}}$, we conclude that each nodal line $\tilde{\Gamma}_{{j}}$ intersects $\overline{{z}_{1}{z}_{4}} \cup \overline{{z}_{1}{z}_{2}}$ only at ${z}_{1}$.

To prove the claim, we argue by contradiction. Suppose that $\mathcal{Z}({w}) \cup \overline{{z}_{1}{z}_{2}}$ contains a loop. Then there exists a nodal domain ${D}$ of ${w}$ with $\partial {D} \subset \mathcal{Z}({w}) \cup \overline{{z}_{1}{z}_{2}}$ (see the shadow domain ${D}_{-1}$ in the left panel of \autoref{fig5}). Inside ${D}$ the function ${w}$ satisfies
\begin{equation*}
\Delta {w} + \mu {w} = 0 \text{ in } {D}, \quad
{w} = 0 \text{ on } \partial {D} \setminus \overline{{z}_{1}{z}_{2}}, \quad
\partial_{\nu}{w} = 0 \text{ on } \partial {D} \cap \overline{{z}_{1}{z}_{2}}.
\end{equation*}
Hence the first mixed eigenvalue $\lambda_{1}({D}, \partial {D} \setminus \overline{{z}_{1}{z}_{2}})$ equals $\mu$. By the variational characterization, one has $\lambda_{1}({D}, \partial {D} \setminus \overline{{z}_{1}{z}_{2}}) > \lambda_{1}({T}, \partial {T} \setminus \overline{{z}_{1}{z}_{2}})$.
However, the mixed-Neumann comparison (see \autoref{lma25LR}) yields $\lambda_{1}({T}, \partial {T} \setminus \overline{{z}_{1}{z}_{2}}) \ge \mu_{2}({T}) = \mu$, which contradicts the strict inequality above. Therefore, $\mathcal{Z}({w}) \cup \overline{{z}_{1}{z}_{2}}$ contains no loop, and the claim follows.

\textbf{Claim 2}: 
The function ${w}$ cannot be positive in ${D}_{0}$. Indeed, we suppose that ${w} > 0$ in ${D}_{0}$ and then ${w} < 0$ in ${D}_{ - 1}$. From \eqref{CGY0506a}, ${w}$ satisfies
\begin{equation} \label{CGY0506b} 
- \Delta {w} = \mu {w} < 0 \text{ in } {D}_{ - }, \quad
\partial_{\nu}{w} \geq 0 \text{ on } \Gamma_{*}, \quad
\partial_{\nu}{w} = 0 \text{ on } \partial {D}_{ - } \setminus \Gamma_{*}, 
\end{equation}
where we set ${D}_{ - } = {D}_{ - 1}$ and $\Gamma_{*} = \partial {D}_{ - } \setminus \partial {K}$; see the shadow domain ${D}_{ - 1}$ in the middle panel of \autoref{fig5}. By \autoref{lma22cMP}, we obtain
\begin{equation*}
\lambda_{1}({D}_{ - }, \Gamma_{*}) \leq \mu. 
\end{equation*}
From Claim 1, ${D}_{ - }$ is contained within ${T}_{ - }$ and $\partial {D}_{ - }$ intersects $\overline{{z}_{1}{z}_{4}}$ only at ${z}_{1}$. By the variational characterization of the first eigenvalue, we obtain
\begin{equation} \label{CGY0507}\begin{aligned}
\lambda_{1}({D}_{ - }, \Gamma_{*})
& = \inf\bigg\{\frac{\int_{{D}_{ - }} |\nabla \varphi|^{2} dx}{\int_{{D}_{ - }} |\varphi|^{2} dx}: \, 0 \neq \varphi\in W^{1, 2}({D}_{ - }), \varphi = 0 \text{ on } \Gamma_{*}\bigg\}
\\ & = 
\inf\bigg\{\frac{\int_{{T}_{ - }} |\nabla \varphi|^{2} dx}{\int_{{T}_{ - }} |\varphi|^{2} dx}: \, 0 \neq \varphi\in W^{1, 2}({T}_{ - }), \varphi = 0 \text{ on } {T}_{ - } \setminus {D}_{ - }\bigg\}
\\ & > 
\inf\bigg\{\frac{\int_{{T}_{ - }} |\nabla \varphi|^{2} dx}{\int_{{T}_{ - }} |\varphi|^{2} dx}: \, 0 \neq \varphi\in W^{1, 2}({T}_{ - }) \text{ and } \varphi = 0 \text{ on } \overline{{z}_{1}{z}_{4}}\bigg\}
\\ & = 
\lambda_{1}({T}_{ - }, \overline{{z}_{1}{z}_{4}}). 
\end{aligned}
\end{equation}
However, from \autoref{lma44}, we have
\begin{equation*}
\lambda_{1}({T}_{ - }, \overline{{z}_{1}{z}_{4}}) \geq \mu_{2}({T}) = \mu. 
\end{equation*}
This leads to a contradiction. Hence, ${w}$ cannot be positive in ${D}_{0}$. 

\textbf{Claim 3}: 
The function ${w}$ cannot be negative in ${D}_{0}$. Indeed, argue by contradiction and suppose that ${w} < 0$ in ${D}_{0}$. By Claim~1, the boundary of ${D}_{0}$ intersects $\overline{{z}_{1}{z}_{2}} \cup \overline{{z}_{1}{z}_{5}}$ only at ${z}_{1}$ (see the shadow domain ${D}_{0}$ in the right panel of \autoref{fig5}). Set ${D}_{-} = {D}_{0}$ and $\Gamma_{*} = \partial {D}_{ - } \setminus \partial {K} = \partial {D}_{ - } \setminus (\overline{{z}_{2}{z}_{4}} \cup \overline{{z}_{5}{z}_{4}})$. By \eqref{CGY0506a}, ${w}$ satisfies \eqref{CGY0506b}. Then, by \autoref{lma22cMP},  
\begin{equation*}
\lambda_{1}({D}_{ - }, \Gamma_{*}) < \mu.
\end{equation*}
By the same reasoning as in \eqref{CGY0507}, we have
\begin{equation*} 
\lambda_{1}({D}_{ - }, \Gamma_{*}) > \lambda_{1}({K}, \overline{{z}_{1}{z}_{2}} \cup \overline{{z}_{1}{z}_{5}}) = \lambda_{1}( {T}_{ - }, \overline{{z}_{1}{z}_{2}}), 
\end{equation*}
where the last equality follows from the symmetry of the kite ${K}$. Consequently, 
\begin{equation}
\lambda_{1}({T}_{ - }, \overline{{z}_{1}{z}_{2}}) < \mu_{2}({T}) = \mu, 
\end{equation}
which contradicts the inequality \eqref{CGY0422} established in \autoref{thm47}. Therefore, ${w}$ cannot be negative in ${D}_{0}$. 

Having established the three assertions above, we reach a contradiction, and hence ${w}\equiv0$ in the kite ${K}$. Consequently, $|\overline{{z}_{3}{z}_{1}}| = |\overline{{z}_{1}{z}_{2}}|$ and ${u}$ is antisymmetric with respect to the internal bisector $\overline{{z}_{1}{z}_{4}}$. 
Moreover, by \autoref{prop27iso}, the second eigenfunction ${u}$ is symmetric with respect to the bisector $\overline{{z}_{1}{z}_{4}}$ if $\angle {z}_{3}{z}_{1}{z}_{2} < \pi/3$. This implies that $\angle {z}_{3}{z}_{1}{z}_{2} < \pi/3$ cannot occur, and hence $\angle {z}_{3}{z}_{1}{z}_{2} \geq \pi/3$. This completes the proof. 
\end{proof}

The direct study of the function ${w}$ in \eqref{CGY0504} is similar in spirit to the proof of symmetry in \cite[Theorem 4.3]{GM18}. In \cite[Theorem 4.3]{GM18}, the second author and Moradifam use the sphere covering inequality (a general and different form of an eigenvalue inequality) to derive a contradiction, whereas here we reach a contradiction by applying specific eigenvalue inequalities in different domains. 


\section{Proof of uniqueness of critical point} \label{Sect6unique}

In this section, we establish the main result of \autoref{thm12}, excluding the exact locations of the global extrema (i.e., Property \ref{CGY0102ite} of \autoref{thm12}). Notably, the proof relies on \autoref{thm38}, \autoref{thm47}, and \autoref{thm51} from previous sections. 

For simplicity, throughout the rest of the paper, we assume that the triangle ${T}$ with vertices ${z}_{1}$, ${z}_{2}$, ${z}_{3}$ satisfies
\begin{equation} \label{CGY0601a}
{z}_{1} = 0 \in \mathbb{C}, \quad {z}_{2} = 1 \in \mathbb{C} \text{ and } {z}_{3} \in \{{z} \in \mathbb{C}: \, \operatorname{Im}(z) > 0\}
\end{equation}
and $|{z}_{1} - {z}_{2}| \leq |{z}_{3} - {z}_{1}| \leq |{z}_{2} - {z}_{3}|$. 
We remark that in this section we only require 
\begin{equation} \label{CGY0601c}
|{z}_{1} - {z}_{2}| < \min\{|{z}_{2} - {z}_{3}|, |{z}_{3} - {z}_{1}|\}, 
\end{equation}
while the condition that $\overline{{z}_{2}{z}_{3}}$ is the longest side is used later in \autoref{Sect7location}. 

The main theorem is as follows. 

\begin{thm} \label{thm42}
Let ${T}$ be an acute triangle with vertices ${z}_{1}$, ${z}_{2}$, and ${z}_{3}$ satisfying conditions \eqref{CGY0601a} and \eqref{CGY0601c}. Let ${u}$ be an eigenfunction corresponding to the second Neumann eigenvalue $\mu$ of ${T}$. Then ${u}$ has exactly one non-vertex critical point. More precisely, ${u}$ has the following properties. 
\begin{enumerate}[label = \rm(\arabic*), start = 1]
\item \label{CGY0601itema}
The second Neumann eigenfunction is unique (up to multiplication by a constant). 
\item \label{CGY0601itemc}
After possibly multiplying by $ - 1$, we have ${u}({z}_{1}) < 0$, ${u}({z}_{2}) < 0$ and ${u}({z}_{3}) > 0$. 
\item \label{CGY0601itemd}
As in item \ref{CGY0601itemc}, the derivative of ${u}$ with respect to ${x}_{2}$ satisfies $\partial_{{x}_{2}}{u} > 0$ in ${T}$. 
\item \label{CGY0601iteme}
The eigenfunction ${u}$ has exactly one non-vertex critical point, which lies in $\operatorname{Int}(\overline{{z}_{1}{z}_{2}})$ and is a saddle point. 
\end{enumerate}
\end{thm}

The proof relies on the continuity method via domain deformation; see \cite{JN00} for an overview of the method. Let ${T} = {T}^{1}$ be any fixed (target) acute triangle with labeled vertices $({z}_{1}^{1}, {z}_{2}^{1}, {z}_{3}^{1})$ satisfying \eqref{CGY0601a} and \eqref{CGY0601c}. Then ${T}^{1}$ is neither an equilateral triangle nor a superequilateral triangle. 

In order to ensure that the eigenfunctions do not vanish at every vertex, we choose the half-equilateral triangle as the initial triangle, i.e., ${T}^{0}$ with labeled vertices $({z}_{1}^{0}, {z}_{2}^{0}, {z}_{3}^{0})$ defined by
\begin{equation} \label{CGY0602a}
{z}_{1}^{0} = 0, \quad {z}_{2}^{0} = 1, \text{ and } {z}_{3}^{0} = \sqrt{ - 3} \in \mathbb{C}. 
\end{equation}
For the half equilateral triangle ${T}^{0}$, the second Neumann eigenvalue is simple and equals $\mu_{0} = 4\pi^{2}/9$, and the corresponding (normalized) eigenfunction is
\begin{equation} \label{CGY0602b}
{u}^{0}({x}_{1}, {x}_{2}) = (16/27)^{1/4} \cdot \Big( \cos\frac{2\pi {x}_{1}}{3} - 2 \cos\frac{\pi {x}_{2}}{\sqrt{3}} \cos\frac{\pi {x}_{1}}{3} \Big). 
\end{equation}

Next, we define the triangle ${T}^{t}$ with labeled vertices $({z}_{1}^{t}, {z}_{2}^{t}, {z}_{3}^{t})$ by
\begin{equation} \label{CGY0603a}
({z}_{1}^{t}, {z}_{2}^{t}, {z}_{3}^{t}) = (1 - {t}) \cdot ({z}_{1}^{0}, {z}_{2}^{0}, {z}_{3}^{0}) + {t} \cdot ({z}_{1}^{1}, {z}_{2}^{1}, {z}_{3}^{1}), 
\end{equation}
for ${t} \in [0, 1]$. Then ${T}^{t}$, ${t} \in (0, 1]$, is a continuous family of acute triangles that joins our target acute triangle ${T}^{1}$ to the initial triangle ${T}^{0}$. Moreover, one can verify that
\begin{equation} \label{CGY0603b}
|{z}_{1}^{t} - {z}_{2}^{t}| < \min \{|{z}_{2}^{t} - {z}_{3}^{t}|, |{z}_{3}^{t} - {z}_{1}^{t}|\}
\end{equation}
for every ${t} \in [0, 1]$. 
Let ${u}^{t}$ be a second eigenfunction associated with the second Neumann eigenvalue $\mu_{t}$ of ${T}^{t}$, and assume that ${u}^{t}$ is normalized, i.e., 
\begin{equation*}
\int_{{T}^{t}} |{u}^{t}|^{2} dx = 1. 
\end{equation*}
Let ${h}_{t}$ be the unique real affine homeomorphism that maps the ordered triple $({z}_{1}^{0}, {z}_{2}^{0}, {z}_{3}^{0})$ to $({z}_{1}^{t}, {z}_{2}^{t}, {z}_{3}^{t})$. It is clear that the eigenvalue $\mu_{t}$ varies continuously with ${t}$ (see, e.g., \cite[Theorem 2.3.25]{Hen06}). 
Following \cite[Section 9]{JM20}, standard perturbation theory implies that ${u}^{t}$ depends continuously on ${t}$ provided that $\mu_{t}$ is simple. More precisely, if $\mu_{\bar{t}}$ is simple for some $\bar{t}$, then there exists an open interval $\mathcal{I}$ containing $\bar{t}$ such that: (i) $\mu_{t}$ remains simple for all ${t} \in \mathcal{I}$, (ii) after possibly multiplying by $ - 1$, the map ${t} \mapsto {u}^{t} \circ {h}_{t}$ is continuous from $\mathcal{I}$ to $C(\overline{{T}^{0}}) \cap C^{k}_{loc}(\overline{{T}^{0}} \setminus \{\text{vertices}\})$ for every $k \in \mathbb{N}$. 

We first recall the number of critical points near a vertex. 

\begin{lma} \label{lma62}
Assume that the second Neumann eigenvalue $\mu_{\bar{t}}$ of the triangle ${T}^{\bar{t}}$ (acute or not) is simple for some $\bar{t} \in [0, 1]$. Then there exist two small positive constants $\bar{\delta} > 0$ and $\bar{\varepsilon} > 0$, depending on $\bar{t}$, such that for any ${j} \in \{1, 2, 3\}$ and $|{t} - \bar{t}| < \bar{\varepsilon}$ (with ${t} \in [0, 1]$), the following conditions are satisfied: 
\begin{enumerate}[label = \rm(\roman*), start = 1]
\item
$\operatorname{crit}_{\mathrm{nv}}({u}^{t}) \cap \operatorname{Int}({T}^{t}) \cap \mathcal{O}_{\bar{\delta}}({z}_{j}^{t}) = \emptyset$. 
\item
$\operatorname{crit}_{\mathrm{nv}}({u}^{t}) \cap \mathcal{O}_{\bar{\delta}}({z}_{j}^{t})$ is contained only on one side of ${T}^{t}$. 
\item
If ${u}^{\bar{t}}({z}_{j}^{\bar{t}}) \neq 0$ and the interior angle of ${T}^{\bar{t}}$ at ${z}_{j}^{\bar{t}}$ is less than $\pi/2$, then
$\operatorname{crit}_{\mathrm{nv}}({u}^{t}) \cap \mathcal{O}_{\bar{\delta}}({z}_{j}^{t}) = \emptyset$. 
\end{enumerate}
Here $\mathcal{O}_{\bar{\delta}}({z}_{j}^{t}) = \{{z}: |{z} - {z}_{j}^{t}| < \bar{\delta}\}$ denotes a neighborhood of radius $\bar{\delta}$ around the vertex ${z}_{j}^{t}$. 
\end{lma}

\begin{proof}
This can be established through local analysis near a vertex. The detailed proof can be found in Lemmas 9.1, 9.2, and 9.3 of \cite{JM20}. 
\end{proof}

The proof of \autoref{thm42} consists of the following two lemmas. 

\begin{lma} \label{lma63open}
Let $\bar{t} \in [0, 1]$ be a fixed number such that properties \ref{CGY0601itema}, \ref{CGY0601itemc}, and \ref{CGY0601itemd} in \autoref{thm42} hold for the triangle ${T}^{\bar{t}}$. Then \autoref{thm42} is valid for any triangle ${T}^{t}$ with $|{t} - \bar{t}|$ sufficiently small and ${t} \in (0, 1]$. 
\end{lma}

\begin{proof}
For a positive constant $\varepsilon$, we set
\begin{equation*}
{J}_{\varepsilon} = \{{t} \in (0, 1]: \, |{t} - \bar{t}| < \varepsilon\}. 
\end{equation*}
Note that properties \ref{CGY0601itema} and \ref{CGY0601itemc} of \autoref{thm42} hold for ${u}^{\bar{t}}$ and ${T}^{\bar{t}}$. By continuity, there exists a sufficiently small positive constant $\varepsilon_{1} > 0$ such that for every ${t} \in {J}_{\varepsilon_{1}}$, the second Neumann eigenvalue $\mu_{t}$ is simple, and the following inequalities hold at the vertices: 
\begin{equation*}
{u}^{t}({z}_{1}^{t}) < 0, \quad {u}^{t}({z}_{2}^{t}) < 0, \quad {u}^{t}({z}_{3}^{t}) > 0. 
\end{equation*}
Building upon these properties and the local behavior near the vertices (see \autoref{lma62}), there exist small positive constants $\bar{\delta}$ and $\varepsilon_{2} \in (0, \varepsilon_{1})$ such that for all ${t} \in {J}_{\varepsilon_{2}}$, the following conditions are satisfied: 
\begin{itemize}
\item
$\operatorname{crit}_{\mathrm{nv}}({u}^{t}) \cap \mathcal{O}_{\bar{\delta}}({z}_{j}^{t}) = \emptyset$ for $j = 2, 3$; 
\item
$\operatorname{crit}_{\mathrm{nv}}({u}^{t}) \cap \mathcal{O}_{\bar{\delta}}({z}_{1}^{t}) = \emptyset$ if the interior angle at ${z}_{1}^{\bar{t}}$ is acute (i.e., $\bar{t} > 0$); 
\item
$\operatorname{crit}_{\mathrm{nv}}({u}^{t}) \cap \mathcal{O}_{\bar{\delta}}({z}_{1}^{t})$ is contained only on one side of the triangle ${T}^{t}$. 
\end{itemize}
Here $\mathcal{O}_{\bar{\delta}}({z}_{j}^{t}) = \{{z}: |{z} - {z}_{j}^{t}| < \bar{\delta}\}$. 
Define the region
\begin{equation*}
\mathcal{V}_{\bar{\delta}}^{t} = \{{x}: \, 0 < \operatorname{dist}({x}, \{{z}_{1}^{t}, {z}_{2}^{t}, {z}_{3}^{t}\}) < \bar{\delta}\}. 
\end{equation*}
Within this region, for ${t} \in {J}_{\varepsilon_{2}}$, we have
\begin{equation} \label{CGY0604a}
|\nabla {u}^{t}| > 0 \text{ and } \partial_{{x}_{2}} {u}^{t} \neq 0 \text{ in } \big(\operatorname{Int}(\overline{{z}_{3}^{t}{z}_{1}^{t}}) \cup \operatorname{Int}(\overline{{z}_{3}^{t}{z}_{2}^{t}})\big) \cap \mathcal{V}_{\bar{\delta}}^{t}
\end{equation}
provided that $\bar{t} > 0$. 

Next, we show that \eqref{CGY0604a} is still valid in the case where $\bar{t} = 0$. 
Using the expression \eqref{CGY0602b} for the eigenfunction in the half-equilateral triangle, we have
\begin{equation*} 
\partial_{{x}_{1}} {u}^{0} < 0 \text{ in the compact set } \overline{{z}_{1}^{0}{z}_{2}^{0}} \setminus \mathcal{V}_{\bar{\delta}}^{0}. 
\end{equation*} 
After possibly shrinking $\varepsilon_{2}$, we deduce by continuity that for ${t} \in (0, \varepsilon_{2})$, 
\begin{equation} \label{CGY0604b}
\partial_{{x}_{1}} {u}^{t} < 0 \text{ in } \overline{{z}_{1}^{t}{z}_{2}^{t}} \setminus \mathcal{V}_{\bar{\delta}}^{t}. 
\end{equation}
Furthermore, from \autoref{lma26bJM} and the condition ${u}^{t}({z}_{1}^{t}) < 0$, the point ${z}_{1}^{t}$ is a local minimum of ${u}^{t}$. Together with \eqref{CGY0604b}, this implies that ${u}^{t}$ has at least one critical point in $\operatorname{Int}(\overline{{z}_{1}^{t}{z}_{2}^{t}}) \cap \mathcal{O}_{\bar{\delta}}({z}_{1}^{t})$. 
Combining this with the second statement of \autoref{lma62}, we conclude that ${u}^{t}$ has no critical points in $\operatorname{Int}(\overline{{z}_{3}^{t}{z}_{1}^{t}}) \cap \mathcal{O}_{\bar{\delta}}({z}_{1}^{t})$ for all ${t} \in (0, \varepsilon_{2})$. Consequently, \eqref{CGY0604a} holds for ${t} \in (0, \varepsilon_{2})$ when $\bar{t} = 0$. 

Now, employing property \ref{CGY0601itemd} and \autoref{thm38} for ${u}^{\bar{t}}$, we have
\begin{equation*}
\partial_{{x}_{2}}{u}^{\bar{t}} > 0 \text{ in } \big(\operatorname{Int}(\overline{{z}_{3}^{\bar{t}}{z}_{1}^{\bar{t}}}) \cup \operatorname{Int}(\overline{{z}_{3}^{\bar{t}}{z}_{2}^{\bar{t}}})\big) \setminus \mathcal{V}_{\bar{\delta}}^{\bar{t}}. 
\end{equation*}
By continuity, there exists a constant $\varepsilon_{3} \in (0, \varepsilon_{2})$ such that for every ${t} \in {J}_{\varepsilon_{3}}$, the following holds: 
\begin{equation*}
\partial_{{x}_{2}}{u}^{t} > 0 \text{ in } \big(\operatorname{Int}(\overline{{z}_{3}^{t}{z}_{1}^{t}}) \cup \operatorname{Int}(\overline{{z}_{3}^{t}{z}_{2}^{t}})\big) \setminus \mathcal{V}_{\bar{\delta}}^{t}. 
\end{equation*}
Combining this with \eqref{CGY0604a}, we obtain
\begin{equation} \label{CGY0606}
\partial_{{x}_{2}}{u}^{t} \geq, \not\equiv 0 \text{ in } \partial {T}^{t}
\end{equation}
for ${t} \in {J}_{\varepsilon_{3}}$. From \autoref{lma24Pol}, $\lambda_{1}({T}^{t}) > \mu_{2}({T}^{t}) = \mu_{t}$. 
Combining this with \eqref{CGY0606}, we can apply the maximum principle as described in \autoref{lma21BNV}. This yields that $\partial_{{x}_{2}}{u}^{t} > 0$ in ${T}^{t}$ for all ${t} \in {J}_{\varepsilon_{3}}$. 
Finally, by \autoref{thm38}, property \ref{CGY0601iteme} remains valid. Therefore, properties \ref{CGY0601itema}, \ref{CGY0601itemc}, \ref{CGY0601itemd}, and \ref{CGY0601iteme} hold for ${t} \in {J}_{\varepsilon_{3}}$. 
This completes the proof. 
\end{proof}

\begin{lma} \label{lma64closed}
Let $\bar{t} \in (0, 1]$ be a limit of a sequence $\{t_{n}\} \subset (0, 1]$ such that \autoref{thm42} is valid for every triangle ${T}^{t_{n}}$. 
Then \autoref{thm42} is valid for the triangle ${T}^{\bar{t}}$. 
\end{lma}

\begin{proof}
The proof has two steps. 

\textbf{Step 1. Existence of a monotone eigenfunction corresponding to $\mu_{\bar{t}}$.}
Indeed, it is clear that the eigenvalue ${\mu}_{t_{n}}$ converges to the second Neumann eigenvalue ${\mu}_{\bar{t}}$ of ${T}^{\bar{t}}$; see, e.g., \cite[Theorem 2.3.25]{Hen06}. Clearly, the corresponding normalized eigenfunctions ${u}^{t_{n}}$ are bounded in the ${H}^{1}$-norm, and hence ${u}^{t_{n}} \circ {h}_{t_{n}}$ is bounded in ${H}^{1}({T}^{0})$. For any compact subset ${K}$ of $\overline{{T}^{0}} \setminus \{\text{vertices}\}$, standard regularity theory for elliptic PDEs implies that ${u}^{t_{n}} \circ {h}_{t_{n}}$ is bounded in $C^{k}({K})$ for each $k$. After passing to a subsequence if necessary, we may assume that ${u}^{t_{n}} \circ {h}_{t_{n}}$ converges to $\bar{u} \circ {h}_{\bar{t}}$ strongly in both $L^{2}({T}^{0})$ and $C_{\text{loc}}^{2}(\overline{{T}^{0}} \setminus \{\text{vertices}\})$, and weakly in ${H}^{1}({T}^{0})$ for some function $\bar{u}$. It follows that $\bar{u}$ is a normalized Neumann eigenfunction associated with the eigenvalue ${\mu}_{\bar{t}}$ in ${T}^{\bar{t}}$. 
The positivity of $\partial_{{x}_{2}} {u}^{t_{n}}$ implies the nonnegativity of $\partial_{{x}_{2}} \bar{u}$, and therefore
\begin{equation*}
\partial_{{x}_{2}} \bar{u} > 0 \text{ in } {T}^{\bar{t}}. 
\end{equation*}
Combining this with \eqref{CGY0603b} and \autoref{thm51}, we deduce that $\bar{u}({z}_{1}^{\bar{t}}) \neq 0$ and $\bar{u}({z}_{2}^{\bar{t}}) \neq 0$.
It follows from \autoref{lma26bJM} that $\bar{u}({z}_{1}^{\bar{t}}) < 0$ and $\bar{u}({z}_{2}^{\bar{t}}) < 0$. Now, according to \autoref{thm38}, $\bar{u}$ satisfies properties \ref{CGY0601itemc}, \ref{CGY0601itemd}, and \ref{CGY0601iteme} of \autoref{thm42} for the triangle ${T}^{\bar{t}}$. 

\textbf{Step 2. Uniqueness of eigenfunction corresponding to $\mu_{\bar{t}}$.}
Suppose that $\hat{{u}}$ is an eigenfunction corresponding to $\mu_{\bar{t}}$ such that $\hat{{u}}$ and $\bar{{u}}$ are linearly independent. Then
\begin{equation}
{w}^{s} = \bar{{u}} + {s} \hat{{u}}
\end{equation}
is an eigenfunction corresponding to $\mu_{\bar{t}}$ for any ${s} \in \R$. Let $\mathbf{S}$ be the collection of ${s} \in \mathbb{R}$ such that ${w}^{s}$ satisfies the properties \ref{CGY0601itemc}, \ref{CGY0601itemd} and \ref{CGY0601iteme} on the acute triangle ${T}^{\bar{t}}$. Step 1 implies that $0 \in \mathbf{S}$. By the continuity of the mapping ${s} \mapsto {w}^{s}$ and by the local behavior near the vertices, for ${s}$ close to $0$, ${w}^{s}$ satisfies property \ref{CGY0601iteme} with ${w}^{s}({z}_{1}^{\bar{t}}) < 0$, ${w}^{s}({z}_{2}^{\bar{t}}) < 0$, and ${w}^{s}({z}_{3}^{\bar{t}}) > 0$. It follows that $\partial_{{x}_{2}} {w}^{s} \geq 0$ and $\partial_{{x}_{2}} {w}^{s} \not\equiv 0$ on $\partial {T}^{\bar{t}}$. Combining this with the fact that $\lambda_{1}({T}^{\bar{t}}) > \mu_{\bar{t}}$ (from \autoref{lma24Pol}) and the maximum principle (\autoref{lma21BNV}), we conclude the positivity of $\partial_{{x}_{2}} {w}^{s}$ in ${T}^{\bar{t}}$. Thus, $\mathbf{S}$ contains a neighborhood of ${s} = 0$. This argument indeed implies that $\mathbf{S}$ is open. By a similar argument as in Step 1, $\mathbf{S}$ is closed. Thus, $\mathbf{S}$ is an open, closed, and non-empty set in $\R$. This implies $\mathbf{S} = \R$. Therefore, for any ${s} \in \R$, 
\begin{equation*}
\partial_{{x}_{2}}\bar{{u}} + {s} \partial_{{x}_{2}}\hat{{u}} = \partial_{{x}_{2}}{w}^{s} > 0 \text{ in } {T}^{\bar{t}}, 
\end{equation*}
which implies that $\partial_{{x}_{2}}\hat{{u}} = 0$ in ${T}^{\bar{t}}$. This is a contradiction, completing Step 2. 
\end{proof}

\begin{proof}[Proof of \autoref{thm42}]
Let $\mathcal{S}$ be the collection of those values ${t} \in (0, 1]$ such that \autoref{thm42} is valid for the triangle ${T}^{t}$. Since the eigenfunction ${u}^{0}$ of half-equilateral triangle ${T}^{0}$ is given in \eqref{CGY0602b}, we know that the properties \ref{CGY0601itema}, \ref{CGY0601itemc} and \ref{CGY0601itemd} are satisfied for the triangle ${T}^{0}$. By \autoref{lma63open}, $\mathcal{S}$ is non-empty and relatively open in $(0, 1]$. By \autoref{lma64closed}, $\mathcal{S}$ is closed. 
Thus, $\mathcal{S}$ is a relatively open, closed and non-empty subset of $(0, 1]$. It follows that $\mathcal{S} = (0, 1]$. In particular, \autoref{thm42} is true. 
\end{proof}

\begin{rmk} \label{rmk65}
For an equilateral triangle ${T}$, let ${u}$ be any eigenfunction associated with the second Neumann eigenvalue $\mu$ of ${T}$. Then:
\begin{enumerate}[label = \rm(\roman*), start = 1]
\item \label{CGY0610itf}
The eigenspace corresponding to $\mu$ is two-dimensional.
\item \label{CGY0610ita}
The non-vertex critical point of ${u}$ (if it exists) is unique and is a saddle point; moreover, the non-vertex critical point of ${u}$ exists if and only if ${u}$ does not vanish at all vertices of ${T}$. 
\item \label{CGY0610itc}
${u}$ vanishes at a vertex if and only if ${u}$ is antisymmetric with respect to one of the symmetry axes of of ${T}$.
\item \label{CGY0610itd}
The global extrema are attained only at the vertices; moreover, all three vertices are global extrema if and only if ${u}$ is symmetric with respect to one of the symmetry axes of ${T}$. 
\item \label{CGY0610ite}
${u}$ is monotone in the direction orthogonal to a side ${e}$, provided that ${u}$ does not take opposite signs at the endpoints of ${e}$. 
\end{enumerate}
\end{rmk}

\begin{proof}
Conclusion \ref{CGY0610itf} is standard. 
By applying a rigid transformation to ${u}$ and multiplying ${u}$ by $ - 1$ if necessary, we may assume that ${u}$ satisfies
\begin{equation*}
{u}({z}_{1}) \leq 0, \quad {u}({z}_{2}) < 0, \quad \text{and} \quad {u}({z}_{3}) > 0. 
\end{equation*}
Let $\bar{u}$ be an eigenfunction associated with $\mu$ that is symmetric with respect to the perpendicular bisector of $\overline{{z}_{1}{z}_{2}}$, and satisfies $\partial_{{x}_{2}} \bar{u} > 0$ in ${T}$. Then $\bar{u}({z}_{1}) < 0$, $\bar{u}({z}_{2}) < 0$ and $\bar{u}({z}_{3}) > 0$. Following the same arguments as in the proofs of \autoref{lma63open} and \autoref{lma64closed}, the eigenfunction $\bar{u} + {s}{u}$ is monotonically increasing in the ${x}_{2}$-direction for every ${s} \geq 0$. Hence 
\begin{equation} \label{CGY0612}
\partial_{{x}_{2}}{u} > 0 \text{ in } {T}. 
\end{equation}
Consequently, conclusion~\ref{CGY0610ite} holds both when ${u}({z}_{1}) < 0$ and, by \autoref{lma37}, when ${u}({z}_{1}) = 0$. Combining \eqref{CGY0612} with \autoref{thm38} yields \ref{CGY0610ita}, and combining \eqref{CGY0612} with \autoref{thm51} yields \ref{CGY0610itc}.

To prove conclusion \ref{CGY0610itd}, assume that ${u}({z}_{1}) = {u}({z}_{2}) < 0$. Define ${w} = \bar{u} + {s}_{0}{u}$, where ${s}_{0} = - \bar{u}({z}_{2})/{u}({z}_{2})$. Then we have ${w}({z}_{1}) = {w}({z}_{2}) = 0$. However, the second Neumann eigenfunction cannot vanish at two vertices of the triangle (cf. \autoref{lma26aJM}). This implies that ${w}\equiv0$. Consequently, ${u} = - \bar{u}/s_{0}$, which shows the symmetry of ${u}$. This completes the proof of conclusion \ref{CGY0610itd}. 
\end{proof}

\begin{thm} 
Let ${u}^{1}$ be a second Neumann eigenfunction of ${T}^{1} = \triangle {z}_{1}^{1}{z}_{2}^{1}{z}_{3}^{1}$ satisfying conditions \eqref{CGY0601a} and
\begin{equation} 
|{z}_{2}^{1} - {z}_{3}^{1}| \geq |{z}_{3}^{1} - {z}_{1}^{1}| > |{z}_{1}^{1} - {z}_{2}^{1}| \text{ and } \operatorname{Re}({z}_{3}^{1}) \leq 0. 
\end{equation}
Then all six properties of \autoref{thm12} hold true. 
\end{thm}

\begin{proof}
Since ${T}^{1}$ is a non-acute triangle, all results are covered in \autoref{prop28obt}, except for the third property of \autoref{thm12}. Let the initial triangle ${T}^{0}$ be the half-equilateral triangle with labeled vertices $({z}_{1}^{0}, {z}_{2}^{0}, {z}_{3}^{0})$ as defined in \eqref{CGY0602a}. We define the vertices ${z}_{1}^{t} = 0$, ${z}_{2}^{t} = 1$, and
\begin{align*}
{z}_{3}^{t} & = (1 - 2{t}) \cdot {z}_{3}^{0} + 2{t} \cdot {z}_{3}^{1/2} \text{ for } {t} \in (0, 1/2), \\
{z}_{3}^{t} & = (2 - 2{t}) \cdot {z}_{3}^{1/2} + (2{t} - 1) \cdot {z}_{3}^{1} \text{ for } {t} \in (1/2, 1), 
\end{align*}
where
\begin{equation*}
\operatorname{Re}({z}_{3}^{t}) = 0 \text{ and } \operatorname{Im}({z}_{3}^{t}) = \max\{\sqrt{3}, \tfrac{|{z}_{3}^{1}|^{2}}{\operatorname{Im}({z}_{3}^{1})}\} \text{ at } {t} = 1/2. 
\end{equation*}
Then the triangle ${T}^{t}$ with labeled vertices $({z}_{1}^{t}, {z}_{2}^{t}, {z}_{3}^{t})$ satisfies
\begin{equation*} 
|{z}_{2}^{t} - {z}_{3}^{t}| > |{z}_{3}^{t} - {z}_{1}^{t}| > |{z}_{1}^{t} - {z}_{2}^{t}| \text{ and } \operatorname{Re}({z}_{3}^{t}) \leq 0
\end{equation*}
for every ${t} \in [0, 1]$. From \autoref{prop28obt}, the second Neumann (normalized) eigenfunction ${u}^{t}$ depends continuously on ${t} \in [0, 1]$ and is monotone in the direction parallel to the longest side. Applying \autoref{thm51}, we deduce that ${u}^{t}$ does not vanish at any vertex of ${T}^{t}$, that is, 
\begin{equation} \label{CGY0613a}
{u}^{t}({z}_{3}^{t}) \neq 0, \quad {u}^{t}({z}_{1}^{t}) \neq 0, \quad {u}^{t}({z}_{2}^{t}) \neq 0 \text{ for every } t \in [0, 1]. 
\end{equation}
Since ${T}^{0}$ is a half-equilateral triangle, we have
\begin{equation} \label{CGY0613b}
{u}^{t}({z}_{3}^{t}) \cdot {u}^{t}({z}_{1}^{t}) < 0, \quad
{u}^{t}({z}_{3}^{t}) \cdot {u}^{t}({z}_{2}^{t}) < 0
\end{equation}
for ${t} = 0$. By continuity, it follows that \eqref{CGY0613b} is valid for all ${t} \in [0, 1]$. This completes the proof. 
\end{proof}


\section{The location of global extremum point} \label{Sect7location}

In this section, we establish that the global extrema occur only at the endpoints of the longest side, a phenomenon that was initially discovered numerically by Terence Tao; see the beginning of \href{https://polymathprojects.org/2012/06/12/polymath7-research-thread-1-the-hot-spots-conjecture/}{\textcolor{black}{Polymath7 research thread 1}} in \cite{Pol12}. This section does not depend on \autoref{thm42} as well as \autoref{thm47} and \autoref{thm51}, but relies on the results presented in \autoref{Sect3mon}. 
We provide alternative proofs of several results, including the monotonicity of the eigenfunction, the uniqueness of the non-vertex critical point, and the fact that the second Neumann eigenfunction is unique up to scalar multiplication. The key observation is that the eigenfunction ${u}$ satisfies a partial monotonicity property, namely, it is monotone in the direction parallel to the shortest side in a certain subdomain of the triangle.

We begin by introducing some notations. Let ${T}$ be a triangle with labeled vertices $({z}_{1}, {z}_{2}, {z}_{3})$ satisfying condition \eqref{CGY0601a} and
\begin{equation} \label{CGY0701b}
|{z}_{2} - {z}_{3}| \geq \max\{|{z}_{3} - {z}_{1}|, |{z}_{1} - {z}_{2}|\}. 
\end{equation}
Let $\Upsilon_{S}$ and $\Upsilon_{M}$ denote the perpendicular bisectors of the lower side $\overline{{z}_{1}{z}_{2}}$ and the left side $\overline{{z}_{3}{z}_{1}}$ of ${T}$, respectively. Let $\color{blue}{D}_{S}$ be the right open cap of ${T}$ cut by $\Upsilon_{S}$, and let $\color{blue}{D}_{M}$ be the reflection of the above open cap of ${T}$ cut by $\Upsilon_{M}$; see \autoref{fig7A}. For any point ${x} \in \R^{2}$, let $\Upsilon_{S}{x}$ and $\Upsilon_{M}{x}$ denote the reflections of ${x}$ with respect to the lines $\Upsilon_{S}$ and $\Upsilon_{M}$, respectively. We adopt the definitions of the normal and tangential directions given below \eqref{CGY0106angle} at the end of \autoref{Sect1intro}, that is, 
\begin{gather*}
\boldsymbol{\tau}_{S} = (1, 0), \quad \boldsymbol{n}_{S} = (0, 1), 
\\
\boldsymbol{\tau}_{M} = ( - \cos\alpha_{1}, - \sin\alpha_{1}), \quad \boldsymbol{n}_{M} = (\sin\alpha_{1}, - \cos\alpha_{1}),
\end{gather*}
where $\alpha_{1}$, $\alpha_{2}$ and $\alpha_{3}$ are defined in \eqref{CGY0106angle}. 

\begin{figure}[htp]\centering  \vspace*{-2ex}
\begin{tikzpicture}[scale = 2.6]
\pgfmathsetmacro\AngleL{76} 
\pgfmathsetmacro\AngleM{56} 
\pgfmathsetmacro\AngleS{180-\AngleL-\AngleM} 
\pgfmathsetmacro\LenBase{1.00} 
\pgfmathsetmacro\LenRight{\LenBase*sin(\AngleL)/sin(\AngleS)}
\pgfmathsetmacro\LenLeft{\LenBase*sin(\AngleM)/sin(\AngleS)}
\pgfmathsetmacro\xSS{\LenLeft*cos(\AngleL)} 
\pgfmathsetmacro\ySS{\LenLeft*sin(\AngleL)} 
\pgfmathsetmacro\xLL{0}                     
\pgfmathsetmacro\yLL{0}                     
\pgfmathsetmacro\xMM{\LenBase}              
\pgfmathsetmacro\yMM{0}                     
\pgfmathsetmacro\xH{0.5*\LenBase}
\pgfmathsetmacro\yH{0}
\pgfmathsetmacro\xI{0.5*\LenBase}
\pgfmathsetmacro\yI{(\LenBase-0.5*\LenBase)*tan(\AngleM)}
\fill[fill = green, fill opacity = 0.3, draw = black, very thin] (\xLL, \yLL) -- (\xH, \yH) -- (\xI, \yI) -- cycle;  
\fill[fill = blue,  fill opacity = 0.3, draw = black, very thin] (\xMM, \yMM) -- (\xH, \yH) -- (\xI, \yI) -- cycle;  
\draw[black] 
(\xLL, \yLL) node [left]  {${z}_{1}$} -- (\xMM, \yMM) node [right] {${z}_{2}$} -- (\xSS, \ySS) node [left]  {${z}_{3}$} -- cycle; 
\draw[red, very thick] 
({\xI+1.05*(\xH-\xI)}, {\yI+1.05*(\yH-\yI)}) -- ({\xI-0.3*(\xH-\xI)}, {\yI-0.3*(\yH-\yI)}) node [above] {\footnotesize $\Upsilon_{S}$}; 
\node at ({\xH*4/3}, {\yI*1/3}) {\small ${D}_{S}$}; 
\draw[->] (\LenBase*0.65, \ySS*0.6) -- ++ ({0} : 0.33*\LenBase) node [above] {\tiny $\mathbf{\tau}_{S}$}; 
\draw[->] (\LenBase*0.65, \ySS*0.6) -- ++ ({90} : 0.33*\LenBase) node [right] {\tiny $\mathbf{n}_{S}$}; 
\end{tikzpicture} \hspace{3em} 
\begin{tikzpicture}[scale = 2.6]
\pgfmathsetmacro\AngleL{76}
\pgfmathsetmacro\AngleM{56}
\pgfmathsetmacro\AngleS{180-\AngleL-\AngleM}
\pgfmathsetmacro\LenBase{1.00}
\pgfmathsetmacro\LenRight{\LenBase*sin(\AngleL)/sin(\AngleS)}
\pgfmathsetmacro\LenLeft{\LenBase*sin(\AngleM)/sin(\AngleS)}
\pgfmathsetmacro\xSS{\LenLeft*cos(\AngleL)}
\pgfmathsetmacro\ySS{\LenLeft*sin(\AngleL)}
\pgfmathsetmacro\xLL{0}
\pgfmathsetmacro\yLL{0}
\pgfmathsetmacro\xMM{\LenBase}
\pgfmathsetmacro\yMM{0}
\pgfmathsetmacro\xH{0.5*\LenLeft*cos(\AngleL)}
\pgfmathsetmacro\yH{0.5*\LenLeft*sin(\AngleL)}
\pgfmathsetmacro\xI{\xSS + (0.5*\LenLeft/cos(\AngleS))*cos(\AngleM)}
\pgfmathsetmacro\yI{\ySS - (0.5*\LenLeft/cos(\AngleS))*sin(\AngleM)}
\fill[fill = blue, fill opacity = 0.3, draw = black, very thin] (\xLL, \yLL) -- (\xH, \yH) -- (\xI, \yI) -- cycle; 
\fill[fill = green, fill opacity = 0.3, draw = black, very thin] (\xSS, \ySS) -- (\xH, \yH) -- (\xI, \yI) -- cycle; 
\draw[black] 
(\xLL, \yLL) node [left]  {${z}_{1}$} -- (\xMM, \yMM) node [right] {${z}_{2}$} -- (\xSS, \ySS) node [left]  {${z}_{3}$} -- cycle; 
\draw[red, very thick] 
({\xI+1.4*(\xH-\xI)}, {\yI+1.4*(\yH-\yI)}) -- ({\xI-0.5*(\xH-\xI)}, {\yI-0.5*(\yH-\yI)}) node [below] {\footnotesize $\Upsilon_{M}$}; 
\node at ({(0+\xH+\xI)/3}, {(0+\yH+\yI)/3}) {\small ${D}_{M}$}; 
\draw[->] (-\LenBase*0.22, \LenBase*0.8) -- ++ ({180+\AngleL} : 0.33*\LenBase) node [right] {\tiny $\mathbf{\tau}_{M}$}; 
\draw[->] (-\LenBase*0.22, \LenBase*0.8) -- ++ ({270+\AngleL} : 0.33*\LenBase) node [above] {\tiny $\mathbf{n}_{M}$}; 
\end{tikzpicture}  \vspace*{-2ex}
\caption{The bisectors $\color{red} \Upsilon_{S}$, $\color{red} \Upsilon_{M}$ and subdomains $\color{blue} {D}_{S}$, $\color{blue} {D}_{M}$}
\label{fig7A}
\end{figure}


In this section we will show that the second Neumann eigenfunction ${u}$ (after possibly a sign change) of ${T}$ has the following properties: 
\begin{enumerate}[label = \rm(P\arabic*), start = 1]
\item \label{CGY0702P1}
${u}$ is monotone in a fixed direction throughout ${T}$, specifically, 
\begin{equation} \label{CGY0704}
\text{either } \nabla {u} \cdot \boldsymbol{n}_{S} > 0 \text{ in } {T} \text{ or } \nabla {u} \cdot \boldsymbol{n}_{M} < 0 \text{ in } {T}. 
\end{equation}
\item \label{CGY0702P2}
The following inequalities hold: 
\begin{subequations} \label{CGY0705}
\begin{align}
\label{CGY0705a} 
{u}(\Upsilon_{S}{x}) \geq {u}({x}), & \quad {x} \in {D}_{S}, \\
\label{CGY0705b} 
{u}(\Upsilon_{M}{x}) \geq {u}({x}), & \quad {x} \in {D}_{M}. 
\end{align}
\end{subequations}
\end{enumerate}

The properties \ref{CGY0702P1} and \ref{CGY0702P2} are ``stable" under domain perturbation, i.e., they are preserved for any triangle sufficiently close to ${T}$, as in the proof of \autoref{thm74}. 
As a direct consequence of the key observation \ref{CGY0702P2}, we have
\begin{equation} \label{CGY0705d}
{u}({z}_{3}) \geq {u}({z}_{1}) \geq {u}({z}_{2}). 
\end{equation}
Furthermore, if ${T}$ is neither superequilateral nor equilateral, then \eqref{CGY0705} and \eqref{CGY0705d} hold strictly (see \autoref{lma73}). Consequently, ${z}_{3}$ and ${z}_{2}$ are always the global extrema, while ${z}_{1}$ is a global extremum if and only if ${T}$ is a subequilateral triangle.

\begin{lma} \label{lma71}
Let ${u}$ be an eigenfunction associated with the second Neumann eigenvalue $\mu$ on a non-obtuse triangle ${T}$ with labeled vertices $({z}_{1}, {z}_{2}, {z}_{3})$ satisfying conditions \eqref{CGY0601a} and $|{z}_{3} - {z}_{1}| \leq |{z}_{3} - {z}_{2}|$. Suppose that ${u}$ satisfies
\begin{subequations}
\begin{align}
\label{CGY0707a}
\nabla {u} \cdot \boldsymbol{n}_{S} > 0 & \text{ in } {T}, \\
\label{CGY0707b}
{u}(\Upsilon_{S}{x}) - {u}({x}) \geq 0 & \text{ for } x \in {D}_{S}. 
\end{align}
\end{subequations}
Then ${u}$ satisfies
\begin{equation} \label{CGY0707c}
\nabla {u} \cdot \boldsymbol{\tau}_{S} < 0 \text{ on } ({T} \cup \operatorname{Int}(\overline{{z}_{1}{z}_{2}})) \cap \{{x}_{1} = \lambda\}
\end{equation}
for $\lambda > \operatorname{Re}({z}_{3})$, and \eqref{CGY0707c} holds for $\lambda = \operatorname{Re}({z}_{3})$ if $|{z}_{3} - {z}_{1}| < |{z}_{3} - {z}_{2}|$. 
\end{lma}

\begin{proof}
We employ the moving plane method to establish the desired inequalities. As a direct consequence of \eqref{CGY0707b}, we have $\partial_{{x}_{1}}{u} \leq 0$ on $\Upsilon_{S} \cap {T}$. From the monotonicity property \eqref{CGY0707a} and the Neumann boundary condition of ${u}$ on the right side $\overline{{z}_{3}{z}_{2}}$, it follows that $\partial_{{x}_{1}}{u} \leq 0$ on $\overline{{z}_{3}{z}_{2}}$. As a result, $\partial_{{x}_{1}}{u} \leq 0$ on $\partial {D}_{S} \setminus \overline{{z}_{1}{z}_{2}}$.
Thus, $\partial_{{x}_{1}}{u}$ satisfies 
\begin{equation*}
\begin{cases}
\Delta \partial_{{x}_{1}}{u} + \mu \partial_{{x}_{1}}{u} = 0 & \text{ in } {D}_{S}, 
\\
\partial_{\nu}(\partial_{{x}_{1}}{u}) = 0 & \text{ on } \partial {D}_{S} \cap \overline{{z}_{1}{z}_{2}}, 
\\
\partial_{{x}_{1}}{u} \leq 0 & \text{ on } \partial {D}_{S} \setminus \overline{{z}_{1}{z}_{2}}. 
\end{cases}
\end{equation*}
Using the variational characterization of eigenvalues, we deduce that $\lambda_{1}({D}_{S}, \partial {D}_{S} \setminus \overline{{z}_{1}{z}_{2}}) > \lambda_{1}({T}, \partial {T} \setminus \overline{{z}_{1}{z}_{2}})$. 
From \autoref{lma25LR}, $\lambda_{1}({T}, \partial {T} \setminus \overline{{z}_{1}{z}_{2}}) \geq \mu_{2}({T}) = \mu$. Thus, 
\begin{equation*}
\lambda_{1}({D}_{S}, \partial {D}_{S} \setminus \overline{{z}_{1}{z}_{2}}) > \mu_{2}({T}) = \mu. 
\end{equation*}
It follows from the maximum principle in \autoref{lma22aMP} that $\partial_{{x}_{1}}{u} \leq 0$ in ${D}_{S}$. Moreover, since $\partial_{{x}_{1}}{u}$ does not vanish identically, the strong maximum principle implies the positivity of $\partial_{{x}_{1}}{u}$ in ${D}_{S}$.

\begin{figure}[htp]\centering
\begin{tikzpicture}[scale = 2.6]
\pgfmathsetmacro\AngleL{76}
\pgfmathsetmacro\AngleM{56}
\pgfmathsetmacro\AngleS{180-\AngleL-\AngleM}
\pgfmathsetmacro\LenBase{1.00}
\pgfmathsetmacro\LenRight{\LenBase*sin(\AngleL)/sin(\AngleS)}
\pgfmathsetmacro\LenLeft{\LenBase*sin(\AngleM)/sin(\AngleS)}
\pgfmathsetmacro\xLL{0}      \pgfmathsetmacro\yLL{0}           
\pgfmathsetmacro\xMM{\LenBase} \pgfmathsetmacro\yMM{0}         
\pgfmathsetmacro\xSS{\LenLeft*cos(\AngleL)}
\pgfmathsetmacro\ySS{\LenLeft*sin(\AngleL)}                    
\pgfmathsetmacro\xH{0.5*\LenBase} \pgfmathsetmacro\yH{0}
\pgfmathsetmacro\xI{0.5*\LenBase}
\pgfmathsetmacro\yI{(\LenBase-0.5*\LenBase)*tan(\AngleM)}
\fill[fill = green, fill opacity = 0.3, draw = black, very thin] (\xLL, \yLL) -- (\xH, \yH) -- (\xI, \yI) -- cycle; 
\fill[fill = blue, fill opacity = 0.3, draw = black, very thin] (\xMM, \yMM) -- (\xH, \yH) -- (\xI, \yI) -- cycle; 
\draw[black] (\xLL, \yLL) node [below] {${z}_{1}$} -- (\xMM, \yMM) node [below] {${z}_{2}$} -- (\xSS, \ySS) node [left] {${z}_{3}$} -- cycle; 
\fill[red] (\xH, \yH) circle (0.3pt) node[below left]{\scriptsize ${H}$}; 
\fill[red, very thick] (\xI, \yI) circle (0.3pt) node[right]{\scriptsize ${I}$}; 
\node at ({\xH*4/3}, {\yI*1/3}) {\small ${D}_{S}^{\lambda}$}; 
\draw[red, thick] ({\xH}, {-\LenBase*0.15}) node[above = 2pt, right=-2pt]{\tiny $\Upsilon_{S}^{\lambda}$} -- ({\xI}, {\ySS*1.1}); 
\end{tikzpicture}
\hspace{1em}
\begin{tikzpicture}[scale = 2.6]
\pgfmathsetmacro\AngleL{76}
\pgfmathsetmacro\AngleM{56}
\pgfmathsetmacro\AngleS{180-\AngleL-\AngleM}
\pgfmathsetmacro\LenBase{1.00}
\pgfmathsetmacro\LenRight{\LenBase*sin(\AngleL)/sin(\AngleS)}
\pgfmathsetmacro\LenLeft{\LenBase*sin(\AngleM)/sin(\AngleS)}
\pgfmathsetmacro\xLL{0}      
\pgfmathsetmacro\yLL{0}
\pgfmathsetmacro\xMM{\LenBase} 
\pgfmathsetmacro\yMM{0}
\pgfmathsetmacro\xSS{\LenLeft*cos(\AngleL)}
\pgfmathsetmacro\ySS{\LenLeft*sin(\AngleL)}
\pgfmathsetmacro\LAM{\xSS+(\LenBase/2-\xSS)*0.52}
\pgfmathsetmacro\xH{\LAM} \pgfmathsetmacro\yH{0}
\pgfmathsetmacro\xI{\LAM}
\pgfmathsetmacro\yI{(\LenBase-\LAM)*tan(\AngleM)}
\pgfmathsetmacro\xR{(2*\LAM*tan(\AngleL)-\LenBase*tan(\AngleM))/(tan(\AngleL)-tan(\AngleM))}
\pgfmathsetmacro\yR{-(\xR-\LenBase)*tan(\AngleM)}
\pgfmathsetmacro\xL{2*\LAM-\xR}
\pgfmathsetmacro\yL{\yR}
\fill[fill = blue, fill opacity = 0.3, draw = black, very thin] (\xH, \yH) -- (\xI, \yI) -- (\xR, \yR) -- ({\LAM*2}, 0) -- cycle; 
\fill[fill = green, fill opacity = 0.3, draw = black, very thin] (\xH, \yH) -- (\xI, \yI) -- (\xL, \yL) -- (\xLL, \yLL) -- cycle; 
\draw[black] (\xLL, \yLL) node [below] {${z}_{1}$} -- (\xMM, \yMM) node [below] {${z}_{2}$} -- (\xSS, \ySS) node [left] {${z}_{3}$} -- cycle; 
\fill (\xR, \yR) circle (0.3pt) node[right]{\scriptsize ${J}$}; 
\fill (\xL, \yL) circle (0.3pt) node[left]{\scriptsize ${L}$}; 
\fill ({\LAM*2}, 0) circle (0.3pt) node[below]{\scriptsize ${K}$}; 
\fill[red] (\xH, \yH) circle (0.4pt) node[below left]{\scriptsize ${H}$}; 
\fill[red, very thick] (\xI, \yI) circle (0.4pt) node[right]{\scriptsize ${I}$}; 
\node at ({\xH*4/3}, {\yI*1/3}) {\small ${D}_{S}^{\lambda}$}; 
\draw[red, thick] ({\xH}, {-\LenBase*0.15}) node[above = 2pt, right=-2pt]{\tiny $\Upsilon_{S}^{\lambda}$} -- ({\xI}, {\ySS*1.1}); 
\end{tikzpicture} \hspace{1em}
\begin{tikzpicture}[scale = 2.6]
\pgfmathsetmacro\AngleL{76}
\pgfmathsetmacro\AngleM{56}
\pgfmathsetmacro\AngleS{180-\AngleL-\AngleM}
\pgfmathsetmacro\LenBase{1.00}
\pgfmathsetmacro\LenRight{\LenBase*sin(\AngleL)/sin(\AngleS)}
\pgfmathsetmacro\LenLeft{\LenBase*sin(\AngleM)/sin(\AngleS)}
\pgfmathsetmacro\xLL{0}      
\pgfmathsetmacro\yLL{0}
\pgfmathsetmacro\xMM{\LenBase} 
\pgfmathsetmacro\yMM{0}
\pgfmathsetmacro\xSS{\LenLeft*cos(\AngleL)}
\pgfmathsetmacro\ySS{\LenLeft*sin(\AngleL)}
\pgfmathsetmacro\LAM{\xSS}
\pgfmathsetmacro\xH{\LAM} \pgfmathsetmacro\yH{0}
\pgfmathsetmacro\xI{\LAM} \pgfmathsetmacro\yI{(\LenBase-\LAM)*tan(\AngleM)}
\fill[fill = blue, fill opacity = 0.3, draw = black, very thin] (\xH, \yH) -- (\xI, \yI) -- ({\LAM*2}, 0) -- cycle; 
\fill[fill = green, fill opacity = 0.3, draw = black, very thin] (\xH, \yH) -- (\xI, \yI) -- (\xLL, \yLL) -- cycle; 
\draw[black] (\xLL, \yLL) node [below] {${z}_{1}$} -- (\xMM, \yMM) node [below] {${z}_{2}$} -- (\xSS, \ySS) node [left] {${z}_{3}$} -- cycle; 
\fill ({\LAM*2}, 0) circle (0.3pt) node[below]{\scriptsize ${K}$}; 
\fill[red] (\xH, \yH) circle (0.3pt) node[below left]{\scriptsize ${H}$}; 
\fill[red, very thick] (\xI, \yI) circle (0.3pt) node[right]{\scriptsize ${I}$}; 
\node at ({\xH*4/3}, {\yI*1/3}) {\small ${D}_{S}^{\lambda}$}; 
\draw[red, thick] ({\xH}, {-\LenBase*0.15}) node[above = 2pt, right=-2pt]{\tiny $\Upsilon_{S}^{\lambda}$} -- ({\xI}, {\ySS*1.1}); 
\end{tikzpicture} \vspace*{-2ex}
\caption{The lines $\color{red} \Upsilon_{S}^{\lambda}$ and the moving domains $\color{blue} {D}_{S}^{\lambda}$}
\label{fig7B}
\end{figure}

Now, consider the case where $|{z}_{3} - {z}_{1}| < |{z}_{3} - {z}_{2}|$ and proceed with the standard moving plane method. Define $\Upsilon_{S}^{\lambda} = \{{x} \in \R^{2}: {x}_{1} = \lambda\}$ and set 
\begin{gather*}
{D}_{S}^{\lambda} = \{{x} \in {T}: \, {x}_{1} > \lambda \text{ and } (2\lambda - {x}_{1}, {x}_{2}) \in {T}\}, 
\\
{w}_{\lambda}({x}_{1}, {x}_{2}) = {u}(2\lambda - {x}_{1}, {x}_{2}) - {u}({x}_{1}, {x}_{2})
\end{gather*}
for $\lambda \in [\lambda_{\star}, \lambda^{\star}]$ where
$\lambda^{\star} = ( \operatorname{Re}({z}_{1}) + \operatorname{Re}({z}_{2}))/2$ and $\lambda_{\star} = \operatorname{Re}({z}_{3})$. 
Let ${I}$ denote the intersection point $\overline{{z}_{3}{z}_{1}} \cap \Upsilon_{S}^{\lambda}$ and ${H}$ the intersection \(\overline{{z}_{1}{z}_{2}} \cap \Upsilon_{S}^{\lambda}\). From the positivity of the angular derivative in \eqref{CGY0305a}, we infer that $\nabla {u} \cdot ( - \sin\alpha_{2}, \cos\alpha_{2}) > 0$ on $\operatorname{Int}(\overline{{L}{I}})$, where $\overline{{L}{I}}$ is defined as the reflection of $\partial {D}_{S}^{\lambda} \cap \overline{{z}_{3}{z}_{2}}$ with respect to $\Upsilon_{S}^{\lambda}$. Hence, 
\begin{equation} \label{CGY0709}
\partial_{\nu} {w}_{\lambda} > 0 \text{ on } \partial {D}_{S}^{\lambda} \cap \operatorname{Int}(\overline{{z}_{3}{z}_{2}})
\text{ for } \lambda \in (\lambda_{\star}, \lambda^{\star}], 
\end{equation}
and $\operatorname{Int}(\overline{{L}{I}})$ is empty when $\lambda = \lambda_{\star}$. Thus, the difference function ${w}_{\lambda}$ satisfies
\begin{equation} \label{CGY0710}
\begin{cases}
\Delta {w}_{\lambda} + \mu {w}_{\lambda} = 0 & \text{ in } {D}_{S}^{\lambda}
\\
{w}_{\lambda} = 0 & \text{ on } \partial {D}_{S}^{\lambda} \cap \Upsilon_{S}^{\lambda}, 
\\
\partial_{\nu} {w}_{\lambda} = 0 & \text{ on } \partial {D}_{S}^{\lambda} \cap \overline{{z}_{1}{z}_{2}}, 
\\
\partial_{\nu} {w}_{\lambda} > 0 & \text{ on } \partial {D}_{S}^{\lambda} \setminus (\Upsilon_{S}^{\lambda} \cup \overline{{z}_{1}{z}_{2}}), 
\end{cases}
\end{equation}
for $\lambda = \lambda^{\star}$. 
By the strong maximum principle and the Hopf Lemma, the nonnegative function ${w}_{\lambda^{\star}}$ is positive in ${D}_{S}$ and at smooth boundary points of $\partial {D}_{S} \setminus \Upsilon_{S}$. Additionally, ${w}_{\lambda^{\star}}$ remains positive at the nonsmooth boundary point ${z}_{2} \in \partial {D}_{S}$ (see \cite[Remark 1]{YCG21} or \cite[Lemma 2.4]{CW92}). Thus, for $\lambda = \lambda^{\star}$, 
\begin{equation} \label{CGY0711}
{w}_{\lambda} > 0 \text{ in } \overline{{D}_{S}^{\lambda}} \setminus \Upsilon_{S}^{\lambda} \text{ and } \partial_{{x}_{1}}{u} < 0 \text{ on } \Upsilon_{S}^{\lambda} \cap {T}. 
\end{equation}
Now we let $\bar{\lambda}$ be the supremum of those $\lambda \in [\lambda_{\star}, \lambda^{\star}]$ such that \eqref{CGY0711} fails to hold, i.e., 
\begin{equation*}
\bar{\lambda} = \inf\{\lambda' \in [\lambda_{\star}, \lambda^{\star}]: \, \eqref{CGY0711}\mbox{ holds for every }\lambda \in [\lambda', \lambda^{\star}]\}. 
\end{equation*}
To show that $\bar{\lambda} = \lambda_{\star}$, assume by contradiction that $\bar{\lambda} > \lambda_{\star}$. By continuity, ${w}_{\bar{\lambda}}$ is nonnegative in ${D}_{S}^{\bar{\lambda}}$, and $\partial_{{x}_{1}}{u} \leq 0$ on $\Upsilon_{S}^{\bar{\lambda}} \cap {T}$. Now for any point ${p} = ({p}_{1}, {p}_{2})$ lying on the line segment $\overline{{z}_{3}{z}_{2}}$ with ${p}_{1} > \bar{\lambda}$, we have
\begin{equation*}
{R}_{{p}}{u} \leq, \not\equiv0 \text{ in } \partial({T} \cap \{{x}_{1} > {p}_{1}\}). 
\end{equation*}
By \autoref{lma24Pol}, $\lambda_{1}({T} \cap \{{x}_{1} > {p}_{1}\}) > \lambda_{1}({T}) > \mu_{2}(T) = \mu$. Applying the maximum principle as in \autoref{lma21BNV} to ${R}_{{p}}{u}$, we obtain that
\begin{equation*} 
{R}_{{p}}{u} < 0 \text{ in } {T} \cap \{{x}_{1} > {p}_{1}\} \text{ whenever } {p} \in \overline{{z}_{3}{z}_{2}} \text{ satisfying }{p}_{1} \geq \lambda
\end{equation*}
for $\lambda > \bar{\lambda}$. Combining this with \eqref{CGY0709}, we deduce that \eqref{CGY0710} holds for both $\lambda = \bar{\lambda}$ and $\lambda \in [\bar{\lambda} - \varepsilon_{1}, \bar{\lambda}]$ for some $\varepsilon_{1} > 0$. Hence, the nonnegative function ${w}_{\bar{\lambda}}$ satisfies \eqref{CGY0711} for $\lambda = \bar{\lambda}$. Combining this with \autoref{lma22cMP}, we know that the first mixed eigenvalue $\lambda_{1}({D}_{S}^{\bar{\lambda}}, \partial {D}_{S}^{\bar{\lambda}} \cap \Upsilon_{S}^{\bar{\lambda}})$ is greater than $\mu$. By continuity, 
\begin{equation*}
\lambda_{1}({D}_{S}^{\lambda}, \partial {D}_{S}^{\lambda} \cap \Upsilon_{S}^{\lambda}) > \mu \text{ for } \lambda \in [\bar{\lambda} - \varepsilon_{2}, \bar{\lambda}]
\end{equation*}
for some $\varepsilon_{2} \in (0, \varepsilon_{1})$. Now applying the maximum principle \autoref{lma22aMP} to ${w}_{\lambda}$-equation, we get that \eqref{CGY0711} is valid for every $\lambda \in [\bar{\lambda} - \varepsilon_{2}, \bar{\lambda}]$. This contradicts the definition of $\bar{\lambda}$. Therefore, $\bar{\lambda} = \lambda_{\star}$, and consequently, \eqref{CGY0711} holds for all $\lambda \in [\lambda_{\star}, \lambda^{\star}]$. 

Finally, we demonstrate the strict monotonicity of ${u}$ along a portion of $\overline{{z}_{1}{z}_{2}}$, i.e., 
\begin{equation} \label{CGY0707d}
\nabla {u} \cdot \boldsymbol{\tau}_{S} < 0 \text{ on } \operatorname{Int}(\overline{{z}_{1}{z}_{2}}) \cap \{{x}_{1} = \lambda\}
\end{equation}
holds for $\lambda > \lambda_{\star}$, and \eqref{CGY0707d} holds for $\lambda = \lambda_{\star}$ if $|{z}_{3} - {z}_{1}| < |{z}_{3} - {z}_{2}|$. 
In fact, since $\partial_{{x}_{1}} {u} < 0$ in ${D}_{S}^{\lambda_{\star}}$, it follows from the proof of the third part of \autoref{lma31} that \eqref{CGY0707d} holds for $\lambda > \lambda_{\star}$. In the case when $|{z}_{3} - {z}_{1}| < |{z}_{3} - {z}_{2}|$, recalling that ${w}_{\lambda_{\star}} < 0$ in ${D}_{S}^{\lambda_{\star}}$, one may apply Serrin's corner lemma to ${w}_{\lambda_{\star}}$ to deduce that $\nabla {u} \cdot \boldsymbol{\tau}_{S} < 0$ at $\overline{{z}_{3}{z}_{1}} \cap \{{x}_{1} = \lambda_{\star}\}$ (see details in \cite[Theorem 2.4]{BP89}). Alternatively, this result can be obtained via the Hopf lemma applied to the linear function of the extended function of ${w}_{\lambda_{\star}}$ with reflection along the side $\overline{{z}_{1}{z}_{2}}$. This completes the proof. 
\end{proof}

\begin{lma} \label{lma72}
Let ${u}$ be an eigenfunction associated with the second Neumann eigenvalue $\mu$ on an acute triangle ${T}$ with labeled vertices $({z}_{1}, {z}_{2}, {z}_{3})$ satisfying conditions \eqref{CGY0601a} and \eqref{CGY0701b}.  Suppose that ${u}$ satisfies properties \ref{CGY0702P1} and \ref{CGY0702P2}. Then ${u}$ has at most one non-vertex critical point, which is a saddle point. 
\end{lma}

\begin{proof}
\textbf{Step 1}. 
We show that ${u}$ has no critical point on $\operatorname{Int}(\overline{{z}_{2}{z}_{3}}) \cup \operatorname{Int}(\overline{{z}_{1}{z}_{2}})$ if 
\begin{equation} \label{CGY0715}
\nabla {u} \cdot \boldsymbol{n}_{M} < 0 \text{ in } {T} \text{ and } {u}(\Upsilon_{M} {x}) - {u}({x}) \geq 0 \text{ for } {x} \in {D}_{M}.
\end{equation}
In fact, combining the assumption $|\overline{{z}_{2}{z}_{3}}| \geq |\overline{{z}_{1}{z}_{2}}|$ with the proof of \autoref{thm38}, one deduces that $\operatorname{crit}_{\mathrm{nv}}({u}) \cap \operatorname{Int}(\overline{{z}_{1}{z}_{2}}) = \emptyset$. Under assumption \eqref{CGY0715}, the same argument as in the proof of \autoref{lma71} implies that
\begin{equation} \label{CGY0716}
\nabla {u} \cdot \boldsymbol{\tau}_{M} < 0 \text{ in } \mathcal{D} \text{ and } \nabla {u} \cdot \boldsymbol{n}_{M} < 0 \text{ in } {T},
\end{equation}
where $\mathcal{D} = {T} \cap \{{x}: ({x} - {z}_{2}) \cdot \boldsymbol{\tau}_{M} < 0\}$. In the remainder of the proof, we present two methods to show that $\operatorname{crit}_{\mathrm{nv}}({u}) \cap \overline{{z}_{2}{z}_{3}} = \emptyset$. 

\textbf{Method 1}. 
Suppose that $|\nabla {u}({p})| = 0$ for some ${p} \in \operatorname{Int}(\overline{{z}_{2}{z}_{3}})$. Note that both $\nabla {u} \cdot \boldsymbol{\tau}_{M}$ and $\nabla {u} \cdot \boldsymbol{n}_{M}$ attain their local maxima at ${p}$. Applying the Hopf lemma to the linear equations satisfied by $\nabla {u} \cdot \boldsymbol{\tau}_{M}$ and $\nabla {u} \cdot \boldsymbol{n}_{M}$, we conclude that
\begin{equation*}
\nabla \left( \nabla {u} \cdot \boldsymbol{\tau}_{M} \right) \cdot \boldsymbol{n}_{M} > 0 \text{ and } \nabla \left( \nabla {u} \cdot \boldsymbol{n}_{M} \right) \cdot \boldsymbol{\tau}_{M} < 0 \text{ at } {p},
\end{equation*}
which yields a contradiction. Therefore, $|\nabla {u}| > 0$ on $\operatorname{Int}(\overline{{z}_{2}{z}_{3}})$. 

\textbf{Method 2}. 
Denote by $\boldsymbol{\tau}_{L} = ( - \cos{\alpha}_{2}, \sin{\alpha}_{2})$ the tangential direction along the longest side $\overline{{z}_{2}{z}_{3}}$. From \eqref{CGY0716}, it follows that $\nabla {u} \cdot \boldsymbol{\tau}_{L} > 0$ in $\mathcal{D}$. Let $\tilde{u}$ and $\tilde{\mathcal{D}}$ be the extensions of ${u}$ and $\mathcal{D}$ along the side $\overline{{z}_{3}{z}_{1}}$ by reflection. Then $\operatorname{Int}(\overline{{z}_{2}{z}_{3}})$ lies in the interior of the kite $\tilde{\mathcal{D}}$, and $\nabla \tilde{u} \cdot \boldsymbol{\tau}_{L} \geq, \not\equiv 0$ in $\tilde{\mathcal{D}}$. Hence, the strong maximum principle implies that $\nabla \tilde{u} \cdot \boldsymbol{\tau}_{L} > 0$ in $\tilde{\mathcal{D}}$. Consequently, $\nabla {u} \cdot \boldsymbol{\tau}_{L} > 0$ on $\operatorname{Int}(\overline{{z}_{2}{z}_{3}})$. 

\textbf{Step 2}. 
One can prove that ${u}$ has no critical point on $\operatorname{Int}(\overline{{z}_{2}{z}_{3}}) \cup \operatorname{Int}(\overline{{z}_{3}{z}_{1}})$ if 
\begin{equation*}
\nabla {u} \cdot \boldsymbol{n}_{S} > 0 \text{ in } {T} \text{ and } {u}(\Upsilon_{S} {x}) - {u}({x}) \geq 0 \text{ for } {x} \in {D}_{S}.
\end{equation*}
The proof is similar to that of Step 1, and we omit the details. 

Combining these two steps with \autoref{lma35}, it follows that the non-vertex critical point (if it exists) is unique, is a saddle point, and does not lie on the longest side of ${T}$. 
\end{proof}

\begin{lma} \label{lma73}
Let ${u}$ be an eigenfunction associated with the second Neumann eigenvalue $\mu$ on an acute triangle ${T}$ with labeled vertices $({z}_{1}, {z}_{2}, {z}_{3})$ satisfying conditions \eqref{CGY0601a} and \eqref{CGY0701b}. Suppose that ${u}$ satisfies properties \ref{CGY0702P1} and \ref{CGY0702P2}. 
Then we have the following strict inequalities: 
\begin{equation} \label{CGY0718a}
\nabla {u} \cdot \boldsymbol{\tau}_{S} < 0 \text{ on } \Upsilon_{S} \cap \overline{{T}}, \thinspace\thinspace {u}(\Upsilon_{S}{x}) > {u}({x}) \text{ for } {x} \in \overline{{D}_{S}} \setminus \Upsilon_{S}, \thinspace\thinspace {u}({z}_{1}) > {u}({z}_{2})
\end{equation}
provided $|{z}_{3} - {z}_{2}| > |{z}_{3} - {z}_{1}|$, and
\begin{equation} \label{CGY0718b}
\nabla {u} \cdot \boldsymbol{\tau}_{M} < 0 \text{ on } \Upsilon_{M} \cap \overline{{T}}, \thinspace\thinspace {u}(\Upsilon_{M}{x}) > {u}({x}) \text{ for } {x} \in \overline{{D}_{M}} \setminus \Upsilon_{M}, \thinspace\thinspace {u}({z}_{3}) > {u}({z}_{1})
\end{equation}
provided $|{z}_{2} - {z}_{3}| > |{z}_{2} - {z}_{1}|$.
\end{lma}

\begin{proof}
From \autoref{lma72}, we observe that
\begin{equation*}
|\nabla {u}| > 0, \quad \nabla {u} = |\nabla {u}|\boldsymbol{\tau}_{L} \text{ on } \operatorname{Int}(\overline{{z}_{3}{z}_{2}}), 
\end{equation*}
where $\boldsymbol{\tau}_{L} = ( - \cos\alpha_{2}, \sin\alpha_{2})$ is the tangential direction along the longest side $\overline{{z}_{2}{z}_{3}}$. 

Now we assume that $|{z}_{2} - {z}_{3}| > |{z}_{2} - {z}_{1}|$. In this case, $\operatorname{Int}(\overline{{I}{z}_{1}})$ is contained in ${T}$, where ${I}$ is the intersection point of the two lines $\Upsilon_{M}$ and $\overline{{z}_{2}{z}_{3}}$; see the right picture in \autoref{fig7A}. 
Under each monotonicity assumption in \eqref{CGY0704}, following the same argument as in the proof of \autoref{lma31}, one can prove that the angular derivative ${R}_{{z}_{1}}{u}$ satisfies 
\begin{equation*}%
{R}_{{z}_{1}}{u} > 0 \text{ in } {T}.
\end{equation*}
It follows that the nonnegative function ${w}({x}) = {u}(\Upsilon_{M}{x}) - {u}({x})$ satisfies
\begin{equation*}\begin{cases}
\Delta {w} + \mu {w} = 0 & \text{ in } {D}_{M}
\\
{w} = 0 & \text{ on } \partial {D}_{M} \cap \Upsilon_{M}, 
\\
\partial_{\nu} {w} = 0 & \text{ on } \partial {D}_{M} \cap \overline{{z}_{3}{z}_{1}}, 
\\
\partial_{\nu} {w} > 0 & \text{ on } \operatorname{Int}(\overline{{I}{z}_{1}}) \subset \partial {D}_{M}. 
\end{cases}\end{equation*}
The strong maximum principle and the Hopf lemma imply that
\begin{gather*}
{w}({x}) = {u}(\Upsilon_{M}{x}) - {u}({x}) > 0 \text{ for } {x} \in \overline{{D}_{M}} \setminus \Upsilon_{M}
\text{ and } 
\nabla {u} \cdot \boldsymbol{\tau}_{M} < 0 \text{ on } \Upsilon_{M} \cap {T}. 
\end{gather*}
The negativity of $\nabla {u} \cdot \boldsymbol{\tau}_{M}$ follows by the same proof as for \eqref{CGY0707d}. Therefore, inequality \eqref{CGY0718b} is established. Similarly, one can deduce inequality \eqref{CGY0718a}. 
\end{proof}

\begin{thm} \label{thm74}
Let ${T} = \triangle {z}_{1}{z}_{2}{z}_{3}$ be an acute triangle satisfying conditions \eqref{CGY0601a} and \eqref{CGY0701b}, and assume ${T}$ is non-equilateral. Let ${u}$ be an eigenfunction associated with the second Neumann eigenvalue $\mu$. Then the following conclusions hold. 
\begin{enumerate}[label = \rm(\roman*), start = 1]
\item 
The second Neumann eigenfunction is unique up to a multiplication by a constant. 
\item 
The global extrema of ${u}$ occur only at the endpoints of the longest side, and the nodal line $\mathcal{Z}({u})$ intersects the interior of the longest side. 
\item 
The non-vertex critical point, if it exists, is unique. It is a saddle point and lies in the interior of either the shortest or intermediate side of the triangle. 
\item 
The non-vertex critical point exists if and only if ${u}$ does not vanish at the vertex with the largest interior angle. 
\item 
After a possible sign change, ${u}$ has the properties \ref{CGY0702P1} and \ref{CGY0702P2}. 
\end{enumerate}
\end{thm}

\begin{proof} 
We begin by selecting the initial triangle ${T}^{0}$ as an isosceles right triangle with labeled vertices $({z}_{1}^{0}, {z}_{2}^{0}, {z}_{3}^{0})$ defined as
\begin{equation}
{z}_{1}^{0} = 0 \in \mathbb{C}, \quad {z}_{2}^{0} = 1 \in \mathbb{C}, \text{ and } {z}_{3}^{0} = \sqrt{ - 1} \in \mathbb{C}. 
\end{equation}
It is well known that the eigenspace associated with the second Neumann eigenvalue $\mu_{0}$ of ${T}^{0}$ is one-dimensional and spanned by the normalized eigenfunction ${u}^{0}$: 
\begin{equation}
{u}^{0} = \sqrt{2}\big(\cos (\pi {x}_{1}) - \cos (\pi {x}_{2})\big), \quad \mu_{0} = \pi^{2}. 
\end{equation}
From the explicit expression of ${u}^{0}$, the theorem holds trivially for the triangle ${T}^{0}$. 

Let ${T}^{1} = {T}$ be the target triangle with labeled vertices $({z}_{1}^{1}, {z}_{2}^{1}, {z}_{3}^{1})$, and let ${T}^{t}$ be a continuous family of triangles connecting ${T}^{1}$ and ${T}^{0}$, as defined in \eqref{CGY0603a}. One can verify that, for every ${t} \in [0, 1)$, the following inequalities hold: 
\begin{equation} \label{CGY0720}
|{z}_{2}^{t} - {z}_{3}^{t}| > \max\{|{z}_{3}^{t} - {z}_{1}^{t}|, |{z}_{1}^{t} - {z}_{2}^{t}|\}. 
\end{equation}
Similar to the notation for the triangle ${T}^{1} = {T}$, let $\Upsilon_{S}^{t}$, $\Upsilon_{M}^{t}$, ${D}_{S}^{t}$, and ${D}_{M}^{t}$ denote the perpendicular bisectors and corresponding subdomains of the triangle ${T}^{t}$. The second Neumann eigenvalue $\mu_{t}$ of ${T}^{t}$ depends continuously on ${t}$. Furthermore, the corresponding normalized eigenfunction ${u}^{t}$ depends continuously on ${t}$ when the eigenvalue is simple. Define
\begin{equation*}
\mathcal{S} = \{{t} \in [0, 1]: \, \text{ \autoref{thm74} holds for the triangle ${T}^{t}$}\}. 
\end{equation*}

\textbf{Step 1}. 
The set $\mathcal{S} \cap [0, 1)$ is a relatively open subset of $[0, 1)$. 
In fact, suppose that \autoref{thm74} holds for the triangle ${T}^{\bar{t}}$ for some $\bar{t} \in [0, 1)$. By continuity, the second Neumann eigenvalue $\mu_{t}$ is simple for $|{t} - \bar{t}| \ll 1$. From \eqref{CGY0720}, both $\Upsilon_{S}^{\bar{t}}$ and $\Upsilon_{M}^{\bar{t}}$ do not pass through any vertex of ${T}^{\bar{t}}$. By \autoref{lma73}, the boundary tangential derivatives of ${u}^{\bar{t}}$ at both $\Upsilon_{S}^{\bar{t}} \cap \partial {T}^{\bar{t}}$ and $\Upsilon_{M}^{\bar{t}} \cap \partial {T}^{\bar{t}}$ are non-zero. Moreover, 
\begin{equation} \label{CGY0721b}
\begin{aligned}
{u}^{t}(\Upsilon_{S}^{t}{x}) > {u}^{t}({x}) \text{ in } \overline{{D}_{S}^{t}} \setminus \Upsilon_{S}^{t}, & \quad 
\nabla {u}^{t} \cdot \boldsymbol{\tau}_{S}^{t} < 0 \text{ on } \Upsilon_{S}^{t} \cap \overline{{T}^{t}}, 
\\
{u}^{t}(\Upsilon_{M}^{t}{x}) > {u}^{t}({x}) \text{ in } \overline{{D}_{M}^{t} } \setminus \Upsilon_{M}^{t}, & \quad 
\nabla {u}^{t} \cdot \boldsymbol{\tau}_{M}^{t} < 0 \text{ on } \Upsilon_{M}^{t} \cap \overline{{T}^{t}}
\end{aligned}
\end{equation}
is valid for ${t} = \bar{t}$. Consequently, \eqref{CGY0721b} holds for all $t$ near $\bar{t}$. It remains to verify the monotonicity property \eqref{CGY0704}. By \eqref{CGY0721b} and \autoref{lma72}, one deduces that the unique non-vertex critical point is a saddle point, and that ${u}^{\bar{t}}$ attains its global extrema only at ${z}_{3}^{\bar{t}}$ and ${z}_{2}^{\bar{t}}$. Hence, 
\begin{equation*}
{u}^{{t}}({z}_{3}^{{t}}) > 0 \text{ and } {u}^{{t}}({z}_{2}^{{t}}) < 0
\end{equation*}
for ${t} = \bar{t}$ and consequently for ${t}$ close to $\bar{t}$. 

\textbf{Case 1.1}. 
${u}^{\bar{t}}({z}_{1}^{\bar{t}}) < 0$ and ${T}^{\bar{t}}$ is acute. In this scenario, by \autoref{lma26bJM}, we have $\nabla {u}^{\bar{t}} \cdot \boldsymbol{n}_{M}^{\bar{t}} > 0$ in a deleted neighborhood of ${z}_{1}^{\bar{t}}$. It then follows from \eqref{CGY0704} that $\nabla {u}^{\bar{t}} \cdot \boldsymbol{n}_{S}^{\bar{t}} > 0$ in ${T}^{\bar{t}}$. Combining this with \eqref{CGY0705a} and \autoref{lma72}, we derive that
\begin{equation} \label{CGY0722a}
\nabla {u}^{{t}} \cdot \boldsymbol{n}_{S}^{{t}} > 0 \text{ in } {T}^{{t}} \cup \operatorname{Int}(\overline{{z}_{3}^{{t}}{z}_{1}^{{t}}}) \cup \operatorname{Int}(\overline{{z}_{3}^{{t}}{z}_{2}^{{t}}})
\end{equation}
for ${t} = \bar{t}$. Combining this with the third conclusion of \autoref{lma62}, we deduce that
\begin{equation*}
\nabla {u}^{{t}} \cdot \boldsymbol{n}_{S}^{{t}} > 0 \text{ on } \operatorname{Int}(\overline{{z}_{3}^{{t}}{z}_{1}^{{t}}}) \cup \operatorname{Int}(\overline{{z}_{3}^{{t}}{z}_{2}^{{t}}}) \text{ and } \nabla {u}^{{t}} \cdot \boldsymbol{n}_{S}^{{t}} = 0 \text{ on } \operatorname{Int}(\overline{{z}_{1}^{{t}}{z}_{2}^{{t}}}). 
\end{equation*}
Thanks to $\lambda_{1}({T}^{t}) > \mu_{2}({T}^{t}) = \mu_{t}$ (cf. \autoref{lma24Pol}) and the maximum principle as described in \autoref{lma21BNV}, we deduce that \eqref{CGY0722a} remains valid for $|{t} - \bar{t}| \ll 1$. 

\textbf{Case 1.2}. 
${u}^{\bar{t}}({z}_{1}^{\bar{t}}) > 0$ and ${T}^{\bar{t}}$ is acute. In this case, by a similar argument as in Case 1.1, one can show that for $|{t} - \bar{t}| \ll 1$, there holds
\begin{equation} \label{CGY0722b}
\nabla {u}^{{t}} \cdot \boldsymbol{n}_{M}^{{t}} < 0 \text{ in } {T}^{{t}} \cup \operatorname{Int}(\overline{{z}_{2}^{{t}}{z}_{1}^{{t}}}) \cup \operatorname{Int}(\overline{{z}_{2}^{{t}}{z}_{3}^{t}}). 
\end{equation}

\textbf{Case 1.3}. 
The case where either ${u}^{\bar{t}}({z}_{1}^{\bar{t}}) = 0$ or ${T}^{\bar{t}}$ is non-acute (i.e., $\bar{t} = 0$). 
In this situation, from \eqref{CGY0704} and \autoref{lma35}, we deduce that ${u}^{\bar{t}}$ does not admit any non-vertex critical points. Consequently, both \eqref{CGY0722a} and \eqref{CGY0722b} are valid for ${t} = \bar{t}$ by \autoref{lma37}. Combining this with \autoref{lma62}, we have for $|{t} - \bar{t}| \ll 1$, 
\begin{equation} \label{CGY0723}
\text{Either } \operatorname{crit}_{\mathrm{nv}}({u}^{t}) \subset \overline{{z}_{1}^{t}{z}_{2}^{t}} \text{ or } \operatorname{crit}_{\mathrm{nv}}({u}^{t}) \subset \overline{{z}_{1}^{t}{z}_{3}^{t}}. 
\end{equation}
Set ${v}^{t} = \nabla {u}^{t} \cdot \boldsymbol{n}_{S}^{t}$ when $\operatorname{crit}_{\mathrm{nv}}({u}^{t}) \subset \overline{{z}_{1}^{t}{z}_{2}^{t}}$, and ${v}^{t} = - \nabla {u}^{t} \cdot \boldsymbol{n}_{M}^{t}$ when $\operatorname{crit}_{\mathrm{nv}}({u}^{t}) \subset \overline{{z}_{1}^{t}{z}_{3}^{t}}$. From \eqref{CGY0723}, we have
\begin{equation*}
\Delta {v}^{t} + \mu_{t}{v}^{t} = 0 \text{ in } {T}^{t} \text{ and } {v}^{t} \geq, \not\equiv 0 \text{ on } \partial {T}^{t}. 
\end{equation*}
From \autoref{lma24Pol}, $\lambda_{1}({T}^{t}) > \mu_{t}$. Applying the maximum principle from \autoref{lma21BNV}, we conclude that ${v}^{t} > 0$ in ${T}^{t}$. Combining this with \autoref{lma72}, either \eqref{CGY0722a} or \eqref{CGY0722b} holds for ${T}^{t}$ when ${t}$ is close to $\bar{t}$. 

By \autoref{lma35} and \eqref{CGY0721b}, all the properties hold for the triangle ${T}^{t}$ when $|{t} - \bar{t}| \ll 1$. 

\textbf{Step 2}. 
The set $\mathcal{S}$ is a closed subset of $[0, 1]$. 
Indeed, suppose that $\lim_{n\to\infty}{t}_{n} = \bar{t}$ with ${t}_{n} \in \mathcal{S}$. That is, \autoref{thm74} holds for the triangle ${T}^{{t}_{n}}$. We proceed with the proof in three substeps. 

\textit{Step 2.1}. 
Existence of an eigenfunction associated with ${\mu}_{\bar{t}}$ in ${T}^{\bar{t}}$ satisfying \ref{CGY0702P1} and \ref{CGY0702P2}. To establish this, following the convergence argument of eigenfunctions in Step 1 of \autoref{lma64closed}, one can find a normalized Neumann eigenfunction $\bar{u}$ associated with the eigenvalue ${\mu}_{\bar{t}}$ in ${T}^{\bar{t}}$. Due to the fact that ${u}^{{t}_{n}}$ satisfies \eqref{CGY0704} and \eqref{CGY0705} in ${T}^{{t}_{n}}$, it follows that the limit function $\bar{u}$ satisfies \eqref{CGY0705} and
\begin{equation*} 
\text{either } \nabla \bar{u} \cdot \boldsymbol{n}_{S}^{\bar{t}} \geq 0 \text{ in } {T}^{\bar{t}}
\text{ or } \nabla \bar{u} \cdot \boldsymbol{n}_{M}^{\bar{t}} \leq 0 \text{ in } {T}^{\bar{t}}. 
\end{equation*}
It follows from \autoref{lma31} that $\bar{u}$ satisfies the conditions \ref{CGY0702P1} and \ref{CGY0702P2} in ${T}^{\bar{t}}$, that is, 
\begin{equation} \label{CGY0726}
\begin{gathered}
\text{either } \nabla \bar{u} \cdot \boldsymbol{n}_{S}^{\bar{t}} > 0 \text{ in } {T}^{\bar{t}}
\text{ or } \nabla \bar{u} \cdot \boldsymbol{n}_{M}^{\bar{t}} < 0 \text{ in } {T}^{\bar{t}}, 
\\
\bar{u}(\Upsilon_{S}^{\bar{t}}{x}) \geq \bar{u}({x}) \text{ for } {x} \in {D}_{S}^{\bar{t}}, \quad
\bar{u}(\Upsilon_{M}^{\bar{t}}{x}) \geq \bar{u}({x}) \text{ for } {x} \in {D}_{M}^{\bar{t}}. 
\end{gathered}
\end{equation}
By \autoref{lma72} and \autoref{lma35}, $\bar{u}$ satisfies all the properties stated in the theorem, except for the simplicity of the eigenvalue and the exact location of the global extrema. 

\textit{Step 2.2}. 
The eigenvalue $\mu_{\bar{t}}$ is simple. 
In fact, if it is false, there exists another eigenfunction ${v}$ corresponding to $\mu_{\bar{t}}$ that is linearly independent of $\bar{u}$. Define
\begin{equation*}
{w}^{s} = \bar{u} + {s}{v}. 
\end{equation*}

Case 1: ${T}^{\bar{t}}$ is not subequilateral. 
By the same argument as in Step 1 and Step 2.1, ${w}^{s}$ satisfies \eqref{CGY0705} and \eqref{CGY0704} for every ${s} \in \R$. In particular, for all ${s} \in \R$, 
\begin{equation*}
\bar{u}({z}_{3}^{\bar{t}}) + {s}{v}({z}_{3}^{\bar{t}}) = {w}^{s}({z}_{3}^{\bar{t}}) = \max_{ \overline{{T}^{\bar{t}}}} {w}^{s} > 0 \text{ and } 
\bar{u}({z}_{2}^{\bar{t}}) + {s}{v}({z}_{2}^{\bar{t}}) = {w}^{s}({z}_{2}^{\bar{t}}) = \min_{ \overline{{T}^{\bar{t}}}} {w}^{s} < 0. 
\end{equation*}
This implies that ${v}({z}_{3}^{\bar{t}}) = {v}({z}_{2}^{\bar{t}}) = 0$, leading to a contradiction (see \autoref{lma26aJM}). 

Case 2: ${T}^{\bar{t}}$ is subequilateral and $\bar{u}$ is antisymmetric. In this case, we may assume that $|{z}_{3}^{\bar{t}} - {z}_{1}^{\bar{t}}| > |{z}_{1}^{\bar{t}} - {z}_{2}^{\bar{t}}|$. After choosing ${v}(\Upsilon_{S}^{\bar{t}}{x}) + {v}({x})$ if necessary, we can assume that ${v}$ is symmetric. It is evident that ${w}^{s}$ satisfies \eqref{CGY0705a} for ${s} \in \R$. By the same argument as in Step 1 and Step 2.1, the eigenfunction ${w}^{s}$ satisfies \eqref{CGY0705b} and \eqref{CGY0704} for ${s} \in \R$. Therefore, ${v}({z}_{3}^{\bar{t}}) = {v}({z}_{2}^{\bar{t}}) = 0$, which is impossible. 

Case 3: ${T}^{\bar{t}}$ is subequilateral and $\bar{u}$ is not antisymmetric. In this case, we may assume that $|{z}_{3}^{\bar{t}} - {z}_{1}^{\bar{t}}| > |{z}_{1}^{\bar{t}} - {z}_{2}^{\bar{t}}|$ and we may assume that ${v}$ is antisymmetric and ${v} < 0$ in ${D}_{S}^{\bar{t}}$. It is clear that ${w}^{s}$ satisfies \eqref{CGY0705a} for ${s} \in \R^{ + }$. Applying the same argument as in Step 1 and Step 2.1, ${w}^{s}$ satisfies \eqref{CGY0705b} and \eqref{CGY0704} for ${s} \in \R^{ + }$. In particular, 
${w}^{s}({z}_{3}^{\bar{t}}) = \max_{ \overline{{T}^{\bar{t}}}} {w}^{s} > 0$, 
\begin{equation*}
\bar{u}({z}_{3}^{\bar{t}}) + {s}{v}({z}_{3}^{\bar{t}}) = {w}^{s}({z}_{3}^{\bar{t}}) > {w}^{s}({z}_{1}^{\bar{t}}) = \bar{u}({z}_{1}^{\bar{t}}) + {s}{v}({z}_{1}^{\bar{t}}) \text{ for all } {s} > 0. 
\end{equation*}
This leads to a contradiction since ${v}({z}_{3}^{\bar{t}}) = 0$ and ${v}({z}_{1}^{\bar{t}}) > 0$. Therefore, the uniqueness of the eigenfunction is established. 

\textit{Step 2.3}. 
The global extrema of $\bar{u}$ occur exclusively at the endpoints of the longest side. 
In fact, from \eqref{CGY0726} and \autoref{lma72}, we deduce that $\bar{u}$ attains its global maximum at ${z}_{3}^{\bar{t}}$ and its global minimum at ${z}_{2}^{\bar{t}}$. When ${T}^{\bar{t}}$ is not subequilateral, the strict inequality \eqref{CGY0721b} implies that ${z}_{1}^{\bar{t}}$ is not a global extremum. When ${T}^{\bar{t}}$ is subequilateral, the fact that $\bar{u}({z}_{3}^{\bar{t}}) > 0 > \bar{u}({z}_{2}^{\bar{t}})$ ensures that $\bar{u}$ cannot be antisymmetric. Due to the symmetry of the triangle and the simplicity of the eigenvalue, $\bar{u}$ must be symmetric, and thus ${z}_{1}^{\bar{t}}$ is a global extremum. This completes Step 2. 

Combining the above two steps, we conclude that the set $\mathcal{S} \setminus \{1\}$ is a relatively open, closed, and non-empty subset. Therefore, $\mathcal{S} \supset [0, 1)$. Following Step 2, we know that $\mathcal{S}$ contains $1$, and hence $\mathcal{S} = [0, 1]$. This completes the proof of the theorem. 
\end{proof}

The symmetry and antisymmetry properties of the second Neumann eigenfunction on isosceles triangles were first proved in \cite{LS10}, while we point out that they are a direct consequence of the simplicity of the eigenvalue. 

\begin{rmk} 
Let ${T}$ be an isosceles triangle that is not equilateral, and let ${u}$ be a second Neumann eigenfunction in ${T}$. Then ${u}$ is symmetric if and only if the apex angle is less than $\pi/3$, while ${u}$ is antisymmetric if and only if the apex angle is greater than $\pi/3$. 
\end{rmk}

\begin{proof}
Consider ${T}_{\alpha}$, an isosceles triangle where the equal sides have length $1$ and the apex angle between them is $\alpha$. Since ${T}_{\alpha}$ is symmetric, the space of Neumann eigenfunctions on ${T}_{\alpha}$ decomposes into symmetric and antisymmetric functions with respect to the axis of symmetry of the triangle. Define $\mu^{s}(\alpha)$ (resp., $\mu^{a}(\alpha)$) as the smallest positive Neumann eigenvalue of ${T}_{\alpha}$ among the symmetric (resp., antisymmetric) functions with respect to the symmetry line of the triangle. Then $\mu_{2}({T}_{\alpha}) = \min\{\mu^{s}(\alpha), \mu^{a}(\alpha)\}$. From \autoref{thm74} and \autoref{prop28obt}, it follows that $\mu_{2}({T}_{\alpha})$ is simple. Therefore, 
\begin{equation} \label{CGY0725a}
\mu^{s}(\alpha) \neq \mu^{a}(\alpha) \text{ for } \alpha \neq \pi/3. 
\end{equation}
It is straightforward to verify the following estimates (e.g., \autoref{lma45} and \autoref{lma46LP}): 
\begin{equation*}
\mu^{s}(\alpha) < (2{j}_{0})^{2} \text{ and } \mu^{a}(\alpha) > \big(\tfrac{(\pi - \alpha) {j}_{0}}{\sin\alpha}\big)^{2}. 
\end{equation*}
Consequently, $\mu^{s}(\alpha) < \mu^{a}(\alpha)$ for sufficiently small $\alpha$. Combining this with $\mu^{s}(\pi/2) = 2\pi^{2}$, $\mu^{a}(\pi/2) = \pi^{2}$ and \eqref{CGY0725a}, we conclude
\begin{equation}
\mu^{s}(\alpha) < \mu^{a}(\alpha) \text{ for } 0 < \alpha < \pi/3 \text{ and } \mu^{s}(\alpha) > \mu^{a}(\alpha) \text{ for } \pi/3 < \alpha < \pi. 
\end{equation}
This establishes the symmetry and antisymmetry properties of the second Neumann eigenfunction on isosceles triangles. The proof is complete. 
\end{proof}

Combining this with \autoref{thm74}, \autoref{thm42}, and \autoref{lma71}, we determine the location of the critical points as follows. 

\begin{rmk}
Let ${T}$ be an acute triangle that is neither superequilateral nor equilateral. Then the non-vertex critical point exists and lies on the line segment joining the midpoint of the shortest side to the vertex with the largest interior angle. 
\end{rmk}


\section*{Acknowledgments}

The authors thank the anonymous reviewers for their careful reading  of the manuscript and valuable comments for improving the paper.
C. Gui is supported by  NSFC Key Program (Grant No. 12531010),  University of Macau research grants CPG2024-00016-FST, CPG2025-00032-FST, SRG2023-00011-FST, MYRG-GRG2023-00139-FST-UMDF, UMDF Professorial Fellowship of Mathematics, Macao SAR FDCT 0003/2023/RIA1 and  Macao SAR FDCT 0024/2023/RIB1.  R. Yao is supported by Guangdong Basic and Applied Basic Research Foundation (Grant No. 2025A1515011856) and the National Natural Science Foundation of China (Grant No. 12001543). 


\bibliographystyle{plain}
\bibliography{BibCGY}


\end{document}